\documentclass[letterpage, 11pt, notitlepage]{article}
\usepackage[margin=1.0in]{geometry}

\usepackage{amsmath, amssymb, natbib, graphicx, url}
\usepackage[utf8]{inputenc} % allow utf-8 input
\usepackage[T1]{fontenc}    % use 8-bit T1 fonts
\usepackage{hyperref}       % hyperlinks
\usepackage{booktabs}       % professional-quality tables
\usepackage{amsfonts}       % blackboard math symbols
\usepackage{nicefrac}       % compact symbols for 1/2, etc.
\usepackage{microtype}      % microtypography
\usepackage{bbold}
\usepackage{capt-of}

\usepackage[center]{caption}

\usepackage{cite,epstopdf,color,soul,mathabx}
\usepackage{makecell}
\usepackage{color}
\usepackage{enumerate}
\usepackage[shortlabels]{enumitem}
\usepackage{multicol}
\usepackage{empheq}
\usepackage{caption} 
\usepackage{algorithm,algorithmic}

\let\underbrace\LaTeXunderbrace

\allowdisplaybreaks

\newcommand{\Eset}{\mathbb{E}}

\newcommand{\Acal}{{\cal A}}

\newcommand{\Dcal}{{\cal D}}

\newcommand{\Hcal}{{\cal H}}

\newcommand{\Mcal}{{\cal M}}

\newcommand{\Ocal}{{\cal O}}
\newcommand{\Pcal}{{\cal P}}

\newcommand{\Scal}{{\cal S}}

\newcommand{\Xcal}{{\cal X}}

\newcommand{\argmin}{\mathop{\rm argmin}}

\newtheorem{lem}{Lemma}
\newtheorem{thm}{Theorem}

\newtheorem{assump}{Assumption}
\newtheorem{remark}{Remark}
\newcommand{\vct}[1]{\boldsymbol{#1}}
\newcommand{\va}{\vct{a}}
\newcommand{\vb}{\vct{b}}

\usepackage{tikz}

\newenvironment{keywords}
{\bgroup\leftskip 20pt\rightskip 20pt \small\noindent{\bfseries
Keywords:} \ignorespaces}%
{\par\egroup\vskip 0.25ex}
\newcommand{\qed}{\hfill $\blacksquare$}

\title{Fast Two-Time-Scale Stochastic Gradient Method with Applications in Reinforcement Learning\footnote{This is the most up-to-date version of the paper. An earlier version makes an incorrect claim on the complexity of the proposed method for convex functions, which is now removed. We note that this revision does not affect the results for strongly convex, Polyak-Lojasiewicz, and smooth non-convex functions. The revision also does not affect the applications of the proposed method discussed in Section~\ref{sec:application}, as the objective functions fall under strong convexity, Polyak-Lojasiewicz condition, and general non-convexity only.}}

\usepackage{times}

\author{
Sihan Zeng\thanks{J.P. Morgan AI Research}
\and Thinh T. Doan \thanks{UT Austin, Department of Aerospace Engineering \& Engineering Mechanics}
}

\begin{document}

\maketitle
\begin{abstract}

Two-time-scale optimization is a framework introduced in \citet{zeng2021two} that abstracts a range of policy evaluation and policy optimization problems in reinforcement learning (RL). 
Akin to bi-level optimization under a particular type of stochastic oracle, the two-time-scale optimization framework has an upper level objective whose gradient evaluation depends on the solution of a lower level problem, 
% which is the root of a strongly monotone operator.
% while the goal of the lower level problem is 
which is to find the root of a strongly monotone operator. 
% Akin to bi-level optimization under a particular type of stochastic oracle, the two-time-scale optimization framework has an upper level objective whose gradient evaluation depends on an auxiliary variable which needs to be obtained by finding the root of a strongly monotone operator (the lower level problem).
In this work, we propose a new method for solving two-time-scale optimization that achieves significantly faster convergence than the prior arts.  
%an accelerated algorithm for two-time-scale optimization, 
The key idea of our approach is to leverage an averaging step to improve the estimates of the operators in both lower and upper levels before using them to update the decision variables. These additional averaging steps eliminate the direct coupling between the main variables,
% and the noise due to the samples of the operators, 
enabling the accelerated performance of our algorithm. 
We characterize the finite-time convergence rates of the proposed algorithm under various conditions of the underlying objective function, including strong convexity, Polyak-Lojasiewicz condition, and general non-convexity. 
These rates significantly improve over the best-known complexity of the standard two-time-scale stochastic approximation algorithm. 
% Our results show that the proposed method significantly outperforms the existing two-time-scale stochastic gradient algorithms. 
% When applied to reinforcement learning, the proposed method presents a novel online sample-based approach that can converge faster than the existing algorithms in some settings. Finally, we support our theoretical results with numerical simulations to illustrate the performance of the proposed methods in solving reinforcement learning problems.  
When applied to RL, we show how the proposed algorithm specializes to novel online sample-based methods that surpass or match the performance of the existing state of the art. Finally, we support our theoretical results with numerical simulations in RL.\looseness=-1

\end{abstract}


\section{Introduction}
Inspired by the two-time-scale optimization framework for reinforcement learning (RL) introduced in \citet{zeng2021two}, we study the following optimization problem
\begin{align}
\theta^{\star}=\textstyle\argmin_{\theta\in\mathbb{R}^d}~h(\theta),
\label{eq:obj}
\end{align}
where $h:\mathbb{R}^d\rightarrow\mathbb{R}$ is a continuously differentiable 
(possibly nonconvex) function. 
% We will be interested in solving problem \eqref{eq:obj} using stochastic approximation methods where 
The gradient of $h$ can only be accessed through a special stochastic oracle 
% $F(\theta,\omega,X)\in\mathbb{R}^d$.
$F:\mathbb{R}^d\times\mathbb{R}^r\times\Xcal\rightarrow\mathbb{R}$.
The three arguments to $F$ are the decision variable $\theta\in\mathbb{R}^d$, an auxiliary variable $\omega\in\mathbb{R}^r$, and a random variable $X$ over the sample space $\Xcal$. Given any $\theta$, $F(\theta,\omega,X)$ returns an unbiased estimate of the gradient $\nabla h(\theta)$ only under a single ``correct'' setting of $\omega$ when the random variable $X$ follows some distribution $\xi$
\begin{align}
\nabla_{\theta} h(\theta) &= \Eset_{X\sim\xi}[F(\theta,\omega^{\star}(\theta),X)],\label{operator:F}
\end{align}
where $\omega^{\star}(\theta)\in\mathbb{R}^r$ is defined through another stochastic sampling operator $G:\mathbb{R}^d\times\mathbb{R}^r\times\Xcal\rightarrow\mathbb{R}$, as the solution to the equation
\begin{align}
    \mathbb{E}_{X\sim\xi}[G(\theta,\omega^{\star}(\theta),X]=0.\label{operator:G}
\end{align}
In this work, the operator $G$ is assumed to be strongly monotone with respect to $\omega$, implying that $\omega^{\star}(\theta)$ exists and is unique for any $\theta$.
Given any $\omega\neq\omega^{\star}(\theta)$, the stochastic gradient  $F(\theta,\omega,X)$ may be arbitrarily biased. 
By the first-order optimality condition, it is obvious that solving \eqref{eq:obj} is equivalent to finding a pair of solutions $(\theta^{\star},\omega^{\star})$ that solve the following two coupled nonlinear equations\looseness=-1
\begin{align}\label{prob:FG=0}
\left\{\begin{array}{ll}
\Eset_{X\sim\xi}[F(\theta^{\star},\omega^{\star},X)] = 0,\\
\Eset_{X\sim\xi}[G(\theta^{\star},\omega^{\star},X)] = 0.
\end{array}\right.    
\end{align}
As detailed in Section \ref{sec:application}, a wide range of RL methods, such as gradient-based temporal difference learning \citep{sutton2008convergent,sutton2009fast} and actor-critic algorithms for policy optimization \citep{yang2019provably,wu2020finite,chen2024finite}, aim to solve problems that are special cases of \eqref{prob:FG=0}. Thus, by developing and analyzing algorithms for the unified framework \eqref{prob:FG=0}, we can apply them to various RL problems and immediately deduce performance guarantees, without having to tailor analysis for each specific problem. \citet{hong2023two,zeng2021two} pioneer the effort to construct such a unified framework and study solving \eqref{prob:FG=0} with the classic two-time-scale stochastic approximation (SA) algorithm, in which the iterates $\theta_k,\omega_k$ are introduced to track $\theta^{\star},\omega^{\star}$ and updated in the following simple manner
\begin{align}
\begin{aligned}
    \theta_{k+1} = \theta_{k} - \alpha_k F(\theta_k,\omega_k,X_{k}),\quad\omega_{k+1} = \omega_k - \beta_k G(\theta_{k},\omega_k, X_k),
\end{aligned}
\label{eq:standard_SA}
\end{align}
where $X_k$ are independently drawn from $\xi$ in every iteration. 
The step sizes $\alpha_k$ and $\beta_k$ need to be chosen carefully according to the structure of $h$, usually one much smaller than the other, hence leading to the name ``two-time-scale''.
Despite the simplicity and popularity of this algorithm, its convergence rate is (at least as known so far) sub-optimal when compared to standard SGD under the most common oracle where we have access to $H(\theta,X)\in\mathbb{R}^d$ such that
\begin{align}
    \mathbb{E}_{X\in\xi}[H(\theta,X)]=\nabla h(\theta).\label{eq:standard_SGD_oracle}
\end{align}
This is due to the complex coupling between $\theta_k,$ $\omega_k,$ and the noise $X_k$, Taking strongly convex function $h$ as an example, standard SGD with $k$ queries of the stochastic gradient $H$ is able to reduce the optimality error to $\Ocal(1/k)$, while the two-time-scale SA algorithm can only achieve $\Ocal(1/k^{2/3})$.
\begin{remark}
Under a Lipschitz/Holder continuous gradient assumption on $\omega^{\star}$, \citet{shen2022single,han2024finite} show that the standard two-time-scale algorithm achieves $\Ocal(1/k)$ for strongly convex functions. We do not make this assumption in our work.
\end{remark}

% provide a unified framework to study different RL methods. In \citet{zeng2021two}, the authors studied a two-time-scale stochastic gradient method to solve problem \eqref{prob:FG=0}, where the decision variable $\theta_{k}$ is updated by the  standard gradient approach while implementing a stochastic approximation step to update the auxiliary variable $\omega_{k}$. While this approach can return a solution of \eqref{prob:FG=0}, its convergence rate is suboptimal when compared to its single-time-scale counterpart. 

% \paragraph{Main Contributions.} 
\subsection{Main Contributions} 

First, our work proposes a new variant of the two-time-scale SA algorithm, presented in Algorithm \ref{alg:main}, which for the first time achieves the optimal convergence rates when solving \eqref{prob:FG=0} under function $h$ satisfying various common structural properties, namely, strong convexity, Polyak-Lojasiewicz (PL) condition, and general non-convexity. 
% This comprehensively cover the different types of objective function $h$ considered in optimization and RL settings.
We summarize the rates as well as the convergence metrics and choices of step sizes in Table \ref{table:result_summary}.
Our key idea behind the innovation is to leverage an averaging technique to estimate the unknown operators $F,G$ through their samples before updating the decision variables. In other words, we ``denoise'' the estimate of the stochastic operators $F,G$, which leads to more stable updates of $\theta_k$ and $\omega_k$.
This seemingly simple modification effectively decouples the direct interaction between $\theta_k$ and $\omega_k$ and allows us to choose more aggressive step sizes, resulting in the accelerated convergence rates as well as simplified and elegant proofs. \looseness=-1

Second, we show how the results derived for the general framework apply to policy evaluation and optimization problem in reinforcement learning. We discuss three specific applications, the first of which is temporal difference learning with gradient correction (TDC) under linear and non-linear function approximation. In the linear setting, the objective is strongly convex and Algorithm~\ref{alg:main} specializes to a variant of the TDC algorithm guaranteed to converge with rate $\widetilde{\Ocal}(1/k)$. In the non-linear setting, the objective is non-convex and our convergence rate is $\widetilde{\Ocal}(1/k^{1/2})$. These results match the best-known rates of TDC in both cases.
% 
% , matching the best-known complexity of TDC. In the non-linear setting, the objective is non-convex and our algorithm converges with rate $\widetilde{\Ocal}(1/k^{1/2})$ to a stationary point, improving over the rate of TDC which is $\widetilde{\Ocal}(1/k^{2/5})$. 
The second application is sample-based policy optimization in the linear-quadratic regulator, which has an objective function satisfying the PL condition. The state-of-the-art online actor-critic algorithm for this problem has a complexity of $\widetilde{\Ocal}(1/k^{2/3})$ while our algorithm achieves the much faster rate $\widetilde{\Ocal}(1/k)$. The final application is sample-based policy optimization in an entropy-regularized MDP, the objective of which also observes the PL condition. The results in \citet{zeng2021two} imply that the classic two-time-scale SA algorithm finds the globally optimal policy with rate $\widetilde{\Ocal}(1/k^{2/3})$. Our algorithm in this context enjoys a finite-time rate of $\widetilde{\Ocal}(1/k)$, again significantly elevating the state of the art.

% We will show that the proposed methods will find the solutions of \eqref{prob:FG=0} with optimal convergence rates, where we will study a comprehensive coverage of different types of objective function $h$ considered in optimization and RL settings. Our results are summarized in Table \ref{table:result_summary}. While our convergence rates are at least as good as the known results in some settings (e.g., actor-critic), we provide some new results for others (e.g., gradient temporal difference learning with nonlinear function approximation). We will also provide numerical simulations to illustrate the proposed method outperforms the existing RL algorithms in some settings.  

% \subsection{Main Contributions}

\begin{table}[!h]
\centering
\caption{\small\sl Summary of finite-time complexity under various structural assumptions on $h$.
}
\setlength{\tabcolsep}{5pt}
\begin{tabular}{ccccc}
        \toprule
        \makecell{Structural \\ property}  &  Metric  &  \makecell{Complexity \\in this paper} & \makecell{Order of \\ step sizes} & \makecell{Best known complexity\\ of standard \\ two-time-scale SA} \\
        \midrule
        Strong Convexity & $\mathbb{E}[\|\theta_k-\theta^{\star}\|^2]$ & $\Ocal(k^{-1})$ & $k^{-1}$ & $\Ocal(k^{-\frac{2}{3}})$  \\
        % Convexity & $\mathbb{E}[h(\widebar{\theta}_k)-\theta^{\star}]$ & $\Ocal(k^{-\frac{1}{2}})$ & $k^{-\frac{1}{2}}$ & $\Ocal(k^{-\frac{1}{4}})$ \\
        PL Condition & $\mathbb{E}[h(\theta_k)-h(\theta^{\star})]$ & $\Ocal(k^{-1})$ & $k^{-1}$ & $\Ocal(k^{-\frac{2}{3}})$ \\
        Non-Convexity & $\min_{t \leq k}\mathbb{E}[\|\nabla h(\theta_t)\|^2]$ & $\Ocal(k^{-\frac{1}{2}})$ & $k^{-\frac{1}{2}}$ & $\Ocal(k^{-\frac{2}{5}})$ \\
        \bottomrule
        \bottomrule
        \end{tabular}
\label{table:result_summary}
\end{table}

\subsection{Related Work}
This paper is closely related to the existing works on two-time-scale stochastic approximation, stochastic bi-level and composite optimization, as well as those that analyze abstracted or specific RL algorithms. We discuss the latest advances in these subjects to give context to our contributions.

\noindent\textbf{Two-time-scale stochastic approximation.}
The aim of two-time-scale SA, pioneered by works including \citet{borkar1997stochastic,konda1999actor,mokkadem2006convergence}, is to solve a system of equations in the form of \eqref{prob:FG=0}. Existing works in this domain primarily focus on the settings where $F$ and $G$ are linear operators \citep{KondaT2004, DalalTSM2018, doan2019linear, Dalal_Szorenyi_Thoppe_2020,gupta2019finite} or strongly monotone operators \citep{MokkademP2006,doan2022nonlinear,shen2022single}. Our work can be seen as a generalization of this line of work to the case where $F$ is not a strongly monotone operator but the gradient map of a function $h$ exhibiting various structure.

\noindent\textbf{Stochastic bi-level and composite optimization.}
Bi-level optimization \citep{colson2007overview,chen2022single,chen2023optimal} studies solve an optimization program of the form
\begin{align}
    \textstyle\min_{x} f_1(x,y^{\star}(x)) \quad\text{subject to  }y^{\star}(x)\in\argmin_y f_2(x,y),\label{eq:bilevel_opt}
\end{align}
or equivalently (by the first-order condition), finding a pair of stationary points $(x^{\star},y^{\star})$ that observes
\begin{align*}
    \nabla_x f_1(x^{\star},y^{\star})=0,\quad \nabla_y f_2(x^{\star},y^{\star})=0.
\end{align*}
This is a special case of our objective \eqref{prob:FG=0} with $G$ being a gradient map. In RL problems $G$ usually abstracts the Bellman backup operator associated with policy evaluation, which is well-known to not be the gradient of any function. In this sense, our framework is more general and suitable for modeling applications in RL. 
% An additional difference between bi-level optimization and our work lies in the stochastic oracle. Our work and those on two-time-scale SA assume the access to $\nabla h(\theta)$, which is equivalent to $\nabla_x f_1(x,y^{\star}(x))$,

As another related subject, stochastic composite optimization \citep{ghadimi2012optimal,wang2017stochastic,chen2021solving} studies problems structured as
\begin{align*}
    \min_{x}g_1(g_2(x)).
\end{align*}
At a first glance, this optimization objective reduces to \eqref{eq:bilevel_opt} by choosing $g_1=f_1$, $g_2(x)=(x,y^{\star}(x))$, and $y^{\star}(x)$ to be the minimizer of $f_2(x,\cdot)$. 
Nevertheless, there is a key difference in the stochastic oracle. The assumption in stochastic composite optimization is that $g_2$ is continuous differentiable and that we can access an oracle that returns a direct stochastic estimate of $\nabla g_2$. This assumption is restrictive in most RL applications, in which there is no guarantee on the differentiability of $g_2$ and only indirect information of $g_2$ is accessible.

\noindent\textbf{Finite-Time Analysis of RL algorithms.}
Last but not least, we note the connection of our paper to the vast volume of recent literature that present finite-time and finite-sample analysis for various RL algorithms including temporal difference learning \citep{wang2021non,haque2023tight}, actor-critic algorithms for standard MDP \citep{wu2020finite,olshevsky2023small,chen2024finite}, and actor-critic algorithms for the linear-quadratic regulator \citep{yang2019provably,ju2022model,zhou2023single}. As samples in RL are sequentially drawn from a Markov chain, many of these works drop the i.i.d. noise assumption and realistically assume that the samples $X_k$ are Markovian. Markovian samples can either be ``time-invariant'' in the sense that the stationary distribution of the Markov chain does not change over time (for example in policy evaluation or off-policy learning) or ``controlled'' if the samples are obtained from an Markov chain whose stationary distribution changes as the policy is updated (for example in actor-critic algorithms). To focus on the main contribution which is the acceleration due to stable estimates of operators $F$ and $G$, our work assumes i.i.d. noise for simplicity. We note that the techniques to handle ``time-invariant'' and ``controlled'' Markovian noise have been well-known since \citet{srikant2019finite} and \citet{zou2019finite}.
% and frequently practiced in a range of recent works \citep{xiong2021non,chen2022finite,zeng2022finite,haque2023tight}.
% 
Extending the analysis of Algorithm~\ref{alg:main} to the Markovian noise setting can be done in a straightforward way and only makes the bounds worse by a $\log(k)$ factor. We discuss this extension in more details in Remark~\ref{rem:Markovian}.

% \tdoan{You can modify the one in our SIOPT paper for this section. I can take a look later.}

% We note that our objective in this work is slightly more general than that in the common bi-level optimization framework
% \begin{align}
% \theta^{\star}=\argmin_{\theta\in\mathbb{R}^d}~h(\theta,\omega^{\star}(\theta))\text{ such that }\omega^{\star}(\theta)=\argmin_{\omega\in\mathbb{R}^r}~\zeta(\theta,\omega),
% \label{eq:bilevel}
% \end{align}
% which can be regarded as a special case of \eqref{prob:FG=0} with $G$ being the stochastic gradient of function $\zeta$. However, as noted in \citet{zeng2021two}, most bi-level flavored problems in reinforcement learning (RL) have a lower level objective that can only be expressed as finding the root of a strongly monotone operator rather than an optimization problem, and therefore cannot be cast as \eqref{eq:bilevel}. We are motivated to study the generalized bi-level optimization formulation \eqref{eq:obj}-\eqref{operator:G} for this reason, which is essential for rigorously modeling applications in RL.

\vspace{5pt}
\noindent\textbf{Outline of the paper.} The rest of the paper is organized as follows. In Section~\ref{sec:algorithm} we present the proposed fast two-time-scale optimization algorithm. In Section~\ref{sec:main} we study the finite-time convergence rates of the algorithm under various functional structure. Section~\ref{sec:application} discusses how the optimization framework \eqref{eq:obj}-\eqref{operator:G} applies to RL problems and how the proposed algorithm solves these problems faster than (or as fast as) the prior arts. Finally, numerical simulations verifying the superiority of the proposed algorithm are presented in Section~\ref{sec:simulations}.

% \section{Formulation}

\section{Faster Two-Time-Scale Stochastic Gradient Method}\label{sec:algorithm}

We formally present the fast two-time-scale stochastic gradient method in Algorithm~\ref{alg:main}, in which we maintain and iteratively update four variables: $\theta_k$ and $\omega_k$ are used to track $\theta^{\star}$ and $\omega^{\star}(\theta^{\star})$, while $f_k$ and $g_k$ estimate the stochastic operators $F$ and $G$ at the current iterates $(\theta_k,\omega_k)$.
An averaging step is carried out on the estimate of the stochastic operators.
Note that Algorithm~\ref{alg:main} reduces to the standard two-time-scale SA algorithm \eqref{eq:standard_SA} if we set $\lambda_k=1$ for all $k$. With a large $\lambda_k$ the stochastic operator $F(\theta_k,\omega_k,X_k)$ estimates the true gradient $\nabla h(\theta_k)$ with high variance. By carefully choosing $\lambda_k$ to decay sufficiently fast, the algorithm updates the main variables with much more stable gradient estimates, which is the key driver behind the accelerated convergence rates.

\begin{algorithm}
\caption{Fast Two-Time-Scale Stochastic Gradient Method}
\label{alg:main}
\begin{algorithmic}[1]
\STATE{\textbf{Initialization:} decision variables $\theta_0\in\mathbb{R}^d$ and $\omega_0$, operator estimation variables $f_0\in\mathbb{R}^d$ and $g_0\in\mathbb{R}^r$, step size sequences $\{\alpha_k,\beta_k,\lambda_k\}$}
% \STATE{Observe a initial sample $X_0$}
\FOR{$k=0,1,\cdots,K-1$}
\STATE{Independent sample draw $X_k\sim\xi$}
\STATE{Decision variable update:
        \begin{align}
        \begin{aligned}
        \theta_{k+1} = \theta_k - \alpha_k f_k,\quad\omega_{k+1} = \omega_k - \beta_k g_k.
        \end{aligned}
        \label{alg:update_decision}
        \end{align}
        }
\STATE{Gradient estimation variable update:
\begin{align}
\begin{aligned}
f_{k+1} = (1 -\lambda_k)f_{ k } + \lambda_k F(\theta_k,\omega_k,X_k),\quad g_{k+1}=\left(1-\lambda_k\right) g_k+\lambda_k G(\theta_k, \omega_k, X_k).
\end{aligned}
\label{eq:update_auxiliary}
\end{align}}
\ENDFOR
% \RETURN $T$
\end{algorithmic}
\end{algorithm}

\section{Finite-Time Complexity}\label{sec:main}

In this section we study the finite-time complexity of Algorithm~\ref{alg:main} under various structural properties of the objective function, namely, strong convexity, PL condition, and general non-convexity. Since the algorithm draws exactly one sample in each iteration, the time complexity translates to an equivalent sample complexity. We first state the main technical assumptions
% to derive our theoretical results in the next section. These  assumptions, which are the same as the ones in \citet{zeng2021two}, are fairly standard in the existing literature
which are all fairly standard in the existing literature.

\begin{assump}[Lipschitz Continuity]\label{assump:Lipschitz}
There exists a positive constant $L$ such that for any $\theta_1,\theta_2\in\mathbb{R}^d$, $\omega_1,\omega_2\in\mathbb{R}^r$, and $X\in\Xcal$
\begin{align*}
    \left\|F(\theta_1,\omega_1,X)-F(\theta_2,\omega_2,X)\right\| &\leq L\left(\|\theta_1-\theta_2\|+\|\omega_1-\omega_2\|\right),\\
    \left\|G(\theta_1,\omega_1,X)-G(\theta_2,\omega_2,X)\right\| & \leq L\left(\|\theta_1-\theta_2\|+\|\omega_1-\omega_2\|\right),\\
    \left\|\nabla h(\theta_1)-\nabla h(\theta_2)\right\| &\leq L\|\theta_1-\theta_2\|,\\
    \|\omega^{\star}(\theta_1)-\omega^{\star}(\theta_2)\| &\leq L\|\theta_1-\theta_2\|.
\end{align*}
\end{assump}

Assumption~\ref{assump:Lipschitz} states that the operators $F$, $G$, and $\omega^{\star}$ are Lipschitz continuous and the function $h$ has Lipschitz gradients. This assumption is very common among works on two-time-scale SA \citep{zeng2021two,doan2022nonlinear,hong2023two}, and holds true in the RL applications to be discussed in Section~\ref{sec:application}.
Without any loss of generality, we assume $L\geq 1$.
We slightly abuse the notation to write
\begin{gather*}
F(\theta,\omega)\triangleq\mathbb{E}_{X\sim\xi}[F(\theta,\omega,X)],\quad G(\theta,\omega)\triangleq\mathbb{E}_{X\sim\xi}[G(\theta,\omega,X)].
\end{gather*}
Assumption~\ref{assump:Lipschitz} obviously implies
\begin{align*}
    &\left\|F(\theta_1,\omega_1)-F(\theta_2,\omega_2)\right\| \leq L\left(\|\theta_1-\theta_2\|+\|\omega_1-\omega_2\|\right),\\
    &\left\|G(\theta_1,\omega_1)-G(\theta_2,\omega_2)\right\| \leq L\left(\|\theta_1-\theta_2\|+\|\omega_1-\omega_2\|\right).
\end{align*}

% \begin{assump}\label{assump:gradh_bounded}
% There exists a positive constant $B$ such that for all $\theta\in\mathbb{R}^d$ 
% \begin{align*}
%     &\|\nabla h(\theta)\|\leq B.
% \end{align*}
% \end{assump}

\begin{assump}[Bounded Variance]\label{assump:noise_bounded}
There exists a positive constant $B$ such that for all $\theta\in\mathbb{R}^d$, $\omega\in\mathbb{R}^r$, and $X\in\Xcal$
\begin{align*}
    &\mathbb{E}_{X\sim \xi}[\|F(\theta,\omega,X)-F(\theta,\omega)\|^2]\leq B,\quad\mathbb{E}_{X\sim \xi}[\|G(\theta,\omega,X)-G(\theta,\omega)\|^2]\leq B.
\end{align*}
\end{assump}

The second assumption states that the operators $F$ and $G$ have bounded variance. This assumption is again standard in the literature on stochastic gradient algorithms under i.i.d. noise \citep{rakhlin2012making,shamir2013stochastic,hazan2014beyond}. 
We use $\Hcal_k=\{X_0,X_1,\cdots,X_k\}$ to denote the randomness information observed until iteration $k$.
\begin{assump}[Strong Monotonicity]\label{assump:stronglymonotone_G}
There exists a constant $\mu_G>0$ such that
\begin{align}
    \left\langle \Eset_{X\sim\xi}[G(\theta,\omega,X)], \omega-\omega^{\star}(\theta)\right\rangle \geq\mu_G \|\omega-\omega^{\star}(\theta)\|^2,\quad \forall \theta \in \mathbb{R}^{d},\omega \in \mathbb{R}^{r}.\label{assump:stronglymonotone_G:eq1}
\end{align}
\end{assump}
Assumption~\ref{assump:stronglymonotone_G} requires the operator $G$ to be strongly monotone in expectation. If $G$ is the gradient of a function $\zeta$, i.e. there exists a function $\zeta:\mathbb{R}^d\times \mathbb{R}^r\rightarrow\mathbb{R}$ such that $\nabla_{\omega}\zeta(\theta,\omega)=G(\theta,\omega)$ for all $\theta,\omega$, then \eqref{assump:stronglymonotone_G:eq1} amounts to the strong convexity of $\zeta$. Again, this assumption can be verified to hold in all RL applications discussed in Section~\ref{sec:application}.

To simplify notations we introduce the residual variables, which will all be driven to 0
% \begin{align}
% & \Delta f_k=f_k-F\left(\theta_k, \omega_k\right),\quad\Delta g_k=g_k-G\left(\theta_k, \omega_k\right), \quad \text {and } \quad y_k=\|\omega_k-\omega^{\star}(\theta_k)\|^2.\label{eq:residuals}
% \end{align}
\begin{align}
\begin{gathered}
\Delta f_k=f_k-F\left(\theta_k, \omega_k\right),\,\,\Delta g_k=g_k-G\left(\theta_k, \omega_k\right), 
% \\
% \quad \text {and } 
\,\, y_k=\|\omega_k-\omega^{\star}(\theta_k)\|^2.
\end{gathered}\label{eq:residuals}
\end{align}

\subsection{Strong Convexity}
We first analyze the convergence of Algorithm~\ref{alg:main} when the upper level objective is strongly convex, under the following step sizes
\begin{align}
\lambda_{k} = \frac{c_{\lambda}}{k+\tau+1},\quad \alpha_{k} = \frac{c_\alpha}{k+\tau+1},\quad \beta_{k} = \frac{c_\beta}{k+\tau+1}.\label{eq:step_sizes:stronglyconvex}
\end{align}
Here $c_{\alpha}$, $c_\beta$, $c_{\alpha}$, and $\tau\geq1$ are properly selected constants depending only on the structure of $h$ and $G$. Their exact choices are presented in Appendix~\ref{sec:proof_stronglyconvex}.

\begin{assump}[Strong Convexity]\label{assump:stronglyconvex}
The upper level objective $h$ is strongly convex, i.e. there exists a constant $\mu_h>0$ such that
\begin{align}
    \left\langle \nabla h(\theta_1)-\nabla h(\theta_2), \theta_1-\theta_2\right\rangle \geq\mu_h\|\theta_1-\theta_2\|^2,\quad \forall \theta_1,\theta_2 \in \mathbb{R}^{d}.\label{eq:stronglyconvex}
\end{align}
% \begin{align}
%     \left\langle \nabla_{\theta}h(\theta_1,\omega^{\star}(\theta_1))-\nabla_{\theta}h(\theta_2,\omega^{\star}(\theta_2)), \theta_1-\theta_2\right\rangle \geq\mu_h\|\theta_1-\theta_2\|^2,\quad \forall \theta_1,\theta_2 \in \mathbb{R}^{d}.\label{eq:stronglyconvex}
% \end{align}
\end{assump}

\begin{thm}\label{thm:stronglyconvex}
We recall the definition of residual variables in \eqref{eq:residuals} and denote $z_k = \|\theta_k-\theta^{\star}\|^2$.
Under Assumptions~\ref{assump:Lipschitz}-\ref{assump:stronglyconvex} and the step sizes in \eqref{eq:step_sizes:stronglyconvex}, the iterates of Algorithm~\ref{alg:main} satisfy for all $k$
\begin{align}
\mathbb{E}[z_{k+1}] \leq \frac{B\tau^2}{k+\tau+1}+\frac{\tau^2}{(k+\tau+1)^2}\left(\|\Delta f_0\|^2+\|\Delta g_0\|^2+z_0+y_0\right).\label{thm:stronglyconvex:eq0}
\end{align}
\end{thm}
The theorem states that when $h$ is strongly convex the iterate of Algorithm~\ref{alg:main} converge to the globally optimal solution $\theta^{\star}$ in the sense of $\|\theta_k-\theta^{\star}\|^2$ with a rate of $\Ocal(1/k)$, 
which matches the complexity of the standard SGD under the oracle $H$ in \eqref{eq:standard_SGD_oracle} \citep{rakhlin2012making}.
% matching the worst-case lower bound \citep{agarwal2009information}. 
This rate is much better than the best-known complexity of the standard two-time-scale SA algorithm established in \citet{hong2023two,zeng2021two}, which is $\Ocal(1/k^{2/3})$. 
\eqref{eq:obj}-\eqref{operator:G} under a strongly convex $h$ abstracts the objective of gradient-based temporal difference learning in RL, which we present in Section~\ref{sec:TDC}.

\begin{remark}\label{rem:Markovian}
We note that Markovian/i.i.d. sampling results in a difference in the analysis of $f_k,g_k$ in Lemma~\ref{lem:Delta_f:stronglyconvex} and \ref{lem:Delta_g:stronglyconvex}. In the proof of Lemma \ref{lem:Delta_f:stronglyconvex} in Appendix~\ref{sec:proof:lem:Delta_f:stronglyconvex} for example, it affects how we handle the cross term
\begin{align*}
b_{k} = \langle\Delta f_k,F(\theta,\omega_k,X_k)-F(\theta,\omega_k)\rangle\Rightarrow \left\{ \begin{array}{ll}
\mathbb{E}[b_{k}] = 0      & \text{under i.i.d samples}  \\
\mathbb{E}[b_{k}] \neq 0  &\text{under Markovian samples}    
\end{array}\right.  
\end{align*}
When $X_k$ is from a Markov chain with $\xi$ as the stationary distribution, $\mathbb{E}[b_k]$ captures the gap between $\xi$ and the distribution of $X_k$. This gap decays exponentially, under the assumption that the Markov chain is geometrically ergodic, which is a standard assumption in the literature \citep{xiong2021non,chen2022finite,zeng2022finite}.
Using techniques first developed in \citet{srikant2019finite}, we can show that $\Eset[\|b_k\|]\leq \Ocal(\alpha_k\log(1/\alpha_k))$. The rest of the proof of Theorem~\ref{thm:stronglyconvex} can proceed without modification, and eventually we have $\mathbb{E}[z_{k+1}]\leq\Ocal(\log(k+1)/(k+1))$, only differing from \eqref{thm:stronglyconvex:eq0} by a logarithm factor. The same extension can be made similarly on the results to be presented later in this section for other types of functions.
\end{remark}

\subsection{Non-Convexity with Polyak-Lojasiewicz Condition}

We now shift to the non-convex regime. The optimal solution to \eqref{eq:obj} may not be unique without the strong convexity assumption. Let $\Theta^{\star}\subseteq\mathbb{R}^{d}$ denote the set of optimal solutions to \eqref{eq:obj} and $h^{\star}$ denote the optimal function value.
We study the complexity of Algorithm~\ref{alg:main} when the upper-level objective function of \eqref{eq:obj} satisfies the PL condition below.

\begin{assump}[Polyak-Lojasiewicz Condition]\label{assump:PL_condition}
There exists a constant $\mu_h>0$ such that
\[\frac{1}{2}\|\nabla h(\theta)\|^2 \geq \mu_h \left(h(\theta)-h^{\star}\right), \quad \forall\theta\in\mathbb{R}^d.\]
\end{assump}
% The PL condition is commonly considered in the optimization literature under which stochastic gradient ascent converges with rate $\Ocal(1/k)$ \citep{karimi2016linear,yue2023lower}.
% 

We study Algorithm~\ref{alg:main} under the following choice of step sizes
\begin{align*}
\lambda_{k} = \frac{c_{\lambda}}{k+\tau+1},\quad \alpha_{k} = \frac{c_\alpha}{k+\tau+1},\quad \beta_{k} = \frac{c_\beta}{k+\tau+1},
\end{align*}
where $c_{\alpha}$, $c_\beta$, $c_{\lambda}$, and $\tau\geq1$ again need to be properly selected according to the structure of $h$ and $G$. Detailed choice of the parameters can be found in Appendix~\ref{sec:proof_PL}.

\begin{thm}\label{thm:PL}
We recall the residual variables defined in \eqref{eq:residuals} and denote $x_k = h(\theta_k)-h^{\star}$. 
Under Assumptions~\ref{assump:Lipschitz}-\ref{assump:stronglymonotone_G}, Assumption~\ref{assump:PL_condition}, and the step sizes above, the iterates of Algorithm~\ref{alg:main} satisfy for all $k$\looseness=-1
\begin{align*}
\mathbb{E}[x_{k+1}]
&\leq \frac{B\tau^2}{k+\tau+1}+\frac{\tau^2}{(k+\tau+1)^2}\left(\|\Delta f_0\|^2+\|\Delta g_0\|^2+x_0+y_0\right).
\end{align*}
\end{thm}
The theorem states that under the PL condition $h(\theta_k)$ converges with rate $\Ocal(1/k)$ to the globally optimal function value $h^{\star}$. This again outperforms the standard two-time-scale SA algorithm which has a complexity of $\Ocal(1/k^{2/3})$. Functions observing the PL condition include the policy optimization objectives in the linear-quadratic regulator and entropy-regularized MDP, which are the subjects of discussion in Section~\ref{sec:LQR} and \ref{sec:RegularizedMDP}.

\subsection{General Nonconvexity}
For non-convex functions without additional structure, we consider step sizes that decay as $1/\sqrt{k}$
\begin{align}
\lambda_{k} = \frac{1}{4(k+1)^{1/2}},\quad \alpha_{k} = \frac{\alpha_0}{(k+1)^{1/2}},\quad \beta_{k} = \frac{\beta_0}{(k+1)^{1/2}},\label{eq:step_sizes:nonconvex}
\end{align}
with $\alpha_0,\beta_0$ chosen such that $\alpha_0\leq\beta_0$, $\beta_0\leq\frac{1}{72L}$, $\beta_0\leq\lambda_0$, $\alpha_0\leq\frac{1}{72L^2}$.

\begin{thm}\label{thm:nonconvex}
Under Assumptions~\ref{assump:Lipschitz}-\ref{assump:stronglymonotone_G} and the step sizes in \eqref{eq:step_sizes:nonconvex}, we have for all $k$
\begin{align*}
\min_{t \leq k} \mathbb{E}\left[\|\nabla h(\theta_t)\|^2\right] &\leq \frac{32}{\alpha_0(k+1)^{1/2}}\left(\|\Delta f_0\|^2\hspace{-1pt}+\hspace{-1pt}\|\Delta g_0\|^2\hspace{-1pt}+\hspace{-1pt}(h(\theta_0)\hspace{-1pt}-\hspace{-1pt}h^{\star})\hspace{-1pt}+\hspace{-1pt}y_0\right)+\frac{4\log(k+2)}{\alpha_0(k+1)^{1/2}}.
\end{align*}
\end{thm}

The convergence in the non-convex case is not global, but rather to a stationary point of $h$. The state-of-the-art convergence rate of the standard two-time-scale SA algorithm is $\Ocal(1/k^{2/5})$ \citep{hong2023two,zeng2021two}, which we improve to $\Ocal(1/k^{1/2})$. Applications of the optimization framework for non-convex functions include gradient-based policy evaluation under non-linear function approximation which we present in Section~\ref{sec:TDC} and sample-based policy optimization under the standard MDP \citep{wu2020finite,olshevsky2023small,chen2024finite}.

\section{Motivating Applications}\label{sec:application}

% \subsection{Temporal-Difference Learning with Gradient Correction}
\subsection{Gradient-Based Policy Evaluation in Discounted-Reward MDP}
\label{sec:TDC}

Temporal difference learning with gradient correction (TDC) is a popular algorithm for policy evaluation in the setting where samples are collected by a behavior policy different from the one to be evaluated. Originally designed to work with linear function approximation \citep{sutton2008convergent,sutton2009fast}, this method has later been extended to the non-linear function approximation setting in \citet{maei2009convergent}.\looseness=-1
% , where the objective of the algorithm is to solve an optimization problem under stochastic oracles that exactly matches \eqref{eq:obj}-\eqref{operator:G}.

Consider a Markov decision process (MDP) characterized by $\Mcal=(\Scal,\Acal,\Pcal,r,\gamma)$, where $\Scal$, $\Acal$ are the state and action space, $\Pcal:\Scal\times\Acal\rightarrow\Delta_{\Scal}$ is the transition kernel, $r:\Scal\times\Acal\rightarrow[0,1]$ is the reward function, and $\gamma\in(0,1)$ is the discount factor.
In the linear function approximation case, each state $s\in\Scal$ is associated with a feature vector $\phi(s)\in\mathbb{R}^d$. We denote by $\Phi$ the stacked feature matrix
\[\Phi=\left[\begin{array}{ccc}
- & \phi(s_1)^T & - \\
\cdots & \cdots & \cdots \\
- & \phi(s_{|\Scal|})^T & -
\end{array}\right]\in\mathbb{R}^{|\Scal|\times d}.\]

Policy evaluation studies finding the value function of a policy $\pi\in\Delta_{\Acal}^{\Scal}$.
Given a parameter $\theta\in\mathbb{R}^d$, the approximated value function is $\Phi\theta\in\mathbb{R}^{|\Scal|}$. Let $\Pi$ denote the orthogonal projection onto the subspace spanned by $\Phi$. 
The objective of TDC is to find a parameter $\theta$ that minimizes the mean square projected Bellman error
% \begin{align}
%     \min_{\theta} J(\theta)\triangleq\mathbb{E}_{\mu_{\pi_b}}[\left(\Phi\theta-\Pi T^{\pi}\Phi\theta\right)^2]\label{eq:TDC_upper}
% \end{align}
\begin{align}
    \min_{\theta} J(\theta)\triangleq\|\Phi\theta-\Pi T^{\pi}\Phi\theta\|_{\mu_{\pi}}^2,\label{eq:TDC_upper}
\end{align}
where $\mu_{\pi}$ denotes the stationary distribution of states under policy $\pi$ and $\|v\|_{\mu_{\pi}}^2=\sum_{s\in\Scal}\mu_{\pi}(s)v(s)^2$ for any $v\in\mathbb{R}^{|\Scal|}$. 

As noted in \citet{sutton2009fast}, unbiased stochastic gradients of $J$ cannot be directly obtained using samples. Instead, we need to introduce an auxiliary variable $\omega\in\mathbb{R}^{d}$ and solve
% \begin{align}
% \begin{gathered}
%     \frac{1}{2}\nabla_{\theta}J(\theta)=
%     % \mathbb{E}_{s\sim\mu_{\pi},a\sim\pi(\cdot\mid s),s'\sim\Pcal(\cdot\mid s,a)}
%     -\mathbb{E}_{\pi}
%     [r(s,a)+\gamma\phi(s')^{\top}\theta\phi(s)-\phi(s)^{\top}\theta\phi(s)]+\gamma\mathbb{E}_{\pi}[\phi(s')\phi(s)^{\top}]\omega=0,\\
%     \mathbb{E}_{\pi}[\phi(s)\phi(s)^{\top}]\omega=\mathbb{E}_{\pi}
%     [r(s,a)+\gamma\phi(s')^{\top}\theta\phi(s)-\phi(s)^{\top}\theta\phi(s)].
% \end{gathered}\label{eq:TDC}
% \end{align}
\begin{gather}
    \nabla_{\theta}J(\theta)=
    % \mathbb{E}_{s\sim\mu_{\pi},a\sim\pi(\cdot\mid s),s'\sim\Pcal(\cdot\mid s,a)}
    -2\mathbb{E}_{\pi}
    \left[\Big(r(s,a)+\gamma\phi(s')^{\top}\theta-\phi(s)^{\top}\theta\Big)\phi(s)\right]+2\gamma\mathbb{E}_{\pi}[\phi(s')\phi(s)^{\top}]\,\omega=0,\label{eq:TDC:eq1}\\
    \mathbb{E}_{\pi}[\phi(s)\phi(s)^{\top}]\,\,\omega=\mathbb{E}_{\pi}
    \left[\Big(r(s,a)+\gamma\phi(s')^{\top}\theta-\phi(s)^{\top}\theta\Big)\phi(s)\right].\label{eq:TDC:eq2}
\end{gather}

It is straightforward to see that the objective \eqref{eq:TDC:eq1}-\eqref{eq:TDC:eq2} and stochastic oracles of the TDC algorithm exactly matches \eqref{eq:obj}-\eqref{operator:G}, where the operator in the lower level problem \eqref{eq:TDC:eq2} is strongly monotone and the upper level objective is strongly convex. As a result, Algorithm~\ref{alg:main} specializes to a variant of the TDC algorithm guaranteed to converge with rate $\widetilde{\Ocal}(1/k)$, which matches the best known convergence rate of TDC under linear function approximation \citep{ma2020variance}. When non-linear function approximation is considered, TDC solves a system of equations that resemble \eqref{eq:TDC:eq1}-\eqref{eq:TDC:eq2}, where the lower level operator is still strongly monotone but the upper level objective is non-convex \citep{maei2009convergent}. The best known convergence rate for TDC under non-linear function approximation is $\widetilde{\Ocal}(1/k^{-1/2})$ to a stationary point of $J$ \citep{wang2021non}, which we again match by our result in Theorem~\ref{thm:nonconvex}.\looseness=-1

\subsection{Policy Optimization for Linear-Quadratic Regulators}\label{sec:LQR}

Linear-quadratic regulator (LQR) studies finding an optimal control sequence $\{u_k\}$ that stabilizes a linear system with the minimum cost
\begin{align}
\begin{aligned}
\{u_k^{\star}\}\,\,=\,\,\min_{\{u_{k}\}}\,&\lim_{T\rightarrow\infty}\frac{1}{T}\mathbb{E}\Big[\sum_{k=0}^{T}\left(x_{k}^{\top} Q x_{k}+u_{k}^{\top} R u_{k}\right)\mid x_{0}\Big] \\
\text{subject to}\quad  &x_{k+1}=A x_{k}+B u_{k}+ w_k,
\end{aligned}
\label{eq:obj_LQR}
\end{align}
where $x_k\in\mathbb{R}^{d_1}$, $u_k\in\mathbb{R}^{d_2}$ are the state and the control variables, $w_k\in\mathbb{R}^{d_1}$ is the system noise, $A\in\mathbb{R}^{d_1\times d_1}$ and $B\in\mathbb{R}^{d_1\times d_2}$ are the transition matrices, and $Q\in\mathbb{R}^{d_1\times d_1}, R\in\mathbb{R}^{d_2\times d_2}$ are positive-definite cost matrices. We assume that $\{w_k\}$ is independently and identically distributed according to $N(0,\Psi)$ for some covariance matrix $\Psi\in\mathbb{R}^{d_1\times d_1}$.

It is a classic result in the control literature that the optimal control sequence $\{u_k^{\star}\}$ is a time-invariant linear function of the state 
\begin{align}
    u_k^\star = -K^{\star} x_k,\label{eq:lqr_u*}
\end{align}
for a matrix $K^{\star}\in\mathbb{R}^{d_2 \times d_1}$ as a function of $A,B,Q,R$ \citep{bertsekas2012dynamic}. 
Exploiting this fact allows us to re-formulate the LQR objective as finding the optimal control matrix $K$, which results in an optimization problem of the following form 
\begin{align}
    \min_K \quad&  \operatorname{trace}(P_{K}\Psi))\notag\\
    \text{subject to }\quad&P_{K}=Q+K^{\top} R K+(A-B K)^{\top} P_{K}(A-B K).\label{eq:obj_K_LQR:constraint}
\end{align}
To encourage exploration, we consider the regularized objective with $\sigma\geq 0$ and $\Psi_{\sigma}=\Psi+\sigma^2 BB^{\top}$
\begin{align}
    \underset{K}{\min} \,\, J(K)\triangleq \operatorname{trace}(P_{K}\Psi_{\sigma})+\sigma^2\operatorname{trace}(R)\quad\text{subject to }\,\,\eqref{eq:obj_K_LQR:constraint},
\label{eq:obj_K_LQR}
\end{align}
which is known to have the same optimal solution as the unregularized problem \citep{gravell2020learning}.

We are interested in finding the solution to \eqref{eq:obj_K_LQR} under unknown transition matrices $A$ and $B$, by taking a gradient-based approach. The natural gradient of $J$ with respect to $K$, which we denote by $\widetilde{\nabla}_K J$ and which can be regarded as $\nabla_{K} J$ on a slightly transformed system of coordinates, has the following closed-form expression
\begin{align*}
    \widetilde{\nabla}_K J=2\left(R+B^{\top} P_{K} B\right) K-2B^{\top} P_{K} A.
\end{align*}

To obtain a reliable stochastic estimate of $\widetilde{\nabla}_K J$ the key is to estimate $R+B^{\top} P_K B$ and $B^{\top} P_K A$. For simplicity of notation, we define
\begin{align}
    \Omega_{K}=\left(\begin{array}{cc}
    \Omega_{K}^{11} & \Omega_{K}^{12} \\
    \Omega_{K}^{21} & \Omega_{K}^{22}
    \end{array}\right)=\left(\begin{array}{cc}
    Q+A^{\top} P_{K} A & A^{\top} P_{K} B \\
    B^{\top} P_{K} A & R+B^{\top} P_{K} B
    \end{array}\right),
\end{align}
of which $R+B^{\top} P_K B$ and $B^{\top} P_K A$ are sub-matrices.
We also define $\text{svec}(\cdot)$ as the operation to vectorize the upper triangular sub-matrix of a given symmetric matrix with the off-diagonal entries weighted by $\sqrt{2}$. We further define $\phi(x, u)=\text{svec}(\left[x^{\top},u^{\top}\right]^{\top}\left[x^{\top},u^{\top}\right])$ for any $x\in\mathbb{R}^{d_1},u\in\mathbb{R}^{d_2}$, which can be seen as the feature vector for the system in state $(x,u)$.
% \[\phi(x, u)=\text{svec}\Big(\left[\begin{array}{l}
% x \\
% u
% \end{array}\right]\left[\begin{array}{l}
% x \\
% u
% \end{array}\right]^{\top}\Big).\]
It is shown in \citet{yang2019provably} that $\Omega_{K}$ and $J(K)$ jointly satisfy 
\begin{align}
\hspace{-0.3cm}\mathbb{E}_{x\sim\mu_K,u\sim N(-K x,\sigma^2 I)}\left[M_{x,u,x',u'}\right]\left[\begin{array}{c}
J(K) \\
\text{svec}(\Omega_K)
\end{array}\right]=\mathbb{E}_{x\sim\mu_K,u\sim N(-K x,\sigma^2 I)}\left[c_{x,u}\right],
\label{eq:lqr_Goperator}
\end{align}
where the matrix $M_{x,u,x',u'}$ and vector $c_{x,u}$ are
\begin{gather*}
M_{x,u,x',u'}\hspace{-2pt}=\hspace{-2pt}\left[\hspace{-4pt}\begin{array}{cc}
1 & 0 \\
\phi(x, u) & \hspace{-4pt}\phi(x, u)\left[\phi(x, u)-\phi\left(x', u'\right)\right]^{\top}
\end{array}\hspace{-6pt}\right]\hspace{-2pt},\,\,c_{x,u}\hspace{-2pt}=\hspace{-2pt}\left[\hspace{-4pt}\begin{array}{c}
x^{\top} Q x+u^{\top} R u \\
\left(x^{\top} Q x+u^{\top} R u\right)\hspace{-2pt}\phi(x, u)
\end{array}\hspace{-7pt}\right]\hspace{-2pt}.
\end{gather*}

Equation~\eqref{eq:lqr_Goperator} maps to the lower level problem \eqref{operator:G}.
Under the assumption that $K$ is a stable controller with respect to $A$ and $B$ \citep{yang2019provably}, the linear operator $\mathbb{E}[M_{x,u,x',u'}]$ is strongly monotone \citep{yang2019provably}. The upper level problem that corresponds to \eqref{operator:F} is to solve 
% \[\widetilde{\nabla}_K J=2\left(R+B^{\top} P_{K} B\right) K-2B^{\top} P_{K} A,\]
$\widetilde{\nabla}_K J=0$,
which is known to observe the PL condition \citep{fazel2018global}. This means that we can apply Algorithm~\ref{alg:main} to this problem, which specializes to an single-looped ``actor-critic'' algorithm with a finite-time convergence rate of $\widetilde{\Ocal}(1/k)$ to the globally optimal solution measured by $J(K_k)-J(\theta^{\star})$. The complexity of the current state-of-the-art single-looped algorithm for LQR is $\widetilde{\Ocal}(1/k^{2/3})$, which we improve over. We also note the recent few works on nested-loop actor-critic algorithms that rely on batched multi-trajectory samples \citep{ju2022model,zhou2023single}. These algorithms achieve the complexity $\widetilde{\Ocal}(1/k)$ as well, but in practice are much less convenient to implement than our algorithm which is single-looped and only requires a single online trajectory of samples. 

\begin{algorithm}[t]
\caption{Fast Single-Loop Actor-Critic Algorithm for LQR}
\label{Alg:AC_LQR}
\begin{algorithmic}[1]
\STATE{\textbf{Initialization:} Control matrix $K_0$, auxiliary variable $\hat{\Omega}_0, \hat{J}_0$, gradient estimation variable $f_0,g_0^J,g_0^{\Omega}$, step size sequences $\{\alpha_k,\beta_k,\lambda_k\}$}
\STATE{Sample an initial state $x_0\sim\Dcal$, takes an initial control $u_0\sim N(-K_0 x_0,\sigma^2 I)$, and observe cost $c_0=x_0^{\top}Qx_0+u_0^{\top}Ru_0$ and next state $x_1\sim N(Ax_0+Bu_0,\Psi)$}
\FOR{$k=0,1,2,...$}
\STATE{Take the control $u_{k+1}\sim N(-K_{k} x_{k+1},\sigma^2 I)$. Observe cost $c_{k+1}=x_{k+1}^{\top}Qx_{k+1}+u_{k+1}^{\top}Ru_{k+1}$ and next state $x_{k+2}\sim N(Ax_{k+1}+Bu_{k+1},\Psi)$}
\STATE{Control matrix update:
        \begin{align*}
        K_{k+1} = K_k-\alpha_k f_k
        \end{align*}}
\STATE{Auxiliary variable update:
        \begin{align*}
        \hat{J}_{k+1} = \hat{J}_k - \beta_k g_k^J, \quad \hat{\Omega}_{k+1} =\hat{\Omega}_{k} - \beta_k g_k^{\Omega}
        \end{align*}
        }
\STATE{Gradient estimation variable update:
        \begin{align*}
        f_{k+1}&=(1-\lambda_k)f_{k}+\lambda_k(\hat{\Omega}_k^{22}K_k-\hat{\Omega}_k^{21})\\
        g_{k+1}^J &= (1-\lambda_k)g_{k}^J +\lambda_k(\hat{J}_k-x_k^{\top} Q x_k - u_k^{\top} R u_k) \\
        g_{k+1}^{\Omega} &=(1-\lambda_k) g_{k}^{\Omega} - \lambda_k\hspace{-2pt}\left[\begin{array}{c}
        x_k \\
        u_k
        \end{array}\right]\hspace{-3pt}\left[\begin{array}{c}
        x_k \\
        u_k
        \end{array}\right]^{\top}\hspace{-7pt}\left(\left[
        x_k^{\top}\,\,\, u_k^{\top}\right]\hat{\Omega}_k\hspace{-2pt}\left[\begin{array}{c}
        x_k \\
        u_k
        \end{array}\right]\hspace{-2pt}+\hspace{-2pt}\hat{J}_k-x_k^{\top} Q x_k - u_k^{\top} R u_k\right)
        \end{align*}
        }
\ENDFOR
% \RETURN $T$
\end{algorithmic}
\end{algorithm}

\subsection{Policy Optimization for Entropy-Regularized MDPs}\label{sec:RegularizedMDP}

The optimization framework \eqref{eq:obj}-\eqref{operator:G} also applies to the policy optimization problem in the entropy-regularized MDP. Consider the MDP $\Mcal=(\Scal,\Acal,\Pcal,r,\gamma)$ introduced in Section~\ref{sec:TDC}. The regularized value function of a policy $\pi\in\Delta_{\Acal}^{\Scal}$ is
\begin{align*}
    V_{\tau}^{\pi}(s)&=\mathbb{E}_{a_k\sim\pi(\cdot\mid s_k),s_{k+1}\sim\Pcal(\cdot\mid s_k,a_k)}\Big[\sum_{k=0}^{\infty}\gamma^k\big(r(s_k,a_k)-\tau\log\pi(a_k\mid s_k)\big)\mid s_0=s\Big],
\end{align*}
where $\tau$ is a positive regularization weight. Under initial state distribution $\rho\in\Delta_{\Scal}$, the expected cumulative reward under policy $\pi$ is $J_{\tau}(\pi)\triangleq\mathbb{E}_{s\sim\rho}[V_{\tau}^{\pi}(s)]$. We consider the tabular setting where the policy is represented through the softmax function by a parameter $\theta\in\mathbb{R}^{|\Scal|\times|\Acal|}$
\begin{align*}
    \pi_{\theta}(a \mid s)&=\frac{\exp \left(\theta_{s, a}\right)}{\sum_{a' \in \Acal} \exp \left(\theta_{s, a'}\right)}.
 \end{align*}

Policy optimization studies finding the parameter $\theta^{\star}$ that maximizes the cumulative reward $J_{\tau}(\pi_{\theta})$.
By the first-order condition, this is equivalent to solving
\begin{align*}
    &\nabla_{\theta} J_{\tau}(\pi_{\theta})=\frac{1}{1-\gamma}\mathbb{E}_{\pi_{\theta}}\big[\left(r(s,a)-\tau\log\pi_{\theta}(a\mid s)+\gamma V_{\tau}^{\pi_{\theta}}(s')-V_{\tau}^{\pi_{\theta}}(s)\right)\nabla_{\theta}\log\pi_{\theta}(a\mid s)\big]=0.
    % \label{eq:regularizedMDP:operatorF}
\end{align*}

Evaluating the gradient requires computing the regularized value function $V_{\tau}^{\pi_{\theta}}\in\mathbb{R}^{|\Scal|}$, which satisfies the regularized Bellman equation
\begin{align}
\mathbb{E}_{\pi_{\theta}}\Big[r(s,a)-\tau\log\pi_{\theta}(a\mid s)+\gamma V_{\tau}^{\pi_{\theta}}(s')-V_{\tau}^{\pi_{\theta}}(s)\Big]=0.
\label{eq:regularizedMDP:operatorG}
\end{align}

This objective again falls under the general two-time-scale optimization framework, with the operator in lower level problem \eqref{eq:regularizedMDP:operatorG} being strongly monotone. The upper level objective function $J_{\tau}$ satisfies the PL condition \citep{mei2020global}. Applying Algorithm~\ref{alg:main} to this problem results in a fast actor-critic algorithm (Algorithm~\ref{Alg:AC_LQR}) that is able to reduce the optimality gap $J_{\tau}(\theta_k)-J_{\tau}(\theta^{\star})$ with rate $\widetilde{\Ocal}(1/k)$. The analysis in \citet{zeng2021two} implies that the standard two-time-scale algorithm \eqref{eq:standard_SA} solves this problem with a complexity of $\widetilde{\Ocal}(1/k^{2/3})$, which we vastly improve over.

% \subsection{Online Actor-Critic Algorithm for Infinite-Horizon Average-Reward MDPs}

% \subsection{Stochastic Minimax Optimization}

% Minimax optimization considers problems of the following form
% \begin{align*}
% \min_{\phi\in\mathbb{R}^a}\max_{\psi\in\mathbb{R}^a} \mathbb{E}_{X\sim\xi}[\ell(\phi,\psi)],
% \end{align*}
% which is equivalent to
% \begin{align}
% \min_{\phi\in\mathbb{R}^a}\mathbb{E}_{X\sim\xi}[\ell(\phi,\psi^{\star}(\phi))]\text{ such that }\psi^{\star}(\phi)=\min_{\psi\in\mathbb{R}^a}-\mathbb{E}_{X\sim\xi}[\ell(\phi,\psi)].\label{eq:minimax_opt}
% \end{align}
% It is obvious that \eqref{eq:minimax_opt} is a special case of the bi-level optimization objective \eqref{eq:obj} with the lower-level objective being the negative of the upper-level one. 

\section{Numerical Simulations}\label{sec:simulations}
We experimentally verify that Algorithm~\ref{alg:main} converges faster than \eqref{eq:standard_SA} in solving 1) a policy evaluation problem under linear function approximation, and 2) policy optimization for the LQR.

\begin{figure}
\centering
\begin{minipage}{.4\textwidth}
  \centering
  \includegraphics[width=\linewidth]{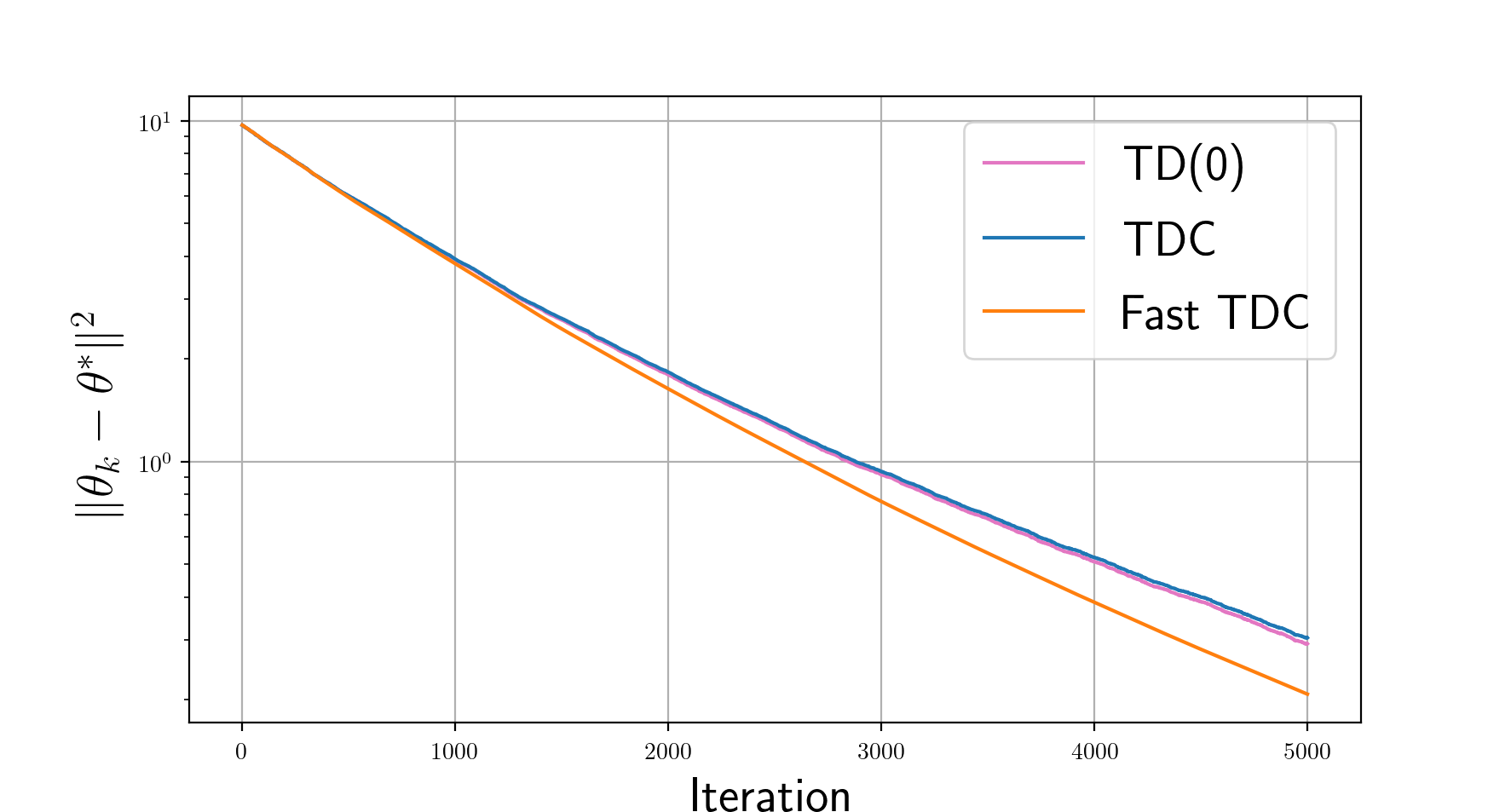}
  \caption{Fast TDC Algorithm for Random Policy Evaluation}
  \label{fig:random_policy_eval}
\end{minipage}%
\begin{minipage}{.02\textwidth}
\hspace{0pt}
\end{minipage}
\begin{minipage}{.4\textwidth}
  \centering
  \includegraphics[width=\linewidth]{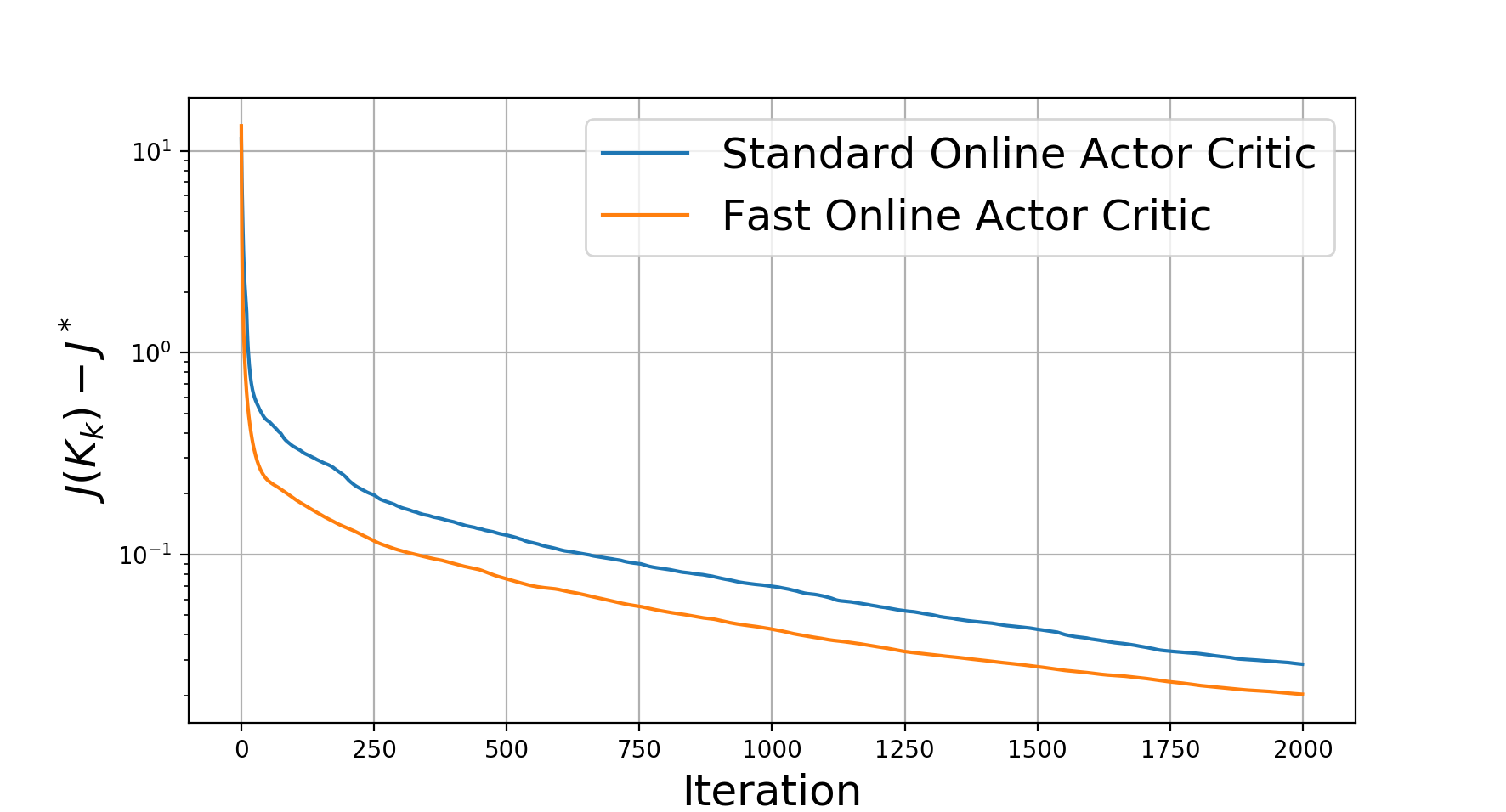}
  \caption{Performance of Fast Actor-Critic Algorithm for LQR}
  \label{fig:policy_optimization_LQR}
\end{minipage}
\end{figure}

\noindent\textbf{Policy evaluation problem under linear function approximation.}
As we have discussed in Section~\ref{sec:TDC}, TDC is the state-of-the-art policy evaluation algorithm under off-policy samples which enjoys a convergence rate of $\widetilde{\Ocal}(1/k)$ under linear function approximation. Numerically we evaluate Algorithm~\ref{alg:main} against TDC on the task of evaluating the value function of a randomly generated policy in a randomly generated environment with $|\Scal|=|\Acal|=50$ and the feature dimension $d=10$.
Figure~\ref{fig:random_policy_eval} shows the empirical performance gap between the standard TDC and Algorithm~\ref{alg:main}, though both algorithm have the same complexity in theory.
For comparison we also include TD(0) learning under linear function approximation. 
% It is known that TD(0) learning under off-policy samples may diverge in the worst case \citep{sutton2009fast}.
The simulation results in Figure~\ref{fig:random_policy_eval} numerically confirms the superiority of the proposed algorithm. Details of the experimental setup can be found in Appendix~\ref{sec:appendix_simulation}.

\noindent\textbf{Policy optimization for LQR.} 
We experiment with the small-dimensional environment considered in \citet{zeng2021two} where the transition and costs matrices are chosen to be
\begin{align}
A=\left[\begin{array}{ccc}
0.5 & 0.01 & 0 \\
0.01 & 0.5 & 0.01 \\
0 & 0.01 & 0.5
\end{array}\right], \quad B=\left[\begin{array}{ll}
1 & 0.1 \\
0 & 0.1 \\
0 & 0.1
\end{array}\right], \quad Q=I_{3},\quad R=I_{2}.
\end{align}

We apply Algorithm~\ref{Alg:AC_LQR} to this problem along with online actor-critic algorithm studied in \citet{zeng2021two}, which is equivalent to \eqref{eq:standard_SA} applied to this context. In Figure \ref{fig:policy_optimization_LQR}, we plot the error decay under the two algorithms, which again shows that the proposed algorithm enjoys faster convergence.

% \begin{figure}[ht]
%   \includegraphics[width=\linewidth]{Figures/random_policy_eval.png}
%   \vspace{-.7cm}
%   \caption{Performance of Fast TDC Algorithm for Random Policy Evaluation in Random MDPs under Linear Function Approximation. Left: MDP with cardinality of state and action space $|\Scal|=|\Acal|=50$ and feature dimension $d=10$. Right: Larger MDP with cardinality of state and action space $|\Scal|=|\Acal|=100$ and feature dimension $d=30$.}
%   \label{fig:random_policy_eval}
% \end{figure}

% \section{Conclusions \& Future Work}\label{sec:conclusion}

\section*{Disclaimer}
This paper was prepared for informational purposes in part by
the Artificial Intelligence Research group of JP Morgan Chase \& Co and its affiliates (``JP Morgan''),
and is not a product of the Research Department of JP Morgan.
JP Morgan makes no representation and warranty whatsoever and disclaims all liability,
for the completeness, accuracy or reliability of the information contained herein.
This document is not intended as investment research or investment advice, or a recommendation,
offer or solicitation for the purchase or sale of any security, financial instrument, financial product or service,
or to be used in any way for evaluating the merits of participating in any transaction,
and shall not constitute a solicitation under any jurisdiction or to any person,
if such solicitation under such jurisdiction or to such person would be unlawful.

\section*{Acknowledgement}
This work was supported by the National Science Foundation under awards ECCS-CAREER-2339509 and CCF-2343599.
\medskip

\bibliographystyle{plainnat} 
\bibliography{references}

@inproceedings{fazel2018global,
  title={Global convergence of policy gradient methods for the linear quadratic regulator},
  author={Fazel, Maryam and Ge, Rong and Kakade, Sham and Mesbahi, Mehran},
  booktitle={International conference on machine learning},
  pages={1467--1476},
  year={2018},
  organization={PMLR}
}

@article{hong2023two,
  title={A two-timescale stochastic algorithm framework for bilevel optimization: Complexity analysis and application to actor-critic},
  author={Hong, Mingyi and Wai, Hoi-To and Wang, Zhaoran and Yang, Zhuoran},
  journal={SIAM Journal on Optimization},
  volume={33},
  number={1},
  pages={147--180},
  year={2023},
  publisher={SIAM}
}

@article{yang2019provably,
  title={Provably global convergence of actor-critic: A case for linear quadratic regulator with ergodic cost},
  author={Yang, Zhuoran and Chen, Yongxin and Hong, Mingyi and Wang, Zhaoran},
  journal={Advances in neural information processing systems},
  volume={32},
  year={2019}
}

@inproceedings{karimi2016linear,
  title={Linear convergence of gradient and proximal-gradient methods under the {P}olyak-{\L}ojasiewicz condition},
  author={Karimi, Hamed and Nutini, Julie and Schmidt, Mark},
  booktitle={Machine Learning and Knowledge Discovery in Databases: European Conference, ECML PKDD 2016, Riva del Garda, Italy, September 19-23, 2016, Proceedings, Part I 16},
  pages={795--811},
  year={2016},
  organization={Springer}
}

@article{zeng2022finite,
  title={Finite-time convergence rates of decentralized stochastic approximation with applications in multi-agent and multi-task learning},
  author={Zeng, Sihan and Doan, Thinh T and Romberg, Justin},
  journal={IEEE Transactions on Automatic Control},
  year={2022},
  publisher={IEEE}
}

@article{zeng2021two,
  title={A two-time-scale stochastic optimization framework with applications in control and reinforcement learning},
  author={Zeng, Sihan and Doan, Thinh T and Romberg, Justin},
  journal={SIAM Journal on Optimization},
  volume={34},
  number={1},
  pages={946--976},
  year={2024},
  publisher={SIAM}
}

@inproceedings{rakhlin2012making,
  title={Making gradient descent optimal for strongly convex stochastic optimization},
  author={Rakhlin, Alexander and Shamir, Ohad and Sridharan, Karthik},
  booktitle={Proceedings of the 29th International Coference on International Conference on Machine Learning},
  pages={1571--1578},
  year={2012}
}

@article{colson2007overview,
  title={An overview of bilevel optimization},
  author={Colson, Beno{\^\i}t and Marcotte, Patrice and Savard, Gilles},
  journal={Annals of operations research},
  volume={153},
  number={1},
  pages={235--256},
  year={2007},
  publisher={Springer}
}

@article{chen2021solving,
  title={Solving stochastic compositional optimization is nearly as easy as solving stochastic optimization},
  author={Chen, Tianyi and Sun, Yuejiao and Yin, Wotao},
  journal={IEEE Transactions on Signal Processing},
  volume={69},
  pages={4937--4948},
  year={2021},
  publisher={IEEE}
}

@article{ghadimi2012optimal,
  title={Optimal stochastic approximation algorithms for strongly convex stochastic composite optimization i: A generic algorithmic framework},
  author={Ghadimi, Saeed and Lan, Guanghui},
  journal={SIAM Journal on Optimization},
  volume={22},
  number={4},
  pages={1469--1492},
  year={2012},
  publisher={SIAM}
}

@article{wang2017stochastic,
  title={Stochastic compositional gradient descent: algorithms for minimizing compositions of expected-value functions},
  author={Wang, Mengdi and Fang, Ethan X and Liu, Han},
  journal={Mathematical Programming},
  volume={161},
  number={1},
  pages={419--449},
  year={2017},
  publisher={Springer}
}

@inproceedings{chen2022single,
  title={A single-timescale method for stochastic bilevel optimization},
  author={Chen, Tianyi and Sun, Yuejiao and Xiao, Quan and Yin, Wotao},
  booktitle={International Conference on Artificial Intelligence and Statistics},
  pages={2466--2488},
  year={2022},
  organization={PMLR}
}

@article{shen2022single,
  title={A single-timescale analysis for stochastic approximation with multiple coupled sequences},
  author={Shen, Han and Chen, Tianyi},
  journal={Advances in Neural Information Processing Systems},
  volume={35},
  pages={17415--17429},
  year={2022}
}

@article{chen2023optimal,
  title={Optimal Algorithms for Stochastic Bilevel Optimization under Relaxed Smoothness Conditions},
  author={Chen, Xuxing and Xiao, Tesi and Balasubramanian, Krishnakumar},
  journal={arXiv preprint arXiv:2306.12067},
  year={2023}
}

@inproceedings{shamir2013stochastic,
  title={Stochastic gradient descent for non-smooth optimization: Convergence results and optimal averaging schemes},
  author={Shamir, Ohad and Zhang, Tong},
  booktitle={International conference on machine learning},
  pages={71--79},
  year={2013},
  organization={PMLR}
}

@article{wu2020finite,
  title={A finite-time analysis of two time-scale actor-critic methods},
  author={Wu, Yue Frank and Zhang, Weitong and Xu, Pan and Gu, Quanquan},
  journal={Advances in Neural Information Processing Systems},
  volume={33},
  pages={17617--17628},
  year={2020}
}

@article{gupta2019finite,
  title={Finite-time performance bounds and adaptive learning rate selection for two time-scale reinforcement learning},
  author={Gupta, Harsh and Srikant, Rayadurgam and Ying, Lei},
  journal={Advances in Neural Information Processing Systems},
  volume={32},
  year={2019}
}

@article{MokkademP2006,
 author = {A. Mokkadem and M. Pelletier},
 journal = {The Annals of Applied Probability},
 number = {3},
 pages = {1671-1702},
title = {Convergence Rate and Averaging of Nonlinear Two-Time-Scale Stochastic Approximation Algorithms},
 volume = {16},
 year = {2006},
}

@article{doan2022nonlinear,
  title={Nonlinear two-time-scale stochastic approximation convergence and finite-time performance},
  author={Doan, Thinh T},
  journal={IEEE Transactions on Automatic Control},
  year={2022},
  publisher={IEEE}
}

@inproceedings{doan2019linear,
  title={Linear two-time-scale stochastic approximation a finite-time analysis},
  author={Doan, Thinh T and Romberg, Justin},
  booktitle={2019 57th Annual Allerton Conference on Communication, Control, and Computing (Allerton)},
  pages={399--406},
  year={2019},
  organization={IEEE}
}

@article{Dalal_Szorenyi_Thoppe_2020, title={A Tale of Two-Timescale Reinforcement Learning with the Tightest Finite-Time Bound}, volume={34}, 
number={04}, 
journal={Proceedings of the AAAI Conference on Artificial Intelligence}, 
author={Dalal, Gal and Szorenyi, Balazs and Thoppe, Gugan}, year={2020}, month={Apr.}, pages={3701-3708}, 
}

@article{KondaT2004,
 author = {V. R. Konda and J. N. Tsitsiklis},
 journal = {The Annals of Applied Probability},
 number = {2},
 pages = {796--819},
 publisher = {Institute of Mathematical Statistics},
 title = {Convergence Rate of Linear Two-Time-Scale Stochastic Approximation},
 volume = {14},
 year = {2004},
}

@inproceedings{DalalTSM2018,
  title={Finite Sample Analysis of Two-Timescale Stochastic   Approximation with Applications to Reinforcement Learning},
  author={G. Dalal and G. Thoppe and B. Sz{\"o}r{\'e}nyi and S. Mannor},
  booktitle={COLT},
  year={2018}
}

@article{chen2024finite,
  title={Finite-time analysis of single-timescale actor-critic},
  author={Chen, Xuyang and Zhao, Lin},
  journal={Advances in Neural Information Processing Systems},
  volume={36},
  year={2024}
}

@article{hazan2014beyond,
  title={Beyond the regret minimization barrier: optimal algorithms for stochastic strongly-convex optimization},
  author={Hazan, Elad and Kale, Satyen},
  journal={The Journal of Machine Learning Research},
  volume={15},
  number={1},
  pages={2489--2512},
  year={2014},
  publisher={JMLR. org}
}

@inproceedings{sutton2009fast,
  title={Fast gradient-descent methods for temporal-difference learning with linear function approximation},
  author={Sutton, Richard S and Maei, Hamid Reza and Precup, Doina and Bhatnagar, Shalabh and Silver, David and Szepesv{\'a}ri, Csaba and Wiewiora, Eric},
  booktitle={Proceedings of the 26th Annual International Conference on Machine Learning},
  pages={993--1000},
  year={2009}
}

@article{sutton2008convergent,
  title={A convergent O (n) algorithm for off-policy temporal-difference learning with linear function approximation},
  author={Sutton, Richard S and Szepesv{\'a}ri, Csaba and Maei, Hamid Reza},
  journal={Advances in neural information processing systems},
  volume={21},
  number={21},
  pages={1609--1616},
  year={2008},
  publisher={MIT Press}
}

@article{maei2009convergent,
  title={Convergent temporal-difference learning with arbitrary smooth function approximation},
  author={Maei, Hamid and Szepesvari, Csaba and Bhatnagar, Shalabh and Precup, Doina and Silver, David and Sutton, Richard S},
  journal={Advances in neural information processing systems},
  volume={22},
  year={2009}
}

@article{gravell2020learning,
  title={Learning optimal controllers for linear systems with multiplicative noise via policy gradient},
  author={Gravell, Benjamin and Esfahani, Peyman Mohajerin and Summers, Tyler},
  journal={IEEE Transactions on Automatic Control},
  volume={66},
  number={11},
  pages={5283--5298},
  year={2020},
  publisher={IEEE}
}

@inproceedings{mei2020global,
  title={On the global convergence rates of softmax policy gradient methods},
  author={Mei, Jincheng and Xiao, Chenjun and Szepesvari, Csaba and Schuurmans, Dale},
  booktitle={International Conference on Machine Learning},
  pages={6820--6829},
  year={2020},
  organization={PMLR}
}

@article{ju2022model,
  title={A model-free first-order method for linear quadratic regulator with O (1/$\varepsilon$) sampling complexity},
  author={Ju, Caleb and Kotsalis, Georgios and Lan, Guanghui},
  journal={arXiv preprint arXiv:2212.00084},
  year={2022}
}

@article{zhou2023single,
  title={Single timescale actor-critic method to solve the linear quadratic regulator with convergence guarantees},
  author={Zhou, Mo and Lu, Jianfeng},
  journal={Journal of Machine Learning Research},
  volume={24},
  number={222},
  pages={1--34},
  year={2023}
}

@article{zou2019finite,
  title={Finite-sample analysis for sarsa with linear function approximation},
  author={Zou, Shaofeng and Xu, Tengyu and Liang, Yingbin},
  journal={Advances in neural information processing systems},
  volume={32},
  year={2019}
}

@book{bertsekas2012dynamic,
  title={Dynamic programming and optimal control: Volume I},
  author={Bertsekas, Dimitri},
  volume={1},
  year={2012},
  publisher={Athena scientific}
}

@article{wang2021non,
  title={Non-asymptotic analysis for two time-scale TDC with general smooth function approximation},
  author={Wang, Yue and Zou, Shaofeng and Zhou, Yi},
  journal={Advances in Neural Information Processing Systems},
  volume={34},
  pages={9747--9758},
  year={2021}
}

@inproceedings{srikant2019finite,
  title={Finite-time error bounds for linear stochastic approximation and {TD} learning},
  author={Srikant, Rayadurgam and Ying, Lei},
  booktitle={Conference on Learning Theory},
  pages={2803--2830},
  year={2019},
  organization={PMLR}
}

@article{chen2022finite,
  title={Finite-sample analysis of nonlinear stochastic approximation with applications in reinforcement learning},
  author={Chen, Zaiwei and Zhang, Sheng and Doan, Thinh T and Clarke, John-Paul and Maguluri, Siva Theja},
  journal={Automatica},
  volume={146},
  pages={110623},
  year={2022},
  publisher={Elsevier}
}

@inproceedings{xiong2021non,
  title={Non-asymptotic convergence of adam-type reinforcement learning algorithms under markovian sampling},
  author={Xiong, Huaqing and Xu, Tengyu and Liang, Yingbin and Zhang, Wei},
  booktitle={Proceedings of the AAAI Conference on Artificial Intelligence},
  volume={35},
  pages={10460--10468},
  year={2021}
}

@article{borkar1997stochastic,
  title={Stochastic approximation with two time scales},
  author={Borkar, Vivek S},
  journal={Systems \& Control Letters},
  volume={29},
  number={5},
  pages={291--294},
  year={1997},
  publisher={Elsevier}
}

@article{konda1999actor,
  title={Actor-critic--type learning algorithms for Markov decision processes},
  author={Konda, Vijaymohan R and Borkar, Vivek S},
  journal={SIAM Journal on control and Optimization},
  volume={38},
  number={1},
  pages={94--123},
  year={1999},
  publisher={SIAM}
}

@article{haque2023tight,
  title={Tight Finite Time Bounds of Two-Time-Scale Linear Stochastic Approximation with Markovian Noise},
  author={Haque, Shaan Ul and Khodadadian, Sajad and Maguluri, Siva Theja},
  journal={arXiv preprint arXiv:2401.00364},
  year={2023}
}

@article{han2024finite,
  title={Finite-Time Decoupled Convergence in Nonlinear Two-Time-Scale Stochastic Approximation},
  author={Han, Yuze and Li, Xiang and Zhang, Zhihua},
  journal={arXiv preprint arXiv:2401.03893},
  year={2024}
}

@article{olshevsky2023small,
  title={A small gain analysis of single timescale actor critic},
  author={Olshevsky, Alex and Gharesifard, Bahman},
  journal={SIAM Journal on Control and Optimization},
  volume={61},
  number={2},
  pages={980--1007},
  year={2023},
  publisher={SIAM}
}

@article{mokkadem2006convergence,
  title={Convergence rate and averaging of nonlinear two-time-scale stochastic approximation algorithms},
  author={Mokkadem, Abdelkader and Pelletier, Mariane},
  journal={Ann. Appl. Probab.},
  volume={16},
  number={1},
  pages={1671--1702},
  year={2006}
}

@article{ma2020variance,
  title={Variance-reduced off-policy TDC learning: Non-asymptotic convergence analysis},
  author={Ma, Shaocong and Zhou, Yi and Zou, Shaofeng},
  journal={Advances in neural information processing systems},
  volume={33},
  pages={14796--14806},
  year={2020}
}
\clearpage
\appendix

\section{Analysis for Strongly Convex Functions}\label{sec:proof_stronglyconvex}

The exact choice of the step size parameters is made such that $\alpha_k\leq\beta_k\leq\lambda_k\leq1/4$, $c_{\alpha}\geq\frac{8}{\mu_h}$, and
\begin{align*}
&\alpha_k\leq\min\{\frac{\lambda_k}{\mu_h},\frac{\mu_h}{6(L+1)^4}, \frac{8\mu_G}{\mu_h}\beta_k, \frac{\mu_h\mu_G}{152(L+1)^6}\beta_k, \frac{\mu_G\beta_k}{8(8L^4+\frac{9L^2}{\mu_G})}, \frac{\lambda_k}{4}(48L^2+\frac{10L}{\mu_G})^{-1}\},\\
&\beta_k\leq \min\{\frac{1}{240(L+1)^6},\frac{1}{\mu_G},\frac{\lambda_k}{4}(56L^2+\frac{9}{2\mu_h})^{-1},\frac{\lambda_k}{4}(48L^2+\frac{10L}{\mu_G})^{-1},\frac{\mu_G}{168L^4}\},
\end{align*}
which requires $\beta_0$ to be selected sufficiently small and $\tau$ sufficiently large. Note that parameters $c_{\alpha}$, $\beta_0$, and $\tau$ always exist and only depend on the structural constants of $h$ and $G$.

% $\alpha_k\leq\beta_k\leq\lambda_k$, $\alpha_k\leq\frac{\lambda_k}{\mu_h}$, $\alpha_k\leq\frac{\mu_h}{6(L+1)^4}$, $\alpha_k\leq\frac{8\mu_G}{\mu_h}\beta_k$ ,$\beta_k\leq \frac{1}{240(L+1)^6}$, $\beta_k\leq\frac{1}{\mu_G}$,

% $\beta_k\leq\frac{1}{\mu_G}$, $\beta_k\leq\frac{\lambda_k}{4}(56L^2+\frac{9}{2\mu_h})^{-1}$, $\beta_k\leq\frac{\lambda_k}{4}(48L^2+\frac{10L}{\mu_G})^{-1}$, $\alpha_k\leq\frac{\mu_h\mu_G}{152(L+1)^6}\beta_k$, $\beta_k\leq \frac{1}{240(L+1)^6}$, $\beta_k\leq\frac{\mu_G}{168L^4}$, and $\alpha_k\leq\frac{\mu_G\beta_k}{8(8L^4+\frac{9L^2}{\mu_G})}$

The proof of Theorem~\ref{thm:stronglyconvex} is built on the lemmas below. Their proofs can be found at the end of the section.

\begin{lem}\label{lem:Lipschitz:stronglyconvex}
Under Assumption~\ref{assump:Lipschitz}-\ref{assump:stronglymonotone_G}, the sequence of variables $\{\theta_k,\omega_k,f_k,g_k\}$ satisfy for all $k$
\begin{align*}
    \|f_k\|&\leq \|\Delta f_k\| + L\sqrt{y_k} + L(L+1)\sqrt{z_k},\\
    \|g_k\|&\leq \|\Delta g_k\|+L\sqrt{y_k},\\
    \|\theta_{k+1}-\theta_k\|&\leq\alpha_k\left(\|\Delta f_k\| + L\sqrt{y_k} + L(L+1)\sqrt{z_k}\right),\\
    \|\omega_{k+1}-\omega_k\| & \leq\beta_k\left(\|\Delta g_k\|+ L\sqrt{y_k}\right).
\end{align*}
\end{lem}

\begin{lem}\label{lem:Delta_f:stronglyconvex}
Suppose that the step size $\lambda_k$ satisfy $\lambda_k \leq 1/4$ for all $k$.
Then, under Assumption~\ref{assump:Lipschitz}-\ref{assump:stronglymonotone_G}, we have the following bound on the iteration-wise convergence of $\mathbb{E}[\|\Delta f_{k}\|^2]$.
\begin{align*}
\mathbb{E}[\|\Delta f_{k+1}\|^2]
&\leq (1-\lambda_k)\mathbb{E}[\|\Delta f_k\|^2]-(\frac{\lambda_k}{4}-\frac{24L^2\beta_k^2}{\lambda_k})\mathbb{E}[\|\Delta f_k\|^2]\notag\\
&\hspace{20pt}+\frac{4}{\lambda_k}\mathbb{E}\left[6L^2\beta_k^2\|\Delta g_k\|^2+10L^4\beta_k^2 y_k + 6(L+1)^6 \alpha_k^2 z_k\right] +2B\lambda_k^2.
\end{align*}
\end{lem}

\begin{lem}\label{lem:Delta_g:stronglyconvex}
Suppose that the step size $\lambda_k$ satisfy $\lambda_k \leq 1/4$ for all $k$.
Then, under Assumption~\ref{assump:Lipschitz}-\ref{assump:stronglymonotone_G}, we have the following bound on the iteration-wise convergence of $\mathbb{E}[\|\Delta g_{k}\|^2]$
\begin{align*}
\mathbb{E}[\|\Delta g_{k+1}\|^2]
&\leq (1-\lambda_k)\mathbb{E}[\|\Delta g_k\|^2]-(\frac{\lambda_k}{4}-\frac{24L^2\beta_k^2}{\lambda_k})\mathbb{E}[\|\Delta g_k\|^2]\notag\\
&\hspace{20pt}+\frac{4}{\lambda_k}\mathbb{E}\left[6L^2\beta_k^2\|\Delta f_k\|^2+10L^4\beta_k^2 y_k + 6(L+1)^6 \alpha_k^2 z_k\right]+2B\lambda_k^2.
\end{align*}
\end{lem}

\begin{lem}\label{lem:x_k:stronglyconvex}
Suppose the step size $\alpha_k$ is chosen such that $\alpha_k\leq \frac{\mu_h}{6(L+1)^4}$ and $\alpha_k\leq\frac{1}{2\mu_h}$ for all $k$.
Under Assumption~\ref{assump:Lipschitz}-\ref{assump:stronglyconvex}, we have
\begin{align*}
z_{k+1}\leq (1-\frac{\mu_h\alpha_k}{4})z_k+\frac{9L^2\alpha_k}{\mu_h}y_k+\frac{9\alpha_k}{2\mu_h}\|\Delta f_k\|^2.
\end{align*}
\end{lem}

\begin{lem}\label{lem:y_k:stronglyconvex}
Suppose that the step sizes satisfy $\alpha_k\leq\beta_k\leq\frac{1}{72(L+1)^6}$ and $\beta_k\leq1/\mu_G$ for all $k$.
Then, under Assumption~\ref{assump:Lipschitz}-\ref{assump:stronglymonotone_G}, we have
\begin{align*}
y_{k+1}
&\leq (1-\mu_G\beta_k)y_k-\left(\frac{\mu_G}{4}\beta_k -8L^4\alpha_k-L^2\beta_k^2\right)y_k\notag\\
&\hspace{20pt}+\left(\frac{19(L+1)^6\alpha_k^2}{\mu_G\beta_k}+3(L+1)^6\alpha_k\beta_k\right)z_k+\left(\frac{8\beta_k}{\mu_G}+2L\beta_k^2\right)\|\Delta g_k\|^2+8L\alpha_k\|\Delta f_k\|^2.
\end{align*}
\end{lem}

\subsection{Proof of Theorem~\ref{thm:stronglyconvex}}

We consider the Lyapunov function 
\begin{align*}
V_{k} \triangleq \mathbb{E}[\|\Delta f_k\|^2 + \|\Delta g_{k}\|^2 + z_k + y_k].   
\end{align*}

Combining Lemma~\ref{lem:Delta_f:stronglyconvex}-\ref{lem:y_k:stronglyconvex}, we get
\begin{align}
V_{k+1}&=\mathbb{E}[\|\Delta f_{k+1}\|^2 + \|\Delta g_{k+1}\|^2 + z_{k+1} + y_{k+1}]\notag\\
&\leq (1-\lambda_k)\mathbb{E}[\|\Delta f_k\|^2]-(\frac{\lambda_k}{4}-\frac{24L^2\beta_k^2}{\lambda_k})\mathbb{E}[\|\Delta f_k\|^2]\notag\\
&\hspace{20pt}+\frac{4}{\lambda_k}\mathbb{E}\left[6L^2\beta_k^2\|\Delta g_k\|^2+10L^4\beta_k^2 y_k + 6(L+1)^6 \alpha_k^2 z_k\right] +2B\lambda_k^2\notag\\
&\hspace{20pt}+(1-\lambda_k)\mathbb{E}[\|\Delta g_k\|^2]-(\frac{\lambda_k}{4}-\frac{24L^2\beta_k^2}{\lambda_k})\mathbb{E}[\|\Delta g_k\|^2]\notag\\
&\hspace{20pt}+\frac{4}{\lambda_k}\mathbb{E}\left[6L^2\beta_k^2\|\Delta f_k\|^2+10L^4\beta_k^2 y_k + 6(L+1)^6 \alpha_k^2 z_k\right]+2B\lambda_k^2\notag\\
&\hspace{20pt}+\mathbb{E}[(1-\frac{\mu_h\alpha_k}{4})z_k+\frac{9L^2\alpha_k}{\mu_h}y_k+\frac{9\alpha_k}{2\mu_h}\|\Delta f_k\|^2]\notag\\
&\hspace{20pt}+\mathbb{E}[(1-\mu_G\beta_k)y_k-\left(\frac{\mu_G}{4}\beta_k -8L^4\alpha_k-L^2\beta_k^2\right)y_k\notag\\
&\hspace{20pt}+\left(\frac{19(L+1)^6\alpha_k^2}{\mu_G\beta_k}+3(L+1)^6\alpha_k\beta_k\right)z_k+\left(\frac{8\beta_k}{\mu_G}+2L\beta_k^2\right)\|\Delta g_k\|^2+8L\alpha_k\|\Delta f_k\|^2]\notag\\
&\leq (1-\lambda_k)\mathbb{E}[\|\Delta f_k\|^2]-(\frac{\lambda_k}{4}-\frac{48L^2\beta_k^2}{\lambda_k}-\frac{9\alpha_k}{2\mu_h}-8L\alpha_k)\mathbb{E}[\|\Delta f_k\|^2]\notag\\
&\hspace{20pt}+(1-\lambda_k)\mathbb{E}[\|\Delta g_k\|^2]-(\frac{\lambda_k}{4}-\frac{48L^2\beta_k^2}{\lambda_k}-\frac{8\beta_k}{\mu_G}-2L\beta_k^2)\mathbb{E}[\|\Delta g_k\|^2]\notag\\
&\hspace{20pt}+(1-\frac{\mu_h\alpha_k}{8})\mathbb{E}[z_k]-(\frac{\mu_h\alpha_k}{8}-15(L+1)^6 \alpha_k\beta_k-\frac{19(L+1)^6\alpha_k^2}{\mu_G\beta_k})\mathbb{E}[z_k]\notag\\
&\hspace{20pt}+(1-\mu_G\beta_k)\mathbb{E}[y_k]-\left(\frac{\mu_G}{4}\beta_k -8L^4\alpha_k-21L^4\beta_k^2-\frac{9L^2\alpha_k}{\mu_h}\right)\mathbb{E}[y_k]+4B\lambda_k^2\notag\\
&\leq (1-\frac{\mu_h\alpha_k}{8})V_k+4B\lambda_k^2-(\frac{\lambda_k}{4}-\frac{48L^2\beta_k^2}{\lambda_k}-\frac{9\alpha_k}{2\mu_h}-8L\alpha_k)\mathbb{E}[\|\Delta f_k\|^2]\notag\\
&\hspace{20pt}-(\frac{\lambda_k}{4}-\frac{48L^2\beta_k^2}{\lambda_k}-\frac{8\beta_k}{\mu_G}-2L\beta_k^2)\mathbb{E}[\|\Delta g_k\|^2]\notag\\
&\hspace{20pt}-(\frac{\mu_h\alpha_k}{8}-15(L+1)^6 \alpha_k\beta_k-\frac{19(L+1)^6\alpha_k^2}{\mu_G\beta_k})\mathbb{E}[z_k]\notag\\
&\hspace{20pt}-\left(\frac{\mu_G}{4}\beta_k -8L^4\alpha_k-21L^4\beta_k^2-\frac{9L^2\alpha_k}{\mu_h}\right)\mathbb{E}[y_k],\label{thm:stronglyconvex:proof_eq1}
\end{align}
where the third inequality simplifies the terms with the conditions $\alpha_k\leq\frac{\lambda_k}{\mu_h}$ and $\alpha_k\leq\frac{8\mu_G\beta_k}{\mu_h}$.

The terms on the right hand side of \eqref{thm:stronglyconvex:proof_eq1} except for the first two are non-positive under the step size conditions $\beta_k\leq\frac{1}{\mu_G}$, $\beta_k\leq\frac{\lambda_k}{4}(56L^2+\frac{9}{2\mu_h})^{-1}$, $\beta_k\leq\frac{\lambda_k}{4}(48L^2+\frac{10L}{\mu_G})^{-1}$, $\alpha_k\leq\frac{\mu_h\mu_G}{152(L+1)^6}\beta_k$, $\beta_k\leq \frac{1}{240(L+1)^6}$, $\beta_k\leq\frac{\mu_G}{168L^4}$, and $\alpha_k\leq\frac{\mu_G\beta_k}{8(8L^4+\frac{9L^2}{\mu_G})}$. This leads to
\begin{align*}
    V_{k+1}\leq(1-\frac{\mu_h\alpha_k}{8})V_k+4B\lambda_k^2\leq(1-\frac{2}{k+\tau+1})V_k+4B\lambda_k^2,
\end{align*}
where the second inequality follows from $c_{\alpha}\geq\frac{16}{\mu_h}$.

Multiplying by $(k+\tau+1)^2$ on both sides, we have
\begin{align*}
    (k+\tau+1)^2 V_{k+1}&\leq(k+\tau+1)(k+\tau-1)V_k+4B(k+\tau+1)^2\lambda_k^2\notag\\
    &\leq (k+\tau)^2 V_k+\frac{B(k+\tau+1)^2}{4(k+1)^2}\notag\\
    &\leq (k+\tau)^2 V_k+(\frac{B}{2}+\frac{B\tau^2}{2})\notag\\
    &\leq \tau^2 V_0+B\tau^2(k+1)
\end{align*}

Dividing by $(k+\tau+1)^2$ now, we have
\begin{align*}
    V_{k+1}\leq\frac{\tau^2 V_0}{(k+\tau+1)^2}+\frac{B\tau^2(k+1)}{(k+\tau+1)^2}\leq\frac{\tau^2 V_0}{(k+\tau+1)^2}+\frac{B\tau^2}{(k+\tau+1)},
\end{align*}
which obviously implies
\begin{align*}
\mathbb{E}[\|\theta_{k+1}-\theta^{\star}\|^2]\leq\frac{\tau^2 V_0}{(k+\tau+1)^2}+\frac{B\tau^2}{(k+\tau+1)}.
\end{align*}

\qed

\subsection{Proof of Lemma~\ref{lem:Lipschitz:stronglyconvex}}

By definition,
\begin{align*}
    F(\theta^{\star}, \omega^{\star}(\theta^{\star}))=\nabla h(\theta^{\star})=0,
\end{align*}
which implies
\begin{align*}
    \|f_k\|&=\|\Delta f_k+F(\theta_k,\omega_k)-F(\theta_k,\omega^{\star}(\theta_k))+F(\theta_k,\omega^{\star}(\theta_k))-F(\theta^{\star},\omega^{\star}(\theta^{\star}))\|\\
    &\leq\|\Delta f_k\|+\|F(\theta_k,\omega_k)-F(\theta_k,\omega^{\star}(\theta_k))\|+\|F(\theta_k,\omega^{\star}(\theta_k))-F(\theta^{\star},\omega^{\star}(\theta^{\star}))\|\\
    &\leq \|\Delta f_k\| + L\|\omega_k-\omega^{\star}(\theta_k)\| + L(L+1)\|\theta_k-\theta^{\star}\|\notag\\
    &\leq \|\Delta f_k\| + L\sqrt{y_k} + L(L+1)\sqrt{z_k},
\end{align*}
where the second inequality follows from the Lipschitz continuity of $F$.

Similarly, since $G(\theta,\omega^{\star}(\theta))=0$ for any $\theta$, we can derive
\begin{align*}
    \|g_k\|&=\|\Delta g_k+G(\theta_k,\omega_k)-G(\theta_k,\omega^{\star}(\theta_k))\|\notag\\
    &\leq\|\Delta g_k\|+\|G(\theta_k,\omega_k)-G(\theta_k,\omega^{\star}(\theta_k))\|\notag\\
    &\leq \|\Delta g_k\|+L\sqrt{y_k}.
\end{align*}

The bounds on $\|\theta_{k+1}-\theta_k\|$ and $\|\omega_{k+1}-\omega_k\|$ easily follow from the two inequalities above and \eqref{alg:update_decision}.

\qed

\subsection{Proof of Lemma~\ref{lem:Delta_f:stronglyconvex}}\label{sec:proof:lem:Delta_f:stronglyconvex}

By the update rule in \eqref{eq:update_auxiliary}, we have
\begin{align*}
\Delta f_{k+1}&=f_{k+1}-F(\theta_{k+1},\omega_{k+1})\notag\\
&=(1-\lambda_k)f_k+\lambda_k F(\theta_k,\omega_k,X_k)-F(\theta_{k+1},\omega_{k+1})\notag\\
&=(1-\lambda_k)f_k+\lambda_k F(\theta_k,\omega_k)-F(\theta_{k+1},\omega_{k+1})+\lambda_k\left(F(\theta_k,\omega_k)-F(\theta_k,\omega_k,X_k)\right)\notag\\
&=(1-\lambda_k)\Delta f_k + F(\theta_k,\omega_k)-F(\theta_{k+1},\omega_{k+1})+\lambda_k\left(F(\theta_k,\omega_k)-F(\theta_k,\omega_k,X_k)\right).
\end{align*}
This leads to 
\begin{align}
&\|\Delta f_{k+1}\|^2\notag\\
&=(1-\lambda_k)^2\|\Delta f_k\|^2 +\Big\|\left(F(\theta_k,\omega_k)-F(\theta_{k+1},\omega_{k+1})\right)+\lambda_k\left(F(\theta_k,\omega_k)-F(\theta_k,\omega_k,X_k)\right)\Big\|^2\notag\\
&\hspace{20pt}+2(1-\lambda_k)\Delta f_k^{\top}\left(F(\theta_k,\omega_k)-F(\theta_{k+1},\omega_{k+1})\right)\notag\\
&\hspace{20pt}+2(1-\lambda_k)\lambda_k\Delta f_k^{\top}\left(F(\theta_k,\omega_k)-F(\theta_k,\omega_k,X_k)\right)\notag\\
&\leq (1-\lambda_k)^2\|\Delta f_k\|^2 +2\|F(\theta_k,\omega_k)-F(\theta_{k+1},\omega_{k+1})\|^2+2\lambda_k^2\|F(\theta_k,\omega_k)-F(\theta_k,\omega_k,X_k)\|^2\notag\\
&\hspace{20pt}+\frac{\lambda_k}{2}\|\Delta f_k\|^2+\frac{2}{\lambda_k}\left\|F(\theta_k,\omega_k)-F(\theta_{k+1},\omega_{k+1})\right\|^2\notag\\
&\hspace{20pt}+2(1-\lambda_k)\lambda_k\Delta f_k^{\top}\left(F(\theta_k,\omega_k)-F(\theta_k,\omega_k,X_k)\right)\notag\\
&\leq (1-\lambda_k)\|\Delta f_k\|^2 +2\lambda_k^2\|F(\theta_k,\omega_k)-F(\theta_k,\omega_k,X_k)\|^2\notag\\
&\hspace{20pt}+(\lambda_k^2-\frac{\lambda_k}{2})\|\Delta f_k\|^2+\frac{4}{\lambda_k}\left\|F(\theta_k,\omega_k)-F(\theta_{k+1},\omega_{k+1})\right\|^2\notag\\
&\hspace{20pt}+2(1-\lambda_k)\lambda_k\Delta f_k^{\top}\left(F(\theta_k,\omega_k)-F(\theta_k,\omega_k,X_k)\right),\label{lem:Delta_f:proof_eq1:stronglyconvex}
\end{align}
where the last inequality follows from $\lambda_k\leq 1$.

The second term on the right hand side of \eqref{lem:Delta_f:proof_eq1:stronglyconvex} is bounded in expectation, i.e.
\[\mathbb{E}_{X_k\sim\xi}[\|F(\theta_k,\omega_k,X_k)-F(\theta_k,\omega_k)\|^2]\leq B.\]
% $\mathbb{E}_{X_k\sim\xi}[\|F(\theta_k,\omega_k,X_k)-F(\theta_k,\omega_k)\|^2]\leq B$.
In addition, the last term of \eqref{lem:Delta_f:proof_eq1:stronglyconvex} is zero in expectation since $\mathbb{E}[\Delta f_k^{\top}(F(\theta_k,\omega_k)-F(\theta_k,\omega_k,X_k))\mid \Hcal_{k-1}]=\Delta f_k^{\top}\mathbb{E}[(F(\theta_k,\omega_k)-F(\theta_k,\omega_k,X_k))\mid \Hcal_{k-1}]=0$. As a result, we have
\begin{align}
\mathbb{E}[\|\Delta f_{k+1}\|^2]&\leq (1-\lambda_k)\mathbb{E}[\|\Delta f_k\|^2]-(\frac{\lambda_k}{2}-\lambda_k^2)\mathbb{E}[\|\Delta f_k\|^2]\notag\\
&\hspace{20pt}+\frac{4}{\lambda_k}\mathbb{E}[\left\|F(\theta_k,\omega_k)-F(\theta_{k+1},\omega_{k+1})\right\|^2]+2B\lambda_k^2.\label{lem:Delta_f:proof_eq2:stronglyconvex}
\end{align}

By Lemma~\ref{lem:Lipschitz:stronglyconvex} and the Lipschitz condition of $F$,
\begin{align*}
    &\|F(\theta_k,\omega_k)-F(\theta_{k+1},\omega_{k+1})\|^2\\
    & \leq 2L^2\|\theta_{k+1}-\theta_k\|^2+ 2L^2\|\omega_{k+1}-\omega_k\|^2\notag\\
    &\leq 6L^2\alpha_k^2(\|\Delta f_k\|^2+L^2 y_k+(L+1)^4 z_k)+4L^2\beta_k^2(\|\Delta g_k\|^2+L^2 y_k)\notag\\
    &\leq 6L^2\beta_k^2\|\Delta f_k\|^2 + 6L^2\beta_k^2\|\Delta g_k\|^2+10L^4\beta_k^2 y_k + 6(L+1)^6 \alpha_k^2 z_k.
\end{align*}
Plugging this bound into \eqref{lem:Delta_f:proof_eq2:stronglyconvex}, we get
\begin{align*}
&\mathbb{E}[\|\Delta f_{k+1}\|^2]\notag\\
&\leq (1-\lambda_k)\mathbb{E}[\|\Delta f_k\|^2]-(\frac{\lambda_k}{2}-\lambda_k^2)\mathbb{E}[\|\Delta f_k\|^2] \notag\\
&\hspace{20pt}+\frac{4}{\lambda_k}\mathbb{E}\left[6L^2\beta_k^2\|\Delta f_k\|^2 + 6L^2\beta_k^2\|\Delta g_k\|^2+10L^4\beta_k^2 y_k + 6(L+1)^6 \alpha_k^2 z_k\right]+ 2B\lambda_k^2\notag\\
&\leq (1-\lambda_k)\mathbb{E}[\|\Delta f_k\|^2]-(\frac{\lambda_k}{4}-\frac{24L^2\beta_k^2}{\lambda_k})\mathbb{E}[\|\Delta f_k\|^2]\notag\\
&\hspace{20pt}+\frac{4}{\lambda_k}\mathbb{E}\left[6L^2\beta_k^2\|\Delta g_k\|^2+10L^4\beta_k^2 y_k + 6(L+1)^6 \alpha_k^2 z_k\right] +2B\lambda_k^2,
\end{align*}
where the last inequality follows from $\lambda_k \leq 1/4$.

\qed

\subsection{Proof of Lemma~\ref{lem:Delta_g:stronglyconvex}}

By the update rule in \eqref{eq:update_auxiliary}, we have
\begin{align*}
\Delta g_{k+1}&=g_{k+1}-G(\theta_{k+1},\omega_{k+1})\notag\\
&=(1-\lambda_k)g_k+\lambda_k G(\theta_k,\omega_k,X_k)-G(\theta_{k+1},\omega_{k+1})\notag\\
&=(1-\lambda_k)g_k+\lambda_k G(\theta_k,\omega_k)-G(\theta_{k+1},\omega_{k+1})+\lambda_k\left(G(\theta_k,\omega_k)-G(\theta_k,\omega_k,X_k)\right)\notag\\
&=(1-\lambda_k)\Delta g_k + G(\theta_k,\omega_k)-G(\theta_{k+1},\omega_{k+1})+\lambda_k\left(G(\theta_k,\omega_k)-G(\theta_k,\omega_k,X_k)\right).
\end{align*}
This leads to 
\begin{align}
&\|\Delta g_{k+1}\|^2\notag\\
&=(1-\lambda_k)^2\|\Delta g_k\|^2 +\Big\|\left(G(\theta_k,\omega_k)-G(\theta_{k+1},\omega_{k+1})\right)+\lambda_k\left(G(\theta_k,\omega_k)-G(\theta_k,\omega_k,X_k)\right)\Big\|^2\notag\\
&\hspace{20pt}+2(1-\lambda_k)\Delta g_k^{\top}\left(G(\theta_k,\omega_k)-G(\theta_{k+1},\omega_{k+1})\right)\notag\\
&\hspace{20pt}+2(1-\lambda_k)\lambda_k\Delta g_k^{\top}\left(G(\theta_k,\omega_k)-G(\theta_k,\omega_k,X_k)\right)\notag\\
&\leq (1-\lambda_k)^2\|\Delta g_k\|^2 +2\|G(\theta_k,\omega_k)-G(\theta_{k+1},\omega_{k+1})\|^2+2\lambda_k^2\|G(\theta_k,\omega_k)-G(\theta_k,\omega_k,X_k)\|^2\notag\\
&\hspace{20pt}+\frac{\lambda_k}{2}\|\Delta g_k\|^2+\frac{2}{\lambda_k}\left\|G(\theta_k,\omega_k)-G(\theta_{k+1},\omega_{k+1})\right\|^2\notag\\
&\hspace{20pt}+2(1-\lambda_k)\lambda_k\Delta g_k^{\top}\left(G(\theta_k,\omega_k)-G(\theta_k,\omega_k,X_k)\right)\notag\\
&\leq (1-\lambda_k)\|\Delta g_k\|^2 +2\lambda_k^2\|G(\theta_k,\omega_k)-G(\theta_k,\omega_k,X_k)\|^2\notag\\
&\hspace{20pt}+(\lambda_k^2-\frac{\lambda_k}{2})\|\Delta g_k\|^2+\frac{4}{\lambda_k}\left\|G(\theta_k,\omega_k)-G(\theta_{k+1},\omega_{k+1})\right\|^2\notag\\
&\hspace{20pt}+2(1-\lambda_k)\lambda_k\Delta g_k^{\top}\left(G(\theta_k,\omega_k)-G(\theta_k,\omega_k,X_k)\right).\label{lem:Delta_g:proof_eq1}
\end{align}

The second term on the right hand side of \eqref{lem:Delta_g:proof_eq1} is bounded in expectation, i.e.
\[\mathbb{E}_{X_k\sim\xi}[\|G(\theta_k,\omega_k,X_k)-G(\theta_k,\omega_k)\|^2]\leq B.\]
% since $\mathbb{E}_{X_k\sim\xi}[\|G(\theta_k,\omega_k,X_k)-G(\theta_k,\omega_k)\|^2]\leq B$.
In addition, the last term of \eqref{lem:Delta_f:proof_eq1:stronglyconvex} is zero in expectation since $\mathbb{E}[\Delta g_k^{\top}(G(\theta_k,\omega_k)-G(\theta_k,\omega_k,X_k))\mid \Hcal_{k-1}]=\Delta g_k^{\top}\mathbb{E}[(G(\theta_k,\omega_k)-G(\theta_k,\omega_k,X_k))\mid \Hcal_{k-1}]=0$. As a result, we have
\begin{align}
\mathbb{E}[\|\Delta g_{k+1}\|^2]&\leq (1-\lambda_k)\mathbb{E}[\|\Delta g_k\|^2]-(\frac{\lambda_k}{2}-\lambda_k^2)\mathbb{E}[\|\Delta g_k\|^2]\notag\\
&\hspace{20pt}+\frac{4}{\lambda_k}\mathbb{E}[\left\|G(\theta_k,\omega_k)-G(\theta_{k+1},\omega_{k+1})\right\|^2]+2B\lambda_k^2.\label{lem:Delta_g:proof_eq2}
\end{align}

By Lemma~\ref{lem:Lipschitz:stronglyconvex} and the Lipschitz condition of $G$,
\begin{align*}
    &\|G(\theta_k,\omega_k)-G(\theta_{k+1},\omega_{k+1})\|^2\notag\\
    & \leq 2L^2\|\theta_{k+1}-\theta_k\|^2+ 2L^2\|\omega_{k+1}-\omega_k\|^2\notag\\
    &\leq 6L^2\alpha_k^2(\|\Delta f_k\|^2+L^2 y_k+(L+1)^4 z_k)+4L^2\beta_k^2(\|\Delta g_k\|^2+L^2 y_k)\notag\\
    &\leq 6L^2\beta_k^2\|\Delta f_k\|^2 + 6L^2\beta_k^2\|\Delta g_k\|^2+10L^4\beta_k^2 y_k + 6(L+1)^6 \alpha_k^2 z_k.
\end{align*}
Plugging this bound into \eqref{lem:Delta_g:proof_eq2}, we get
\begin{align*}
\mathbb{E}[\|\Delta g_{k+1}\|^2]&\leq (1-\lambda_k)\mathbb{E}[\|\Delta g_k\|^2]-(\frac{\lambda_k}{2}-\lambda_k^2)\mathbb{E}[\|\Delta g_k\|^2]+2B\lambda_k^2\notag\\
&\hspace{20pt}+\frac{4}{\lambda_k}\mathbb{E}\left[6L^2\beta_k^2\|\Delta f_k\|^2 + 6L^2\beta_k^2\|\Delta g_k\|^2+10L^4\beta_k^2 y_k + 6(L+1)^6 \alpha_k^2 z_k\right]\notag\\
&\leq (1-\lambda_k)\mathbb{E}[\|\Delta g_k\|^2]-(\frac{\lambda_k}{4}-\frac{24L^2\beta_k^2}{\lambda_k})\mathbb{E}[\|\Delta g_k\|^2]\notag\\
&\hspace{20pt}+\frac{4}{\lambda_k}\mathbb{E}\left[6L^2\beta_k^2\|\Delta f_k\|^2+10L^4\beta_k^2 y_k + 6(L+1)^6 \alpha_k^2 z_k\right]+2B\lambda_k^2,
\end{align*}
where the last inequality follows from $\lambda_k \leq 1/4$.

\qed

\subsection{Proof of Lemma~\ref{lem:x_k:stronglyconvex}}

We have by the update rule in \eqref{alg:update_decision} and the definition of the residual variables in \eqref{eq:residuals}
\begin{align}
    z_{k+1}&=\|\theta_{k+1}-\theta^{\star}\|^2\notag\\
    &=\|\theta_{k}-\alpha_k f_k-\theta^{\star}\|^2\notag\\
    &=\|\theta_{k} -\alpha_k F(\theta_k,\omega_k)-\alpha_k \Delta f_k-\theta^{\star}\|^2\notag\\
    &=\|\theta_{k}-\theta^{\star} -\alpha_k F(\theta_k,\omega_k)\|^2-2\alpha_k\langle \theta_{k}-\theta^{\star} -\alpha_k F(\theta_k,\omega_k), \Delta f_k\rangle+\alpha_k^2 \|\Delta f_k\|^2.\label{lem:x_k:stronglyconvex:proof_eq1}
\end{align}

To bound the first term on the right hand side of \eqref{lem:x_k:stronglyconvex:proof_eq1},
\begin{align*}
&\|\theta_{k}-\theta^{\star} -\alpha_k F(\theta_k,\omega_k)\|^2\notag\\
&=\|\theta_{k}-\theta^{\star} -\alpha_k \left(F(\theta_k,\omega^{\star}(\theta_k))-F(\theta^{\star},\omega^{\star}(\theta^{\star}))\right)+\alpha_k\left(F(\theta_k,\omega^{\star}(\theta_k))-F(\theta_k,\omega_k)\right)\|^2\notag\\
&=z_k+\alpha_k^2\|F(\theta_k,\omega^{\star}(\theta_k))-F(\theta^{\star},\omega^{\star}(\theta^{\star}))\|^2+\alpha_k^2\|F(\theta_k,\omega^{\star}(\theta_k))-F(\theta_k,\omega_k)\|^2\notag\\
&\hspace{20pt}-2\alpha_k\langle\theta_{k}-\theta^{\star},F(\theta_k,\omega^{\star}(\theta_k))-F(\theta^{\star},\omega^{\star}(\theta^{\star}))\rangle\notag\\
&\hspace{20pt}+2\alpha_k\langle\theta_{k}-\theta^{\star},F(\theta_k,\omega^{\star}(\theta_k))-F(\theta_k,\omega_k)\rangle
\notag\\
&\hspace{20pt}-2\alpha_k^2\langle F(\theta_k,\omega^{\star}(\theta_k))-F(\theta^{\star},\omega^{\star}(\theta^{\star}),F(\theta_k,\omega^{\star}(\theta_k))-F(\theta_k,\omega_k)\rangle\notag\\
&\leq z_k+\alpha_k^2\|F(\theta_k,\omega^{\star}(\theta_k))-F(\theta^{\star},\omega^{\star}(\theta^{\star}))\|^2+\alpha_k^2\|F(\theta_k,\omega^{\star}(\theta_k))-F(\theta_k,\omega_k)\|^2\notag\\
&\hspace{20pt}-2\mu_h\alpha_k z_k+\mu_h\alpha_k z_k+\frac{\alpha_k}{\mu_h}\|F(\theta_k,\omega^{\star}(\theta_k))-F(\theta_k,\omega_k)\|^2\notag\\
&\hspace{20pt}+\alpha_k^2 z_k+\alpha_k^2\|F(\theta_k,\omega^{\star}(\theta_k))-F(\theta_k,\omega_k)\|^2\notag\\
&\leq (1-\mu_h\alpha_k+\alpha_k^2)z_k+2L^2(L+1)^2\alpha_k^2 z_k+(2L^2\alpha_k^2+\frac{L^2\alpha_k}{\mu_h}) y_k\notag\\
&\leq (1-\mu_h\alpha_k+3(L+1)^4\alpha_k^2)z_k+(2L^2\alpha_k^2+\frac{L^2\alpha_k}{\mu_h}) y_k,
\end{align*}
where the first equality follows from $F(\theta^{\star},\omega^{\star}(\theta^{\star}))=0$, the first inequality plugs in \eqref{eq:stronglyconvex} and uses the fact that $2\|\va\|\|\vb\|\leq c\|\va\|^2+\frac{1}{c}\|\vb\|^2$ for any vector $\va,\vb$ and scalar $c>0$, and the second inequality follows from the Lipschitz continuity of $F$ and $\omega^{\star}$. Choosing the step size such that $\alpha_k\leq \frac{\mu_h}{6(L+1)^4}$ and $\alpha_k\leq\frac{1}{2\mu_h}$ results in
\begin{align}
\|\theta_{k}-\theta^{\star} -\alpha_k F(\theta_k,\omega_k)\|^2\leq (1-\frac{\mu_h\alpha_k}{2})z_k+\frac{2L^2\alpha_k}{\mu_h}y_k.\label{lem:x_k:stronglyconvex:proof_eq2}
\end{align}

The second term on the right hand side of \eqref{lem:x_k:stronglyconvex:proof_eq1} can be treated based on \eqref{lem:x_k:stronglyconvex:proof_eq2}
\begin{align}
&-2\alpha_k\langle \theta_{k}-\theta^{\star} -\alpha_k F(\theta_k,\omega_k), \Delta f_k\rangle\notag\\
&\leq \frac{\mu_h\alpha_k}{4}\|\theta_{k}-\theta^{\star} -\alpha_k F(\theta_k,\omega_k)\|^2+\frac{4\alpha_k}{\mu_h}\|\Delta f_k\|^2\notag\\
&\leq\frac{\mu_h\alpha_k}{4}(1-\frac{\mu_h\alpha_k}{2})z_k+\frac{\mu_h\alpha_k}{4}\cdot\frac{2L^2\alpha_k}{\mu_h}y_k+\frac{4\alpha_k}{\mu_h}\|\Delta f_k\|^2\notag\\
&\leq(\frac{\mu_h\alpha_k}{4}-\frac{\mu_h^2\alpha_k^2}{8})z_k+\frac{L^2\alpha_k^2}{2} y_k+\frac{4\alpha_k}{\mu_h}\|\Delta f_k\|^2.\label{lem:x_k:stronglyconvex:proof_eq3}
\end{align}

Plugging \eqref{lem:x_k:stronglyconvex:proof_eq2} and \eqref{lem:x_k:stronglyconvex:proof_eq3} into \eqref{lem:x_k:stronglyconvex:proof_eq1}, we get
\begin{align*}
z_{k+1}&\leq(1-\frac{\mu_h\alpha_k}{2})z_k+\frac{2L^2\alpha_k}{\mu_h}y_k\notag\\
&\hspace{20pt}+(\frac{\mu_h\alpha_k}{4}-\frac{\mu_h^2\alpha_k^2}{8})z_k+\frac{L^2\alpha_k^2}{2}y_k+\frac{4\alpha_k}{\mu_h}\|\Delta f_k\|^2+\alpha_k^2\|\Delta f_k\|^2\notag\\
&\leq (1-\frac{\mu_h\alpha_k}{4})z_k+\frac{9L^2\alpha_k}{\mu_h}y_k+\frac{9\alpha_k}{2\mu_h}\|\Delta f_k\|^2,
\end{align*}
where again we use $\alpha_k\leq\frac{1}{2\mu_h}$ to simplify the terms.

\qed

\subsection{Proof of Lemma~\ref{lem:y_k:stronglyconvex}}

By the update rule in \eqref{alg:update_decision},
\begin{align*}
\omega_{k+1}-\omega^{\star}(\theta_{k+1})&=\omega_{k+1}-\omega^{\star}(\theta_{k+1})\notag\\
&=\left(\omega_{k}-\omega^{\star}(\theta_{k})\right)-\beta_k g_k - \left(\omega^{\star}(\theta_{k+1})-\omega^{\star}(\theta_k)\right)\notag\\
&=\left(\omega_{k}-\omega^{\star}(\theta_{k})\right)-\beta_k G(\theta_k,\omega_k)-\beta_k \Delta g_k - \left(\omega^{\star}(\theta_{k+1})-\omega^{\star}(\theta_k)\right).
\end{align*}
Taking the norm yields
\begin{align}
y_{k+1}&=\|\omega_{k+1}-\omega^{\star}(\theta_{k+1})\|^2\notag\\
&= \underbrace{\|\omega_{k}-\omega^{\star}(\theta_{k})-\beta_k G(\theta_k,\omega_k)\|^2}_{T_1}+\beta_k^2\|\Delta g_k\|^2+\underbrace{\|\omega^{\star}(\theta_{k+1})-\omega^{\star}(\theta_k)\|^2}_{T_2}\notag\\
&\hspace{20pt} \underbrace{-2\left(\omega_{k}-\omega^{\star}(\theta_{k})-\beta_k G(\theta_k,\omega_k)\right)^{\top}\left(\omega^{\star}(\theta_{k+1})-\omega^{\star}(\theta_k)\right)}_{T_3}\notag\\
&\hspace{20pt}\underbrace{-2\beta_k\left(\omega_{k}-\omega^{\star}(\theta_{k})-\beta_k G(\theta_k,\omega_k)\right)^{\top}\Delta g_k}_{T_4} + \underbrace{2\beta_k \Delta g_k^{\top}\left(\omega^{\star}(\theta_{k+1})-\omega^{\star}(\theta_k)\right)}_{T_5}.\label{lem:y_k:proof_eq1:stronglyconvex}
\end{align}

We bound each term of \eqref{lem:y_k:proof_eq1:stronglyconvex} individually. First, since by definition $G(\theta,\omega^{\star}(\theta))=0$ for all $\theta$, we have
\begin{align*}
T_1&=\|\omega_{k}-\omega^{\star}(\theta_{k})\|^2 - 2\beta_k\left(\omega_{k}-\omega^{\star}(\theta_{k})\right)^{\top}G(\theta_k,\omega_k)+\beta_k^2\|G(\theta_k,\omega_k)\|^2\notag\\
&=y_k - 2\beta_k\left(\omega_{k}-\omega^{\star}(\theta_{k})\right)^{\top}G(\theta_k,\omega_k)+\beta_k^2\|G(\theta_k,\omega_k)-G(\theta_k,\omega^{\star}(\theta_k))\|^2\notag\\
&\leq y_k - 2\mu_G\beta_k\|\omega_k-\omega^{\star}(\theta_k)\|^2+L^2\beta_k^2\|\omega_k-\omega^{\star}(\theta_k)\|^2\notag\\
&=(1-2\mu_G\beta_k+L^2\beta_k^2)y_k,
\end{align*}
where the first inequality follows from Assumption~\ref{assump:stronglymonotone_G} and the Lipschitz continuity of $G$.

The term $T_2$ can be simply treated using the Lipschitz condition of $\omega^{\star}$
\begin{align*}
T_2 
&=\|\omega^{\star}(\theta_{k+1})-\omega^{\star}(\theta_k)\|^2
\leq L^2\|\theta_{k+1}-\theta_k\|^2\leq L^2\alpha_k^2\left(\|\Delta f_k\| + L\sqrt{y_k} + L(L+1) \sqrt{z_k}\right)^2\notag\\
&\leq 3L^2\alpha_k^2\left(\|\Delta f_k\|^2 + L^2 y_k + (L+1)^4 z_k\right).
\end{align*}

By the Cauchy-Schwarz inequality, 
\begin{align*}
&T_3\notag\\
&\leq 2\|\omega_{k}-\omega^{\star}(\theta_{k})-\beta_k G(\theta_k,\omega_k)\|\|\omega^{\star}(\theta_{k+1})-\omega^{\star}(\theta_k)\|\notag\\
&\leq 2\|\omega_{k}-\omega^{\star}(\theta_{k})\|\|\omega^{\star}(\theta_{k+1})-\omega^{\star}(\theta_k)\|+2\beta_k\|G(\theta_k,\omega_k)-G(\theta_k, \omega^{\star}(\theta_k))\|\|\omega^{\star}(\theta_{k+1})\hspace{-2pt}-\hspace{-2pt}\omega^{\star}(\theta_k)\|\notag\\
&\leq (2L+2L^2\beta_k)\sqrt{y_k}\|\theta_{k+1}-\theta_k\|\notag\\
&\leq 4L\alpha_k\sqrt{y_k}\left(\|\Delta f_k\| + L\sqrt{y_k} + L(L+1) \sqrt{z_k}\right).
\end{align*}
where the second inequality follows from $G(\theta,\omega^{\star}(\theta))=0$ for all $\theta$, and the fourth inequality applies Lemma~\ref{lem:Lipschitz:stronglyconvex} and the step size rule $\beta_k\leq1/L$.

We can bound $T_4$ following a similar line of analysis
\begin{align*}
T_4 &\leq 2\beta_k\|\omega_{k}-\omega^{\star}(\theta_{k})-\beta_k G(\theta_k,\omega_k)\|\|\Delta g_k\|\leq 4\beta_k\sqrt{y_k}\|\Delta g_k\|.
\end{align*}

Finally, we treat $T_5$ again by the Cauchy-Schwarz inequality
\begin{align*}
T_5&\leq2\beta_k\|\Delta g_k\|\|\omega^{\star}(\theta_{k+1})-\omega^{\star}(\theta_k)\|\notag\\
&\leq 2L\alpha_k\beta_k\|\Delta g_k\| \Big(\|\Delta f_k\| + L\sqrt{y_k} + L(L+1)\sqrt{z_k}\Big),
\end{align*}
where the second inequality again uses Lemma~\ref{lem:Lipschitz:stronglyconvex}.

Applying the bounds on $T_1$-$T_5$ to \eqref{lem:y_k:proof_eq1:stronglyconvex}, we get
\begin{align*}
y_{k+1}& \leq (1-2\mu_G\beta_k+L^2\beta_k^2)y_k+\beta_k^2\|\Delta g_k\|^2+3L^2\alpha_k^2\left(\|\Delta f_k\|^2 + L^2 y_k + (L+1)^4 z_k\right)\notag\\
&\hspace{20pt}+4L\alpha_k\sqrt{y_k}\left(\|\Delta f_k\| + L\sqrt{y_k} + L(L+1)\sqrt{z_k}\right) + 4\beta_k\sqrt{y_k}\|\Delta g_k\|\notag\\
&\hspace{20pt}+2L\alpha_k\beta_k\|\Delta g_k\| \Big(\|\Delta f_k\| + L\sqrt{y_k} + L(L+1)\sqrt{z_k}\Big)\\
&\leq (1-2\mu_G\beta_k+L^2\beta_k^2)y_k+\beta_k^2\|\Delta g_k\|^2+3L^2\alpha_k^2\left(\|\Delta f_k\|^2 + L^2 y_k + (L+1)^4 z_k\right)\notag\\
&\hspace{20pt}+2L\alpha_k y_k+2L\alpha_k\|\Delta f_k\|^2+4L^2\alpha_k y_k+\frac{\mu_G}{4}\beta_k y_k+\frac{16(L+1)^6\alpha_k^2}{\mu_G\beta_k}z_k+ \frac{\mu_G\beta_k}{2} y_k\notag\\
&\hspace{20pt}+\frac{8\beta_k}{\mu_G}\|\Delta g_k\|^2+L\alpha_k\beta_k\|\Delta g_k\|^2+3L\alpha_k\beta_k\left(\|\Delta f_k\|^2 + L^2 y_k + (L+1)^4 z_k\right)\\
&\leq (1-\mu_G\beta_k)y_k-\left(\frac{\mu_G}{4}\beta_k -8L^4\alpha_k-L^2\beta_k^2\right)y_k\notag\\
&\hspace{20pt}+\left(\frac{19(L+1)^6\alpha_k^2}{\mu_G\beta_k}+3(L+1)^6\alpha_k\beta_k\right)z_k+\left(\frac{8\beta_k}{\mu_G}+2L\beta_k^2\right)\|\Delta g_k\|^2+8L\alpha_k\|\Delta f_k\|^2,
\end{align*}
where the second inequality follows from the fact that $2\|\va\|\|\vb\|\leq c\|\va\|^2+\frac{1}{c}\|\vb\|^2$ for any vector $\va,\vb$ and scalar $c>0$. We have also simplified terms using the conditions $\alpha_k\leq\beta_k\leq\frac{1}{72(L+1)^6}$ and $\beta_k\leq1/\mu_G$.

\qed
\section{Analysis for Functions Satisfying Polyak-Lojasiewicz Condition}\label{sec:proof_PL}

The exact requirement on step sizes is
\begin{align}
\lambda_{k} = \frac{c_{\lambda}}{k+\tau+1},\quad \alpha_{k} = \frac{c_\alpha}{k+\tau+1},\quad \beta_{k} = \frac{c_\beta}{k+\tau+1},\label{eq:step_sizes:PL}
\end{align}
where $c_{\alpha},c_\beta,c_{\lambda}$ need to be carefully balanced and $\tau$ selected sufficiently large such that $\alpha_k\leq\beta_k\leq\lambda_k\leq1/4$, $c_\beta\leq1/L$, $c_\beta\leq \frac{\mu_G}{320L^4}c_\lambda$, $c_\beta\leq\frac{\mu_h^3}{480(L+1)^6}$, $c_\beta\leq\frac{\mu_G}{8L^2}$, $\beta_k\leq\frac{\lambda_k}{4}\min\{(\frac{2L^2}{\mu_h^3}+58L^2)^{-1},(50L^2+\frac{8}{\mu_G})^{-1}\}$, $\alpha_k\leq\frac{\mu_G\beta_k}{4\left(14L^3+\frac{132(L+1)^6}{\mu_h^3}\right)}$, $\alpha_k\leq \frac{\mu_h^2}{3072(L+1)^6}\lambda_k$, $\alpha_k\leq\beta_k$, and $c_{\alpha}\geq\frac{2}{\min\{\mu_h,\mu_G\}}$.

The PL condition implies the following quadratic growth inequality, where $\text{Proj}_{\Theta^{\star}}(\theta)$ denotes the projection of $\theta$ to $\Theta^{\star}$ \citep{karimi2016linear}.
\begin{align}
\frac{2}{\mu_h}\left(h(\theta)-h^{\star}\right)\geq\|\text{Proj}_{\Theta^{\star}}(\theta)-\theta\|^2.\label{eq:quadratic_growth}
\end{align}
Equation~\eqref{eq:quadratic_growth} plays an important role in the proof of a lemma that we now introduce.

\begin{lem}\label{lem:Lipschitz:PL}
Under Assumption~\ref{assump:Lipschitz}-\ref{assump:stronglymonotone_G} and Assumption~\ref{assump:PL_condition}, the sequence of variables $\{\theta_k,\omega_k,f_k,g_k\}$ satisfy for all $k$
\begin{align*}
    \|f_k\|&\leq \|\Delta f_k\| + L\sqrt{y_k} + \frac{2L(L+1)}{\mu_h}\sqrt{x_k},\\
    \|g_k\|&\leq \|\Delta g_k\|+L\sqrt{y_k},\\
    \|\theta_{k+1}-\theta_k\|&\leq\alpha_k\left(\|\Delta f_k\| + L\sqrt{y_k} + \frac{2L(L+1)}{\mu_h}\sqrt{x_k}\right),\\
    \|\omega_{k+1}-\omega_k\| & \leq\beta_k\left(\|\Delta g_k\|+L\sqrt{y_k}\right).
\end{align*}
\end{lem}

Our analysis in this section also relies on a few more lemmas which we present below. The proofs of all lemmas can be found at the end of the section.

\begin{lem}\label{lem:Delta_f:PL}
Under Assumption~\ref{assump:Lipschitz}-\ref{assump:stronglymonotone_G} and Assumption~\ref{assump:PL_condition} and the step sizes given in \eqref{eq:step_sizes:PL}, we have the following bound on $\|\Delta f_{k+1}\|^2$
\begin{align*}
\mathbb{E}[\|\Delta f_{k+1}\|^2]&\leq(1-\lambda_k)\mathbb{E}[\|\Delta f_k\|^2]-(\frac{\lambda_k}{4}-\frac{24L^2\beta_k^2}{\lambda_k})\mathbb{E}[\|\Delta f_k\|^2]\notag\\
&\hspace{20pt}+\frac{4}{\lambda_k}\mathbb{E}\left[6L^2\beta_k^2\|\Delta g_k\|^2+10L^4\beta_k^2 y_k + \frac{24(L+1)^6}{\mu_h}\alpha_k^2 x_k\right]+2B\lambda_k^2.
\end{align*}
\end{lem}

\begin{lem}\label{lem:Delta_g:PL}
Under Assumption~\ref{assump:Lipschitz}-\ref{assump:stronglymonotone_G} and Assumption~\ref{assump:PL_condition} and the step sizes given in \eqref{eq:step_sizes:PL}, we have the following bound on $\|\Delta g_{k+1}\|^2$
\begin{align*}
\mathbb{E}[\|\Delta g_{k+1}\|^2]
&\leq (1-\lambda_k)\mathbb{E}[\|\Delta g_k\|^2]-(\frac{\lambda_k}{4}-\frac{24L^2\beta_k^2}{\lambda_k})\mathbb{E}[\|\Delta g_k\|^2]\notag\\
&\hspace{20pt}+\frac{4}{\lambda_k}\mathbb{E}\left[6L^2\beta_k^2\|\Delta f_k\|^2+10L^4\beta_k^2 y_k + \frac{24(L+1)^6}{\mu_h}\alpha_k^2 x_k\right]+2B\lambda_k^2.
\end{align*}
\end{lem}

\begin{lem}\label{lem:x_k:PL}
Under Assumption~\ref{assump:Lipschitz}-\ref{assump:stronglymonotone_G} and Assumption~\ref{assump:PL_condition} and the step sizes given in \eqref{eq:step_sizes:PL}, we have
\begin{align*}
x_{k+1} &\leq (1-\mu_h\alpha_k)x_k-(\frac{\mu_h}{4}\alpha_k-\frac{6(L+1)^5}{\mu_h^2}\alpha_k^2) x_k \notag\\
&\hspace{20pt}+ (\frac{4L^4}{\mu_h^3}+2L)\alpha_k y_k +(\frac{2L^2\alpha_k}{\mu_h^3}+2L\alpha_k^2)\|\Delta f_k\|^2.
\end{align*}
\end{lem}

\begin{lem}\label{lem:y_k:PL}
Under Assumption~\ref{assump:Lipschitz}-\ref{assump:stronglymonotone_G} and Assumption~\ref{assump:PL_condition} and the step sizes given in \eqref{eq:step_sizes:PL}, we have
\begin{align*}
y_{k+1}
&\leq (1-\mu_G\beta_k)y_k-\left(\frac{\mu_G}{2}\beta_k -\left(12L^3+\frac{128(L+1)^6}{\mu_h^3}\right)\alpha_k-L^2\beta_k^2\right)y_k\notag\\
&\hspace{20pt}+\left(\frac{\mu_h}{8}\alpha_k+\frac{24(L+1)^6}{\mu_h^2}\alpha_k\beta_k\right)x_k+\left(\frac{8\beta_k}{\mu_G}+2L\beta_k^2\right)\|\Delta g_k\|^2+8L\alpha_k\|\Delta f_k\|^2.
\end{align*}
\end{lem}

\subsection{Proof of Theorem~\ref{thm:PL}}

Combining the results from Lemma~\ref{lem:Delta_f:PL}-\ref{lem:y_k:PL} and simplifying the terms with the choice of step sizes given in \eqref{eq:step_sizes:PL},
\begin{align*}
&\mathbb{E}[V_{k+1}]\notag\\
&= \mathbb{E}[\|\Delta f_{k+1}\|^2 + \|\Delta g_{k+1}\|^2 + x_{k+1} + y_{k+1}]\notag\\
&\leq (1-\lambda_k)\mathbb{E}[\|\Delta f_k\|^2]-(\frac{\lambda_k}{4}-\frac{24L^2\beta_k^2}{\lambda_k})\mathbb{E}[\|\Delta f_k\|^2]\notag\\
&\hspace{20pt}+\frac{4}{\lambda_k}\mathbb{E}\left[6L^2\beta_k^2\|\Delta g_k\|^2+10L^4\beta_k^2 y_k + \frac{24(L+1)^6}{\mu_h}\alpha_k^2 x_k\right]+2B\lambda_k^2\notag\\
&\hspace{20pt}+(1-\lambda_k)\mathbb{E}[\|\Delta g_k\|^2]-(\frac{\lambda_k}{4}-\frac{24L^2\beta_k^2}{\lambda_k})\mathbb{E}[\|\Delta g_k\|^2]\notag\\
&\hspace{20pt}+\frac{4}{\lambda_k}\mathbb{E}\left[6L^2\beta_k^2\|\Delta f_k\|^2+10L^4\beta_k^2 y_k + \frac{24(L+1)^6}{\mu_h}\alpha_k^2 x_k\right]+2B\lambda_k^2\notag\\
&\hspace{20pt}+\mathbb{E}\Big[(1\hspace{-2pt}-\hspace{-2pt}\mu_h\alpha_k)x_k\hspace{-2pt}-\hspace{-2pt}(\frac{\mu_h}{4}\alpha_k\hspace{-2pt}-\hspace{-2pt}\frac{6(L\hspace{-2pt}+\hspace{-2pt}1)^5}{\mu_h^2}\alpha_k^2) x_k \hspace{-2pt}+\hspace{-2pt} (\frac{4L^4}{\mu_h^3}\hspace{-2pt}+\hspace{-2pt}2L)\alpha_k y_k\hspace{-2pt} +\hspace{-2pt}(\frac{2L^2\alpha_k}{\mu_h^3}\hspace{-2pt}+\hspace{-2pt}2L\alpha_k^2)\|\Delta f_k\|^2\notag\\
&\hspace{20pt}+(1-\mu_G\beta_k)y_k-\left(\frac{\mu_G}{2}\beta_k -\left(12L^3+\frac{128(L+1)^6}{\mu_h^3}\right)\alpha_k-L^2\beta_k^2\right)y_k\notag\\
&\hspace{20pt}+\left(\frac{\mu_h}{8}\alpha_k+\frac{24(L+1)^6}{\mu_h^2}\alpha_k\beta_k\right)x_k+\left(\frac{8\beta_k}{\mu_G}+2L\beta_k^2\right)\|\Delta g_k\|^2+8L\alpha_k\|\Delta f_k\|^2\Big]\notag\\
&\leq (1-\lambda_k)\mathbb{E}[\|\Delta f_k\|^2]-(\frac{\lambda_k}{4}-\frac{48L^2\beta_k^2}{\lambda_k}-(\frac{2L^2}{\mu_h^3}+10L)\alpha_k)\mathbb{E}[\|\Delta f_k\|^2]\notag\\
&\hspace{20pt}+(1-\lambda_k)\mathbb{E}[\|\Delta g_k\|^2]-(\frac{\lambda_k}{4}-(50L^2+\frac{8}{\mu_G})\beta_k)\mathbb{E}[\|\Delta g_k\|^2]\notag\\
&\hspace{20pt}+\frac{8}{\lambda_k}\mathbb{E}\left[10L^4\beta_k^2 y_k + \frac{24(L+1)^6}{\mu_h}\alpha_k^2 x_k\right]+4B\lambda_k^2\notag\\
&\hspace{20pt}+(1-\mu_h\alpha_k)\mathbb{E}[x_k]-(\frac{\mu_h}{4}\alpha_k-\frac{6(L+1)^5}{\mu_h^2}\alpha_k^2) \mathbb{E}[x_k] + (\frac{4L^4}{\mu_h^3}+2L)\alpha_k\mathbb{E}[y_k] \notag\\
&\hspace{20pt}+(1-\mu_G\beta_k)\mathbb{E}[y_k]-\left(\frac{\mu_G}{2}\beta_k -\left(12L^3+\frac{128(L+1)^6}{\mu_h^3}\right)\alpha_k-L^2\beta_k^2\right)\mathbb{E}[y_k]\notag\\
&\hspace{20pt}+\left(\frac{\mu_h}{8}\alpha_k+\frac{24(L+1)^6}{\mu_h^2}\alpha_k\beta_k\right)\mathbb{E}[x_k]\notag\\
% &\leq (1-\lambda_k)\mathbb{E}[\|\Delta f_k\|^2]-(\frac{\lambda_k}{4}-(\frac{2L^2}{\mu_h^3}+58L^2)\beta_k)\mathbb{E}[\|\Delta f_k\|^2]\notag\\
% &\hspace{20pt}+(1-\lambda_k)\mathbb{E}[\|\Delta g_k\|^2]-(\frac{\lambda_k}{4}-(50L^2+\frac{8}{\mu_G})\beta_k)\mathbb{E}[\|\Delta g_k\|^2]\notag\\
% &\hspace{20pt}+(1-\mu_h\alpha_k)x_k-(\frac{\mu_h}{16}\alpha_k-\frac{30(L+1)^6}{\mu_h^2}\alpha_k\beta_k) x_k\notag\\
% &\hspace{20pt}+(1-\mu_G\beta_k)y_k-\left(\frac{\mu_G}{4}\beta_k -\left(14L^3+\frac{132(L+1)^6}{\mu_h^3}\right)\alpha_k-L^2\beta_k^2\right)y_k+4B\lambda_k^2
&\leq (1-\lambda_k)\mathbb{E}[\|\Delta f_k\|^2]\hspace{-2pt}+\hspace{-2pt}(1-\lambda_k)\mathbb{E}[\|\Delta g_k\|^2]+(1-\mu_h\alpha_k)\mathbb{E}[x_k]\hspace{-2pt}+\hspace{-2pt}(1-\mu_G\beta_k)\mathbb{E}[y_k]\hspace{-2pt}+\hspace{-2pt}4B\lambda_k^2\notag\\
&\hspace{20pt}-(\frac{\lambda_k}{4}-(\frac{2L^2}{\mu_h^3}+58L^2)\beta_k)\mathbb{E}[\|\Delta f_k\|^2]-(\frac{\lambda_k}{4}-(50L^2+\frac{8}{\mu_G})\beta_k)\mathbb{E}[\|\Delta g_k\|^2]\notag\\
&\hspace{20pt}-(\frac{\mu_h}{16}\alpha_k-\frac{30(L\hspace{-2pt}+\hspace{-2pt}1)^6}{\mu_h^2}\alpha_k\beta_k) \mathbb{E}[x_k]-\left(\frac{\mu_G}{4}\beta_k \hspace{-2pt}-\hspace{-2pt}\left(14L^3\hspace{-2pt}+\hspace{-2pt}\frac{132(L\hspace{-2pt}+\hspace{-2pt}1)^6}{\mu_h^3}\right)\alpha_k\hspace{-2pt}-\hspace{-2pt}L^2\beta_k^2\right)\mathbb{E}[y_k],
\end{align*}
% 
% $\frac{80L^4\beta_k^2}{\lambda_k}\leq\frac{\mu_G}{4}\beta_k$ by the step size rule $\beta_k\leq \frac{\mu_G}{320L^4}\lambda_k$
% 
% $\frac{192(L+1)^6}{\mu_h\lambda_k}\alpha_k^2\leq\frac{\mu_h}{16}\alpha_k$ guaranteed by the step size rule $\alpha_k\leq \frac{\mu_h^2}{3072(L+1)^6}\lambda_k$
%
where we use the step size conditions $\beta_k\leq \frac{\mu_G}{320L^4}\lambda_k$ and $\alpha_k\leq \frac{\mu_h^2}{3072(L+1)^6}\lambda_k$.

Again, due to carefully selected step sizes, (specifically, $\beta_k\leq\frac{\mu_h^3}{480(L+1)^6}$, $\beta_k\leq\frac{\mu_G}{8L^2}$, $\beta_k\leq\frac{\lambda_k}{4}\min\{(\frac{2L^2}{\mu_h^3}+58L^2)^{-1},(50L^2+\frac{8}{\mu_G})^{-1}\}$, and $\alpha_k\leq\frac{\mu_G\beta_k}{4\left(14L^3+\frac{132(L+1)^6}{\mu_h^3}\right)}$), we can make the last four terms of the inequality above non-positive. As a result, we have
\begin{align*}
&\mathbb{E}[V_{k+1}]\notag\\
&\leq(1-\lambda_k)\mathbb{E}[\|\Delta f_k\|^2]+(1-\lambda_k)\mathbb{E}[\|\Delta g_k\|^2]+(1-\mu_h\alpha_k)\mathbb{E}[x_k]+(1-\mu_G\beta_k)\mathbb{E}[y_k]+4B\lambda_k^2\notag\\
&\leq (1-\min\{\mu_h,\mu_G\}\alpha_k)\mathbb{E}[V_k]+4B\lambda_k^2\notag\\
&= (1-\frac{\min\{\mu_h,\mu_G\}c_{\alpha}}{k+\tau+1})\mathbb{E}[V_k]+\frac{B}{4(k+1)^2}\notag\\
&\leq (1-\frac{2}{k+\tau+1})\mathbb{E}[V_k]+\frac{B}{4(k+1)^2}\notag\\
&\leq \frac{k+\tau-1}{k+\tau+1}\mathbb{E}[V_k]+\frac{B}{4(k+1)^2}.
\end{align*}

Multiplying by $(k+\tau+1)^2$ on both sides of the inequality, we have
\begin{align*}
(k+\tau+1)^2\mathbb{E}[V_{k+1}]&\leq (k+\tau+1)(k+\tau-1)\mathbb{E}[V_{k}]+\frac{B(k+\tau+1)^2}{4(k+1)^2}\notag\\
&\leq (k+\tau)^2\mathbb{E}[V_{k}]+\left(\frac{B}{2}+\frac{B\tau^2}{2}\right)\notag\\
&\leq \tau^2 V_0+\left(\frac{B}{2}+\frac{B\tau^2}{2}\right)(k+1).
\end{align*}

As a result,
\begin{align*}
\mathbb{E}[V_{k+1}]&\leq \frac{B\tau^2(k+1)}{(k+\tau+1)^2}+\frac{\tau^2V_0}{(k+\tau+1)^2}\notag\\
&\leq \frac{B\tau^2}{k+\tau+1}+\frac{\tau^2 V_0}{(k+\tau+1)^2},
\end{align*}
which obviously leads to the claimed result.

\qed

\subsection{Proof of Lemma~\ref{lem:Lipschitz:PL}}

By definition,
\begin{align*}
    F(\theta, \omega^{\star}(\theta))=\nabla h(\theta)=0,\quad\forall \theta\in\Theta^{\star}
\end{align*}
which implies
\begin{align*}
    &\|f_k\|\notag\\
    &=\|\Delta f_k+F(\theta_k,\omega_k)-F(\theta_k,\omega^{\star}(\theta_k))+F(\theta_k,\omega^{\star}(\theta_k))-F(\text{Proj}_{\Theta^{\star}}(\theta_k),\omega^{\star}(\text{Proj}_{\Theta^{\star}}(\theta_k)))\|\\
    &\leq\|\Delta f_k\|+\|F(\theta_k,\omega_k)-F(\theta_k,\omega^{\star}(\theta_k))\|+\|F(\theta_k,\omega^{\star}(\theta_k))-F(\text{Proj}_{\Theta^{\star}}(\theta_k),\omega^{\star}(\text{Proj}_{\Theta^{\star}}(\theta_k)))\|\\
    &\leq \|\Delta f_k\| + L\|\omega_k-\omega^{\star}(\theta_k)\| + L(L+1)\|\theta_k-\text{Proj}_{\Theta^{\star}}(\theta_k)\|\notag\\
    &\leq \|\Delta f_k\| + L\sqrt{y_k} + \frac{2L(L+1)}{\mu_h}\sqrt{x_k},
\end{align*}
where the second inequality follows from the Lipschitz continuity of $F$ and $\omega^{\star}$ and the final inequality applies \eqref{eq:quadratic_growth}.

Similarly, since $G(\theta,\omega^{\star}(\theta))=0$ for any $\theta$, we can derive
\begin{align*}
    \|g_k\|&=\|\Delta g_k+G(\theta_k,\omega_k)-G(\theta_k,\omega^{\star}(\theta_k))\|\notag\\
    &\leq\|\Delta g_k\|+\|G(\theta_k,\omega_k)-G(\theta_k,\omega^{\star}(\theta_k))\|\notag\\
    &\leq \|\Delta g_k\|+L\sqrt{y_k}.
\end{align*}

The bounds on $\|\theta_{k+1}-\theta_k\|$ and $\|\omega_{k+1}-\omega_k\|$ easily follow from the two inequalities above and \eqref{eq:update_auxiliary}.

\qed

\subsection{Proof of Lemma~\ref{lem:Delta_f:PL}}

By the update rule in \eqref{eq:update_auxiliary}, we have
\begin{align*}
\Delta f_{k+1}&=f_{k+1}-F(\theta_{k+1},\omega_{k+1})\notag\\
&=(1-\lambda_k)f_k+\lambda_k F(\theta_k,\omega_k,X_k)-F(\theta_{k+1},\omega_{k+1})\notag\\
&=(1-\lambda_k)f_k+\lambda_k F(\theta_k,\omega_k)-F(\theta_{k+1},\omega_{k+1})+\lambda_k\left(F(\theta_k,\omega_k)-F(\theta_k,\omega_k,X_k)\right)\notag\\
&=(1-\lambda_k)\Delta f_k + F(\theta_k,\omega_k)-F(\theta_{k+1},\omega_{k+1})+\lambda_k\left(F(\theta_k,\omega_k)-F(\theta_k,\omega_k,X_k)\right).
\end{align*}
This leads to 
\begin{align}
&\|\Delta f_{k+1}\|^2\notag\\
&=(1-\lambda_k)^2\|\Delta f_k\|^2 +\Big\|\left(F(\theta_k,\omega_k)-F(\theta_{k+1},\omega_{k+1})\right)+\lambda_k\left(F(\theta_k,\omega_k)-F(\theta_k,\omega_k,X_k)\right)\Big\|^2\notag\\
&\hspace{20pt}+2(1-\lambda_k)\Delta f_k^{\top}\left(F(\theta_k,\omega_k)-F(\theta_{k+1},\omega_{k+1})\right)\notag\\
&\hspace{20pt}+2(1-\lambda_k)\lambda_k\Delta f_k^{\top}\left(F(\theta_k,\omega_k)-F(\theta_k,\omega_k,X_k)\right)\notag\\
&\leq (1-\lambda_k)^2\|\Delta f_k\|^2 +2\|F(\theta_k,\omega_k)-F(\theta_{k+1},\omega_{k+1})\|^2+2\lambda_k^2\|F(\theta_k,\omega_k)-F(\theta_k,\omega_k,X_k)\|^2\notag\\
&\hspace{20pt}+\frac{\lambda_k}{2}\|\Delta f_k\|^2+\frac{2}{\lambda_k}\left\|F(\theta_k,\omega_k)-F(\theta_{k+1},\omega_{k+1})\right\|^2\notag\\
&\hspace{20pt}+2(1-\lambda_k)\lambda_k\Delta f_k^{\top}\left(F(\theta_k,\omega_k)-F(\theta_k,\omega_k,X_k)\right)\notag\\
&\leq (1-\lambda_k)\|\Delta f_k\|^2 +2\lambda_k^2\|F(\theta_k,\omega_k)-F(\theta_k,\omega_k,X_k)\|^2\notag\\
&\hspace{20pt}+(\lambda_k^2-\frac{\lambda_k}{2})\|\Delta f_k\|^2+\frac{4}{\lambda_k}\left\|F(\theta_k,\omega_k)-F(\theta_{k+1},\omega_{k+1})\right\|^2\notag\\
&\hspace{20pt}+2(1-\lambda_k)\lambda_k\Delta f_k^{\top}\left(F(\theta_k,\omega_k)-F(\theta_k,\omega_k,X_k)\right),\label{lem:Delta_f:proof_eq1:PL}
\end{align}
where the last inequality follows from $\lambda_k\leq 1$.

The second term on the right hand side of \eqref{lem:Delta_f:proof_eq1:PL} is bounded in expectation, i.e.
\[\mathbb{E}_{X_k\sim\xi}[\|F(\theta_k,\omega_k,X_k)-F(\theta_k,\omega_k)\|^2]\leq B.\]
% since $\mathbb{E}_{X_k\sim\xi}[\|F(\theta_k,\omega_k,X_k)-F(\theta_k,\omega_k)\|^2]\leq B$.
In addition, the last term of \eqref{lem:Delta_f:proof_eq1:PL} is zero in expectation since $\mathbb{E}[\Delta f_k^{\top}(F(\theta_k,\omega_k)-F(\theta_k,\omega_k,X_k))\mid \Hcal_{k-1}]=\Delta f_k^{\top}\mathbb{E}[(F(\theta_k,\omega_k)-F(\theta_k,\omega_k,X_k))\mid \Hcal_{k-1}]=0$. As a result, we have
\begin{align}
\mathbb{E}[\|\Delta f_{k+1}\|^2]&\leq (1-\lambda_k)\mathbb{E}[\|\Delta f_k\|^2]-(\frac{\lambda_k}{2}-\lambda_k^2)\mathbb{E}[\|\Delta f_k\|^2]\notag\\
&\hspace{20pt}+\frac{4}{\lambda_k}\mathbb{E}[\left\|F(\theta_k,\omega_k)-F(\theta_{k+1},\omega_{k+1})\right\|^2]+2B\lambda_k^2.\label{lem:Delta_f:proof_eq2:PL}
\end{align}

By Lemma~\ref{lem:Lipschitz:PL} and the Lipschitz condition of $F$,
\begin{align*}
    &\|F(\theta_k,\omega_k)-F(\theta_{k+1},\omega_{k+1})\|^2\notag\\
    & \leq 2L^2\|\theta_{k+1}-\theta_k\|^2+ 2L^2\|\omega_{k+1}-\omega_k\|^2\notag\\
    &\leq 6L^2\alpha_k^2(\|\Delta f_k\|^2+L^2 y_k+\frac{4(L+1)^4}{\mu_h^2}x_k)+4L^2\beta_k^2(\|\Delta g_k\|^2+L^2 y_k)\notag\\
    &\leq 6L^2\beta_k^2\|\Delta f_k\|^2 + 6L^2\beta_k^2\|\Delta g_k\|^2+10L^4\beta_k^2 y_k + \frac{24(L+1)^6}{\mu_h}\alpha_k^2 x_k.
\end{align*}
Plugging this bound into \eqref{lem:Delta_f:proof_eq2:PL}, we get
\begin{align*}
\mathbb{E}[\|\Delta f_{k+1}\|^2]&\leq (1-\lambda_k)\mathbb{E}[\|\Delta f_k\|^2]-(\frac{\lambda_k}{2}-\lambda_k^2)\mathbb{E}[\|\Delta f_k\|^2]+2B\lambda_k^2\notag\\
&\hspace{20pt}+\frac{4}{\lambda_k}\mathbb{E}\left[6L^2\beta_k^2\|\Delta f_k\|^2 + 6L^2\beta_k^2\|\Delta g_k\|^2+10L^4\beta_k^2 y_k + \frac{24(L+1)^6}{\mu_h}\alpha_k^2 x_k\right]\notag\\
&\leq (1-\lambda_k)\mathbb{E}[\|\Delta f_k\|^2]-(\frac{\lambda_k}{4}-\frac{24L^2\beta_k^2}{\lambda_k})\mathbb{E}[\|\Delta f_k\|^2]\notag\\
&\hspace{20pt}+\frac{4}{\lambda_k}\mathbb{E}\left[6L^2\beta_k^2\|\Delta g_k\|^2+10L^4\beta_k^2 y_k + \frac{24(L+1)^6}{\mu_h}\alpha_k^2 x_k\right]+2B\lambda_k^2,
\end{align*}
where the last inequality follows from $\lambda_k \leq 1/4$.

\qed

\subsection{Proof of Lemma~\ref{lem:Delta_g:PL}}

By the update rule in \eqref{eq:update_auxiliary}, we have
\begin{align*}
\Delta g_{k+1}&=g_{k+1}-G(\theta_{k+1},\omega_{k+1})\notag\\
&=(1-\lambda_k)g_k+\lambda_k G(\theta_k,\omega_k,X_k)-G(\theta_{k+1},\omega_{k+1})\notag\\
&=(1-\lambda_k)g_k+\lambda_k G(\theta_k,\omega_k)-G(\theta_{k+1},\omega_{k+1})+\lambda_k\left(G(\theta_k,\omega_k)-G(\theta_k,\omega_k,X_k)\right)\notag\\
&=(1-\lambda_k)\Delta g_k + G(\theta_k,\omega_k)-G(\theta_{k+1},\omega_{k+1})+\lambda_k\left(G(\theta_k,\omega_k)-G(\theta_k,\omega_k,X_k)\right).
\end{align*}
This leads to 
\begin{align}
&\|\Delta g_{k+1}\|^2\notag\\
&=(1-\lambda_k)^2\|\Delta g_k\|^2 +\Big\|\left(G(\theta_k,\omega_k)-G(\theta_{k+1},\omega_{k+1})\right)+\lambda_k\left(G(\theta_k,\omega_k)-G(\theta_k,\omega_k,X_k)\right)\Big\|^2\notag\\
&\hspace{20pt}+2(1-\lambda_k)\Delta g_k^{\top}\left(G(\theta_k,\omega_k)-G(\theta_{k+1},\omega_{k+1})\right)\notag\\
&\hspace{20pt}+2(1-\lambda_k)\lambda_k\Delta g_k^{\top}\left(G(\theta_k,\omega_k)-G(\theta_k,\omega_k,X_k)\right)\notag\\
&\leq (1-\lambda_k)^2\|\Delta g_k\|^2 +2\|G(\theta_k,\omega_k)-G(\theta_{k+1},\omega_{k+1})\|^2+2\lambda_k^2\|G(\theta_k,\omega_k)-G(\theta_k,\omega_k,X_k)\|^2\notag\\
&\hspace{20pt}+\frac{\lambda_k}{2}\|\Delta g_k\|^2+\frac{2}{\lambda_k}\left\|G(\theta_k,\omega_k)-G(\theta_{k+1},\omega_{k+1})\right\|^2\notag\\
&\hspace{20pt}+2(1-\lambda_k)\lambda_k\Delta g_k^{\top}\left(G(\theta_k,\omega_k)-G(\theta_k,\omega_k,X_k)\right)\notag\\
&\leq (1-\lambda_k)\|\Delta g_k\|^2 +2\lambda_k^2\|G(\theta_k,\omega_k)-G(\theta_k,\omega_k,X_k)\|^2\notag\\
&\hspace{20pt}+(\lambda_k^2-\frac{\lambda_k}{2})\|\Delta g_k\|^2+\frac{4}{\lambda_k}\left\|G(\theta_k,\omega_k)-G(\theta_{k+1},\omega_{k+1})\right\|^2\notag\\
&\hspace{20pt}+2(1-\lambda_k)\lambda_k\Delta g_k^{\top}\left(G(\theta_k,\omega_k)-G(\theta_k,\omega_k,X_k)\right).\label{lem:Delta_g:proof_eq1:PL}
\end{align}

The second term on the right hand side of \eqref{lem:Delta_g:proof_eq1:PL} is bounded in expectation, i.e.
\[\mathbb{E}_{X_k\sim\xi}[\|G(\theta_k,\omega_k,X_k)-G(\theta_k,\omega_k)\|^2]\leq B.\]
% since $\mathbb{E}_{X_k\sim\xi}[\|G(\theta_k,\omega_k,X_k)-G(\theta_k,\omega_k)\|^2]\leq B$.
In addition, the last term of \eqref{lem:Delta_f:proof_eq1:PL} is zero in expectation since $\mathbb{E}[\Delta g_k^{\top}(G(\theta_k,\omega_k)-G(\theta_k,\omega_k,X_k))\mid \Hcal_{k-1}]=\Delta g_k^{\top}\mathbb{E}[(G(\theta_k,\omega_k)-G(\theta_k,\omega_k,X_k))\mid \Hcal_{k-1}]=0$. As a result, we have
\begin{align}
\mathbb{E}[\|\Delta g_{k+1}\|^2]&\leq (1-\lambda_k)\mathbb{E}[\|\Delta g_k\|^2]-(\frac{\lambda_k}{2}-\lambda_k^2)\mathbb{E}[\|\Delta g_k\|^2]\notag\\
&\hspace{20pt}+\frac{4}{\lambda_k}\mathbb{E}[\left\|G(\theta_k,\omega_k)-G(\theta_{k+1},\omega_{k+1})\right\|^2]+2B\lambda_k^2.\label{lem:Delta_g:proof_eq2:PL}
\end{align}

By Lemma~\ref{lem:Lipschitz:PL} and the Lipschitz condition of $G$,
\begin{align*}
    &\|G(\theta_k,\omega_k)-G(\theta_{k+1},\omega_{k+1})\|^2\notag\\
    & \leq 2L^2\|\theta_{k+1}-\theta_k\|^2+ 2L^2\|\omega_{k+1}-\omega_k\|^2\notag\\
    &\leq 6L^2\alpha_k^2(\|\Delta f_k\|^2+L^2 y_k+\frac{4(L+1)^4}{\mu_h^2}x_k)+4L^2\beta_k^2(\|\Delta g_k\|^2+L^2 y_k)\notag\\
    &\leq 6L^2\beta_k^2\|\Delta f_k\|^2 + 6L^2\beta_k^2\|\Delta g_k\|^2+10L^4\beta_k^2 y_k + \frac{24(L+1)^6}{\mu_h}\alpha_k^2 x_k.
\end{align*}
Plugging this bound into \eqref{lem:Delta_g:proof_eq2:PL}, we get
\begin{align*}
&\mathbb{E}[\|\Delta g_{k+1}\|^2]\notag\\
&\leq (1-\lambda_k)\mathbb{E}[\|\Delta g_k\|^2]-(\frac{\lambda_k}{2}-\lambda_k^2)\mathbb{E}[\|\Delta g_k\|^2]\notag\\
&\hspace{20pt}+\frac{4}{\lambda_k}\mathbb{E}\left[6L^2\beta_k^2\|\Delta f_k\|^2 + 6L^2\beta_k^2\|\Delta g_k\|^2+10L^4\beta_k^2 y_k + \frac{24(L+1)^6}{\mu_h}\alpha_k^2 x_k\right]+2B\lambda_k^2\notag\\
&\leq (1-\lambda_k)\mathbb{E}[\|\Delta g_k\|^2]-(\frac{\lambda_k}{4}-\frac{24L^2\beta_k^2}{\lambda_k})\mathbb{E}[\|\Delta g_k\|^2]\notag\\
&\hspace{20pt}+\frac{4}{\lambda_k}\mathbb{E}\left[6L^2\beta_k^2\|\Delta f_k\|^2+10L^4\beta_k^2 y_k + \frac{24(L+1)^6}{\mu_h}\alpha_k^2 x_k\right]+2B\lambda_k^2,
\end{align*}
where the last inequality follows from $\lambda_k \leq 1/4$.

\qed

\subsection{Proof of Lemma~\ref{lem:x_k:PL}}

From the update rule \eqref{alg:update_decision} and the Lipschitz gradient condition in Assumption~\ref{assump:Lipschitz}, we have
\begin{align*}
    h(\theta_{k+1})&\leq h(\theta_{k})+\langle\nabla h(\theta_k),\theta_{k+1}-\theta_k\rangle+\frac{L}{2}\|\theta_{k+1}-\theta_k\|^2\notag\\
    &= h(\theta_{k})-\alpha_k\langle\nabla h(\theta_k),f_k\rangle+\frac{L \alpha_k^2}{2}\|f_k\|^2\notag\\
    &= h(\theta_{k})-\alpha_k\langle\nabla h(\theta_k),\Delta f_k\rangle-\alpha_k\langle\nabla h(\theta_k),F(\theta_k,\omega_k)\rangle+\frac{L \alpha_k^2}{2}\|f_k\|^2\notag\\
    &\leq h(\theta_{k})-\alpha_k\langle\nabla h(\theta_k),\Delta f_k\rangle-\alpha_k\langle\nabla h(\theta_k),F(\theta_k,\omega^{\star}(\theta_k))\rangle\notag\\
    &\hspace{20pt}+\alpha_k\langle\nabla h(\theta_k),F(\theta_k,\omega^{\star}(\theta_k))-F(\theta_k,\omega_k)\rangle\notag\\
    &\hspace{20pt}+\frac{3L \alpha_k^2}{2}\left(\|\Delta f_k\|^2 + L y_k + \frac{4(L+1)^4}{\mu_h^2}x_k\right)\notag\\
    &=h(\theta_k)-\alpha_k\|\nabla h(\theta_k)\|^2-\alpha_k\langle\nabla h(\theta_k),\Delta f_k\rangle\notag\\
    &\hspace{20pt}+\alpha_k\langle\nabla h(\theta_k),F(\theta_k,\omega^{\star}(\theta_k))-F(\theta_k,\omega_k)\rangle\notag\\
    &\hspace{20pt}+\frac{3L \alpha_k^2}{2}\left(\|\Delta f_k\|^2 + L y_k + \frac{4(L+1)^4}{\mu_h^2}x_k\right),
\end{align*}
where the second inequality applies Lemma~\ref{lem:Lipschitz:PL}.

As a result of Assumption~\ref{assump:PL_condition}, the inequality implies
\begin{align}
x_{x+1}&=h(\theta_{k+1})-h(\theta^{\star})\notag\\
&\leq h(\theta_k)-h(\theta^{\star})-\alpha_k\|\nabla h(\theta_k)\|^2-\alpha_k\langle\nabla h(\theta_k),\Delta f_k\rangle\notag\\
&\hspace{20pt}+\alpha_k\langle\nabla h(\theta_k),F(\theta_k,\omega^{\star}(\theta_k))-F(\theta_k,\omega_k)\rangle+\frac{3L \alpha_k^2}{2}\left(\|\Delta f_k\|^2 + L y_k + \frac{4(L+1)^4}{\mu_h^2}x_k\right)\notag\\
&\leq h(\theta_k)-h(\theta^{\star})-2\mu_h\alpha_k\left(h(\theta_k)-h(\theta^{\star})\right)-\alpha_k\langle\nabla h(\theta_k),\Delta f_k\rangle\notag\\
&\hspace{20pt}+\alpha_k\langle\nabla h(\theta_k),F(\theta_k,\omega^{\star}(\theta_k))-F(\theta_k,\omega_k)\rangle+\frac{3L \alpha_k^2}{2}\left(\|\Delta f_k\|^2 + L y_k + \frac{4(L+1)^4}{\mu_h^2}x_k\right)\notag\\
&\leq (1-2\mu_h\alpha_k)x_k-\alpha_k\langle\nabla h(\theta_k),\Delta f_k\rangle\notag\\
&\hspace{20pt}+\alpha_k\langle\nabla h(\theta_k),F(\theta_k,\omega^{\star}(\theta_k))-F(\theta_k,\omega_k)\rangle+\frac{3L \alpha_k^2}{2}\left(\|\Delta f_k\|^2 \hspace{-2pt} + \hspace{-2pt}L y_k \hspace{-2pt}+\hspace{-2pt} \frac{4(L+1)^4}{\mu_h^2}x_k\right).
\label{lem:x_k:proof_eq1:PL}
\end{align}

To bound the second term on the right hand side of \eqref{lem:x_k:proof_eq1:PL},
\begin{align}
-\alpha_k\langle\nabla h(\theta_k),\Delta f_k\rangle
&\leq \alpha_k\|\nabla h(\theta_k)\|\|\Delta f_k\|\notag\\
&=\alpha_k\|\nabla h(\theta_k)-\nabla h(\text{Proj}_{\Theta^{\star}}(\theta_k))\|\|\Delta f_k\|\notag\\
&\leq L\alpha_k\|\theta_k-\text{Proj}_{\Theta^{\star}}(\theta_k)\|\|\Delta f_k\|\notag\\
&\leq\frac{2L\alpha_k}{\mu_h}\sqrt{x_k}\|\Delta f_k\|\notag\\
&\leq \frac{\mu_h}{2}\alpha_k x_k+\frac{2L^2\alpha_k}{\mu_h^3}\|\Delta f_k\|^2,\label{lem:x_k:proof_eq2:PL}
\end{align}
where the third inequality plugs in the quadratic growth inequality \eqref{eq:quadratic_growth} and the final inequality follows from the fact that $2\|\va\|\|\vb\|\leq c\|\va\|^2+\frac{1}{c}\|\vb\|^2$ for any vector $\va,\vb$ and scalar $c>0$.

The third term on the right hand side of \eqref{lem:x_k:proof_eq1:PL} can be treated in a similar manner
\begin{align}
\alpha_k\langle\nabla h(\theta_k),F(\theta_k,\omega^{\star}(\theta_k))-F(\theta_k,\omega_k)\rangle &\leq \alpha_k\|\nabla h(\theta_k)\|\|F(\theta_k,\omega^{\star}(\theta_k))-F(\theta_k,\omega_k)\|\notag\\
&\leq\frac{2L\alpha_k}{\mu_h}\sqrt{x_k}\cdot\left(L\|\omega^{\star}(\theta_k)-\omega_k\|\right)\notag\\
&=\frac{2L^2\alpha_k}{\mu_h}\sqrt{x_k}\sqrt{y_k}\notag\\
&\leq\frac{\mu_h}{4}\alpha_k x_k+\frac{4L^4\alpha_k}{\mu_h^3}y_k.
\label{lem:x_k:proof_eq3:PL}
\end{align}

Plugging \eqref{lem:x_k:proof_eq2:PL} and \eqref{lem:x_k:proof_eq3:PL} into \eqref{lem:x_k:proof_eq1:PL}, we get
\begin{align*}
&x_{k+1}\notag\\
&\leq (1-2\mu_h\alpha_k)x_k-\alpha_k\langle\nabla h(\theta_k),\Delta f_k\rangle\notag\\
&\hspace{20pt}+\alpha_k\langle\nabla h(\theta_k),F(\theta_k,\omega^{\star}(\theta_k))-F(\theta_k,\omega_k)\rangle+\frac{3L \alpha_k^2}{2}\left(\|\Delta f_k\|^2 + L y_k + \frac{4(L+1)^4}{\mu_h^2}x_k\right)\notag\\
&\leq (1-2\mu_h\alpha_k)x_k+\frac{\mu_h}{2}\alpha_k x_k+\frac{2L^2\alpha_k}{\mu_h^3}\|\Delta f_k\|^2\notag\\
&\hspace{20pt}+\frac{\mu_h}{4}\alpha_k x_k+\frac{4L^4\alpha_k}{\mu_h^3}y_k+\frac{3L \alpha_k^2}{2}\left(\|\Delta f_k\|^2 + L y_k + \frac{4(L+1)^4}{\mu_h^2}x_k\right)\notag\\
&\leq (1-\mu_h\alpha_k)x_k-(\frac{\mu_h}{4}\alpha_k-\frac{6(L+1)^5}{\mu_h^2}\alpha_k^2) x_k + (\frac{4L^4}{\mu_h^3}+2L)\alpha_k y_k +(\frac{2L^2\alpha_k}{\mu_h^3}+2L\alpha_k^2)\|\Delta f_k\|^2.
\end{align*}

\qed

\subsection{Proof of Lemma~\ref{lem:y_k:PL}}

By the update rule in \eqref{alg:update_decision},
\begin{align*}
\omega_{k+1}-\omega^{\star}(\theta_{k+1})&=\omega_{k+1}-\omega^{\star}(\theta_{k+1})\notag\\
&=\left(\omega_{k}-\omega^{\star}(\theta_{k})\right)-\beta_k g_k - \left(\omega^{\star}(\theta_{k+1})-\omega^{\star}(\theta_k)\right)\notag\\
&=\left(\omega_{k}-\omega^{\star}(\theta_{k})\right)-\beta_k G(\theta_k,\omega_k)-\beta_k \Delta g_k - \left(\omega^{\star}(\theta_{k+1})-\omega^{\star}(\theta_k)\right).
\end{align*}
Taking the norm yields
\begin{align}
y_{k+1}&=\|\omega_{k+1}-\omega^{\star}(\theta_{k+1})\|^2\notag\\
&= \underbrace{\|\omega_{k}-\omega^{\star}(\theta_{k})-\beta_k G(\theta_k,\omega_k)\|^2}_{T_1}+\beta_k^2\|\Delta g_k\|^2+\underbrace{\|\omega^{\star}(\theta_{k+1})-\omega^{\star}(\theta_k)\|^2}_{T_2}\notag\\
&\hspace{20pt} \underbrace{-2\left(\omega_{k}-\omega^{\star}(\theta_{k})-\beta_k G(\theta_k,\omega_k)\right)^{\top}\left(\omega^{\star}(\theta_{k+1})-\omega^{\star}(\theta_k)\right)}_{T_3}\notag\\
&\hspace{20pt}\underbrace{-2\beta_k\left(\omega_{k}-\omega^{\star}(\theta_{k})-\beta_k G(\theta_k,\omega_k)\right)^{\top}\Delta g_k}_{T_4} + \underbrace{2\beta_k \Delta g_k^{\top}\left(\omega^{\star}(\theta_{k+1})-\omega^{\star}(\theta_k)\right)}_{T_5}.\label{lem:y_k:proof_eq1:PL}
\end{align}

We bound each term of \eqref{lem:y_k:proof_eq1:PL} individually. First, since by definition $G(\theta,\omega^{\star}(\theta))=0$ for all $\theta$, we have
\begin{align*}
T_1&=\|\omega_{k}-\omega^{\star}(\theta_{k})\|^2 - 2\beta_k\left(\omega_{k}-\omega^{\star}(\theta_{k})\right)^{\top}G(\theta_k,\omega_k)+\beta_k^2\|G(\theta_k,\omega_k)\|^2\notag\\
&=y_k - 2\beta_k\left(\omega_{k}-\omega^{\star}(\theta_{k})\right)^{\top}G(\theta_k,\omega_k)+\beta_k^2\|G(\theta_k,\omega_k)-G(\theta_k,\omega^{\star}(\theta_k))\|^2\notag\\
&\leq y_k - 2\mu_G\beta_k\|\omega_k-\omega^{\star}(\theta_k)\|^2+L^2\beta_k^2\|\omega_k-\omega^{\star}(\theta_k)\|^2\notag\\
&=(1-2\mu_G\beta_k+L^2\beta_k^2)y_k,
\end{align*}
where the first inequality follows from Assumption~\ref{assump:stronglymonotone_G} and the Lipschitz continuity of $G$.

The term $T_2$ can be simply treated using the Lipschitz condition of $\omega^{\star}$
\begin{align*}
T_2 
&=\|\omega^{\star}(\theta_{k+1})-\omega^{\star}(\theta_k)\|^2
\leq L^2\|\theta_{k+1}-\theta_k\|^2\leq L^2\alpha_k^2\left(\|\Delta f_k\| + L\sqrt{y_k} + \frac{2L(L+1)}{\mu_h}\sqrt{x_k}\right)^2\notag\\
&\leq 3L^2\alpha_k^2\left(\|\Delta f_k\|^2 + L^2 y_k + \frac{4(L+1)^4}{\mu_h^2}x_k\right).
\end{align*}

By the Cauchy-Schwarz inequality, 
\begin{align*}
&T_3\notag\\
&\leq 2\|\omega_{k}-\omega^{\star}(\theta_{k})-\beta_k G(\theta_k,\omega_k)\|\|\omega^{\star}(\theta_{k+1})-\omega^{\star}(\theta_k)\|\notag\\
&\leq 2\|\omega_{k}\hspace{-2pt}-\hspace{-2pt}\omega^{\star}(\theta_{k})\|\|\omega^{\star}(\theta_{k+1})\hspace{-2pt}-\hspace{-2pt}\omega^{\star}(\theta_k)\|+2\beta_k\|G(\theta_k,\omega_k)-G(\theta_k, \omega^{\star}(\theta_k))\|\|\omega^{\star}(\theta_{k+1})-\omega^{\star}(\theta_k)\|\notag\\
&\leq (2L+2L^2\beta_k)\sqrt{y_k}\|\theta_{k+1}-\theta_k\|\notag\\
&\leq 4L\alpha_k\sqrt{y_k}\left(\|\Delta f_k\| + L\sqrt{y_k} + \frac{2L(L+1)}{\mu_h}\sqrt{x_k}\right).
\end{align*}
where the second inequality follows from $G(\theta,\omega^{\star}(\theta))=0$ for all $\theta$, and the fourth inequality applies Lemma~\ref{lem:Lipschitz:PL} and the step size rule $\beta_k\leq1/L$.

We can bound $T_4$ following a similar line of analysis
\begin{align*}
T_4 &\leq 2\beta_k\|\omega_{k}-\omega^{\star}(\theta_{k})-\beta_k G(\theta_k,\omega_k)\|\|\Delta g_k\|\leq 4\beta_k\sqrt{y_k}\|\Delta g_k\|.
\end{align*}

Finally, we treat $T_5$ again by the Cauchy-Schwarz inequality
\begin{align*}
T_5&\leq2\beta_k\|\Delta g_k\|\|\omega^{\star}(\theta_{k+1})-\omega^{\star}(\theta_k)\|\notag\\
&\leq 2L\alpha_k\beta_k\|\Delta g_k\| \Big(\|\Delta f_k\| + L\sqrt{y_k} + \frac{2L(L+1)}{\mu_h}\sqrt{x_k}\Big),
\end{align*}
where the second inequality again uses Lemma~\ref{lem:Lipschitz:PL}.

Applying the bounds on $T_1$-$T_5$ to \eqref{lem:y_k:proof_eq1:PL}, we get
\begin{align*}
y_{k+1}& \leq (1-2\mu_G\beta_k+L^2\beta_k^2)y_k+\beta_k^2\|\Delta g_k\|^2+3L^2\alpha_k^2\left(\|\Delta f_k\|^2 + L^2 y_k + \frac{4(L+1)^4}{\mu_h^2}x_k\right)\notag\\
&\hspace{20pt}+4L\alpha_k\sqrt{y_k}\left(\|\Delta f_k\| + L\sqrt{y_k} + \frac{2L(L+1)}{\mu_h}\sqrt{x_k}\right) + 4\beta_k\sqrt{y_k}\|\Delta g_k\|\notag\\
&\hspace{20pt}+2L\alpha_k\beta_k\|\Delta g_k\| \Big(\|\Delta f_k\| + L\sqrt{y_k} + \frac{2L(L+1)}{\mu_h}\sqrt{x_k}\Big)\\
&\leq (1-2\mu_G\beta_k+L^2\beta_k^2)y_k+\beta_k^2\|\Delta g_k\|^2+3L^2\alpha_k^2\left(\|\Delta f_k\|^2 + L^2 y_k + \frac{4(L+1)^4}{\mu_h^2}x_k\right)\notag\\
&\hspace{20pt}+2L\alpha_k y_k+2L\alpha_k\|\Delta f_k\|^2+4L^2\alpha_k y_k+\frac{128(L+1)^6 \alpha_k}{\mu_h^3}y_k+\frac{\mu_h\alpha_k}{8}x_k + \frac{\mu_G\beta_k}{2} y_k\notag\\
&\hspace{20pt}+\frac{8\beta_k}{\mu_G}\|\Delta g_k\|^2 +L\alpha_k\beta_k\|\Delta g_k\|^2+3L\alpha_k\beta_k\left(\|\Delta f_k\|^2 + L^2 y_k + \frac{4(L+1)^4}{\mu_h^2}x_k\right)\\
&\leq (1-\mu_G\beta_k)y_k-\left(\frac{\mu_G}{2}\beta_k -\left(12L^3+\frac{128(L+1)^6}{\mu_h^3}\right)\alpha_k-L^2\beta_k^2\right)y_k\notag\\
&\hspace{20pt}+\left(\frac{\mu_h}{8}\alpha_k+\frac{24(L+1)^6}{\mu_h^2}\alpha_k\beta_k\right)x_k+\left(\frac{8\beta_k}{\mu_G}+2L\beta_k^2\right)\|\Delta g_k\|^2+8L\alpha_k\|\Delta f_k\|^2
\end{align*}
where the second inequality follows from the fact that $2\|\va\|\|\vb\|\leq c\|\va\|^2+\frac{1}{c}\|\vb\|^2$ for any vector $\va,\vb$ and scalar $c>0$.

\qed

\section{Analysis for Non-Convex Functions}\label{sec:proof_nonconvex}

Our analysis in this section relies on the following lemmas, the proofs of which can be found at the end of the section.

\begin{lem}\label{lem:Lipschitz:nonconvex}
Under Assumption~\ref{assump:Lipschitz}-\ref{assump:stronglymonotone_G}, the sequence of variables $\{\theta_k,\omega_k,f_k,g_k\}$ satisfy for all $k$
\begin{align*}
    \|f_k\|&\leq \|\Delta f_k\| + L\sqrt{y_k} + \|\nabla h(\theta_k)\|,\\
    \|g_k\|&\leq \|\Delta g_k\|+L\sqrt{y_k},\\
    \|\theta_{k+1}-\theta_k\|&\leq\alpha_k\left(\|\Delta f_k\| + L\sqrt{y_k} + \|\nabla h(\theta_k)\|\right),\\
    \|\omega_{k+1}-\omega_k\| & \leq\beta_k\left(\|\Delta g_k\|+L\sqrt{y_k}\right).
\end{align*}
\end{lem}

\begin{lem}\label{lem:Delta_f:nonconvex}
Under Assumption~\ref{assump:Lipschitz}-\ref{assump:stronglymonotone_G} and the step sizes given in \eqref{eq:step_sizes:nonconvex}, we have the following bound on $\|\Delta f_{k+1}\|^2$
\begin{align*}
\mathbb{E}[\|\Delta f_{k+1}\|^2]&\leq(1-\lambda_k)\mathbb{E}[\|\Delta f_k\|^2]-(\frac{\lambda_k}{4}-\frac{24L^2\beta_k^2}{\lambda_k})\mathbb{E}[\|\Delta f_k\|^2]\notag\\
&\hspace{20pt}+\frac{4}{\lambda_k}\mathbb{E}\left[6L^2\beta_k^2\|\Delta g_k\|^2+10L^4\beta_k^2 y_k + 6L^2\alpha_k^2 \|\nabla h(\theta_k)\|^2\right] +2B\lambda_k^2.
\end{align*}
\end{lem}

\begin{lem}\label{lem:Delta_g:nonconvex}
Under Assumption~\ref{assump:Lipschitz}-\ref{assump:stronglymonotone_G} and the step sizes given in \eqref{eq:step_sizes:nonconvex}, we have the following bound on $\|\Delta g_{k+1}\|^2$
\begin{align*}
\mathbb{E}[\|\Delta g_{k+1}\|^2]
&\leq (1-\lambda_k)\mathbb{E}[\|\Delta g_k\|^2]-(\frac{\lambda_k}{4}-\frac{24L^2\beta_k^2}{\lambda_k})\mathbb{E}[\|\Delta g_k\|^2]\notag\\
&\hspace{20pt}+\frac{4}{\lambda_k}\mathbb{E}\left[6L^2\beta_k^2\|\Delta f_k\|^2+10L^4\beta_k^2 y_k + 6L^2\alpha_k^2 \|\nabla h(\theta_k)\|^2\right]+2B\lambda_k^2.
\end{align*}
\end{lem}

\begin{lem}\label{lem:x_k:nonconvex}
Under Assumption~\ref{assump:Lipschitz}-\ref{assump:stronglymonotone_G} and the step sizes in \eqref{eq:step_sizes:nonconvex},
we have
\begin{align*}
h(\theta_{k+1})\leq h(\theta_k)-\frac{\alpha_k}{4}\|\nabla h(\theta_k)\|^2+\frac{5\alpha_k}{4}\|\Delta f_k\|^2+\frac{5L^2\alpha_k}{4} y_k.
\end{align*}
\end{lem}

\begin{lem}\label{lem:y_k:nonconvex}
Under Assumption~\ref{assump:Lipschitz}-\ref{assump:stronglymonotone_G} and the step sizes in \eqref{eq:step_sizes:nonconvex}, we have
\begin{align*}
y_{k+1}
&\leq (1-\mu_G\beta_k)y_k-\left(\frac{\mu_G}{2}\beta_k -108L^4\alpha_k-L^2\beta_k^2\right)y_k\notag\\
&\hspace{20pt}+\frac{\alpha_k}{8}\|\nabla h(\theta_k)\|^2+\left(\frac{8\beta_k}{\mu_G}+2L\beta_k^2\right)\|\Delta g_k\|^2+8L\alpha_k\|\Delta f_k\|^2.
\end{align*}
\end{lem}

\subsection{Proof of Theorem~\ref{thm:nonconvex}}

Combining the results from Lemma~\ref{lem:Delta_f:nonconvex}-\ref{lem:y_k:nonconvex} and simplifying the terms with the choice of step sizes given in \eqref{eq:step_sizes:nonconvex}, we have
\begin{align*}
&\mathbb{E}[\|\Delta f_{k+1}\|^2 + \|\Delta g_{k+1}\|^2 + h(\theta_{k+1}) + y_{k+1}]\notag\\
&\leq (1-\lambda_k)\mathbb{E}[\|\Delta f_k\|^2]-(\frac{\lambda_k}{4}-\frac{24L^2\beta_k^2}{\lambda_k})\mathbb{E}[\|\Delta f_k\|^2]\notag\\
&\hspace{20pt}+\frac{4}{\lambda_k}\mathbb{E}\left[6L^2\beta_k^2\|\Delta g_k\|^2+10L^4\beta_k^2 y_k + 6L^2\alpha_k^2 \|\nabla h(\theta_k)\|^2\right]+2B\lambda_k^2\notag\\
&\hspace{20pt}+(1-\lambda_k)\mathbb{E}[\|\Delta g_k\|^2]-(\frac{\lambda_k}{4}-\frac{24L^2\beta_k^2}{\lambda_k})\mathbb{E}[\|\Delta g_k\|^2]\notag\\
&\hspace{20pt}+\frac{4}{\lambda_k}\mathbb{E}\left[6L^2\beta_k^2\|\Delta f_k\|^2+10L^4\beta_k^2 y_k + 6L^2\alpha_k^2 \|\nabla h(\theta_k)\|^2\right]+2B\lambda_k^2\notag\\
&\hspace{20pt}+\mathbb{E}\Big[h(\theta_k)-\frac{\alpha_k}{4}\|\nabla h(\theta_k)\|^2+\frac{5\alpha_k}{4}\|\Delta f_k\|^2+\frac{5L^2\alpha_k}{4}y_k\notag\\
&\hspace{20pt}+(1-\mu_G\beta_k)y_k-\left(\frac{\mu_G}{2}\beta_k -108L^4\alpha_k-L^2\beta_k^2\right)y_k\notag\\
&\hspace{20pt}+\frac{\alpha_k}{8}\|\nabla h(\theta_k)\|^2+\left(\frac{8\beta_k}{\mu_G}+2L\beta_k^2\right)\|\Delta g_k\|^2+8L\alpha_k\|\Delta f_k\|^2\Big]\notag\\
&\leq (1-\lambda_k)\mathbb{E}[\|\Delta f_k\|^2]-(\frac{\lambda_k}{4}-\frac{48L^2\beta_k^2}{\lambda_k}-\frac{37L}{4}\alpha_k)\mathbb{E}[\|\Delta f_k\|^2]\notag\\
&\hspace{20pt}+(1-\lambda_k)\mathbb{E}[\|\Delta g_k\|^2]-(\frac{\lambda_k}{4}-(50L^2+\frac{8}{\mu_G})\beta_k)\mathbb{E}[\|\Delta g_k\|^2]\notag\\
&\hspace{20pt}+\mathbb{E}[h(\theta_k)]-(\frac{\alpha_k}{8}-\frac{48L^2\alpha_k^2}{\lambda_k}) \mathbb{E}[\|\nabla h(\theta_k)\|^2] \notag\\
&\hspace{20pt}+(1-\mu_G\beta_k)\mathbb{E}[y_k]-\left(\frac{\mu_G}{2}\beta_k -110L^4\alpha_k-L^2\beta_k^2-\frac{80L^4\beta_k^2}{\lambda_k}\right)\mathbb{E}[y_k]+4B\lambda_k^2\notag\\
&\leq (1-\lambda_k)\mathbb{E}[\|\Delta f_k\|^2]+(1-\lambda_k)\mathbb{E}[\|\Delta g_k\|^2]+\mathbb{E}[h(\theta_k)]-\frac{\alpha_k}{16}\mathbb{E}[\|\nabla h(\theta_k)\|^2] \notag\\
&\hspace{20pt}+(1-\mu_G\beta_k)\mathbb{E}[y_k]+4B\lambda_k^2,
\end{align*}
where we use the step size conditions $\beta_k\leq\frac{\lambda_k}{4}\min\{(\frac{37L}{4}+48L^2)^{-1},(50L^2+\frac{8}{\mu_G})^{-1}\}$, $\alpha_k\leq\frac{1}{768L^2}\lambda_k$, $\beta_k\leq \frac{\mu_G}{320L^4}\lambda_k$, $\beta_k\leq\frac{\mu_G}{8L^2}$, and $\alpha_k\leq\frac{\mu_G}{880L^4}\beta_k$.

This inequality obviously implies
\begin{align*}
\alpha_k\mathbb{E}[\|\nabla h(\theta_k)\|^2]&\leq 16\mathbb{E}[\|\Delta f_k\|^2-\|\Delta f_{k+1}\|^2]+16\mathbb{E}[\|\Delta g_k\|^2-\|\Delta g_{k+1}\|^2]\notag\\
&\hspace{20pt}+16\mathbb{E}[h(\theta_k)-h(\theta_{k+1})]+16\mathbb{E}[y_{k}-y_{k+1}]+64B\lambda_k^2.
\end{align*}

Since $\|\Delta f_k\|^2$, $\|\Delta g_k\|^2$, and $y_k$ are non-negative and $h(\theta_k)\geq h(\theta^{\star})$, we have the following inequality by summing up over $k$
\begin{align*}
\sum_{t=0}^{k}\alpha_t\mathbb{E}[\|\nabla h(\theta_t)\|^2]&\leq 16\left(\|\Delta f_0\|^2+\|\Delta g_0\|^2+(h(\theta_0)-h(\theta^{\star}))+y_0\right)+64B\sum_{t=0}^{k}\lambda_t^2
\end{align*}

It is a standard result that (see, for example, Lemma 3 of \citet{zeng2021two})
\begin{align*}
\sum_{t=0}^{k}\alpha_t=\alpha_0\sum_{t=0}^{k}\frac{1}{(t+1)^{1/2}}\geq\frac{\alpha_0(k+1)^{1/2}}{2},\\
\sum_{t=0}^{k}\lambda_t^2=\lambda_0^2\sum_{t=0}^{k}\frac{1}{t+1}\leq2\lambda_0^2\log(k+2)
\end{align*}

Finally,
\begin{align*}
&\min_{t \leq k} \mathbb{E}\left[\|\nabla h(\theta_t)\|^2\right]\\
&\leq \frac{1}{\sum_{t'=0}^k \alpha_{t'}} \sum_{t=0}^k \alpha_t \mathbb{E}\left[\|\nabla h(\theta_t)\|^2\right]\notag\\
&\leq \frac{32}{\alpha_0(k+1)^{1/2}}\left(\|\Delta f_0\|^2+\|\Delta g_0\|^2+(h(\theta_0)-h(\theta^{\star}))+y_0\right)+\frac{4\log(k+2)}{\alpha_0(k+1)^{1/2}}.
\end{align*}

\qed

\subsection{Proof of Lemma~\ref{lem:Lipschitz:nonconvex}}

By the update rule,
\begin{align*}
    \|f_k\|&=\|\Delta f_k+F(\theta_k,\omega_k)-F(\theta_k,\omega^{\star}(\theta_k))+F(\theta_k,\omega^{\star}(\theta_k))\|\\
    &\leq\|\Delta f_k\|+\|F(\theta_k,\omega_k)-F(\theta_k,\omega^{\star}(\theta_k))\|+\|F(\theta_k,\omega^{\star}(\theta_k))\|\\
    &=\|\Delta f_k\|+\|F(\theta_k,\omega_k)-F(\theta_k,\omega^{\star}(\theta_k))\|+\|\nabla h(\theta_k)\|\\
    &\leq \|\Delta f_k\| + L\|\omega_k-\omega^{\star}(\theta_k)\| + \|\nabla h(\theta_k)\|\notag\\
    &\leq \|\Delta f_k\| + L\sqrt{y_k} + \|\nabla h(\theta_k)\|,
\end{align*}
where the second inequality follows from the Lipschitz continuity of $F$.

Similarly, since $G(\theta,\omega^{\star}(\theta))=0$ for any $\theta$, we can derive
\begin{align*}
    \|g_k\|&=\|\Delta g_k+G(\theta_k,\omega_k)-G(\theta_k,\omega^{\star}(\theta_k))\|\notag\\
    &\leq\|\Delta g_k\|+\|G(\theta_k,\omega_k)-G(\theta_k,\omega^{\star}(\theta_k))\|\notag\\
    &\leq \|\Delta g_k\|+L\sqrt{y_k}.
\end{align*}

The bounds on $\|\theta_{k+1}-\theta_k\|$ and $\|\omega_{k+1}-\omega_k\|$ easily follow from the two inequalities above and \eqref{eq:update_auxiliary}.

\qed

\subsection{Proof of Lemma~\ref{lem:Delta_f:nonconvex}}

By the update rule in \eqref{eq:update_auxiliary}, we have
\begin{align*}
\Delta f_{k+1}&=f_{k+1}-F(\theta_{k+1},\omega_{k+1})\notag\\
&=(1-\lambda_k)f_k+\lambda_k F(\theta_k,\omega_k,X_k)-F(\theta_{k+1},\omega_{k+1})\notag\\
&=(1-\lambda_k)f_k+\lambda_k F(\theta_k,\omega_k)-F(\theta_{k+1},\omega_{k+1})+\lambda_k\left(F(\theta_k,\omega_k)-F(\theta_k,\omega_k,X_k)\right)\notag\\
&=(1-\lambda_k)\Delta f_k + F(\theta_k,\omega_k)-F(\theta_{k+1},\omega_{k+1})+\lambda_k\left(F(\theta_k,\omega_k)-F(\theta_k,\omega_k,X_k)\right).
\end{align*}
This leads to 
\begin{align}
&\|\Delta f_{k+1}\|^2\notag\\
&=(1-\lambda_k)^2\|\Delta f_k\|^2 +\Big\|\left(F(\theta_k,\omega_k)-F(\theta_{k+1},\omega_{k+1})\right)+\lambda_k\left(F(\theta_k,\omega_k)-F(\theta_k,\omega_k,X_k)\right)\Big\|^2\notag\\
&\hspace{20pt}+2(1-\lambda_k)\Delta f_k^{\top}\left(F(\theta_k,\omega_k)-F(\theta_{k+1},\omega_{k+1})\right)\notag\\
&\hspace{20pt}+2(1-\lambda_k)\lambda_k\Delta f_k^{\top}\left(F(\theta_k,\omega_k)-F(\theta_k,\omega_k,X_k)\right)\notag\\
&\leq (1-\lambda_k)^2\|\Delta f_k\|^2 +2\|F(\theta_k,\omega_k)-F(\theta_{k+1},\omega_{k+1})\|^2+2\lambda_k^2\|F(\theta_k,\omega_k)-F(\theta_k,\omega_k,X_k)\|^2\notag\\
&\hspace{20pt}+\frac{\lambda_k}{2}\|\Delta f_k\|^2+\frac{2}{\lambda_k}\left\|F(\theta_k,\omega_k)-F(\theta_{k+1},\omega_{k+1})\right\|^2\notag\\
&\hspace{20pt}+2(1-\lambda_k)\lambda_k\Delta f_k^{\top}\left(F(\theta_k,\omega_k)-F(\theta_k,\omega_k,X_k)\right)\notag\\
&\leq (1-\lambda_k)\|\Delta f_k\|^2 +2\lambda_k^2\|F(\theta_k,\omega_k)-F(\theta_k,\omega_k,X_k)\|^2\notag\\
&\hspace{20pt}+(\lambda_k^2-\frac{\lambda_k}{2})\|\Delta f_k\|^2+\frac{4}{\lambda_k}\left\|F(\theta_k,\omega_k)-F(\theta_{k+1},\omega_{k+1})\right\|^2\notag\\
&\hspace{20pt}+2(1-\lambda_k)\lambda_k\Delta f_k^{\top}\left(F(\theta_k,\omega_k)-F(\theta_k,\omega_k,X_k)\right),\label{lem:Delta_f:proof_eq1:nonconvex}
\end{align}
where the last inequality follows from $\lambda_k\leq 1$.

The second term on the right hand side of \eqref{lem:Delta_f:proof_eq1:nonconvex} is bounded in expectation, i.e.
\[\mathbb{E}_{X_k\sim\xi}[\|F(\theta_k,\omega_k,X_k)-F(\theta_k,\omega_k)\|^2]\leq B.\]
% since $\mathbb{E}_{X_k\sim\xi}[\|F(\theta_k,\omega_k,X_k)-F(\theta_k,\omega_k)\|^2]\leq B$.
In addition, the last term of \eqref{lem:Delta_f:proof_eq1:nonconvex} is zero in expectation since $\mathbb{E}[\Delta f_k^{\top}(F(\theta_k,\omega_k)-F(\theta_k,\omega_k,X_k))\mid \Hcal_{k-1}]=\Delta f_k^{\top}\mathbb{E}[(F(\theta_k,\omega_k)-F(\theta_k,\omega_k,X_k))\mid \Hcal_{k-1}]=0$. As a result, we have
\begin{align}
\mathbb{E}[\|\Delta f_{k+1}\|^2]&\leq (1-\lambda_k)\mathbb{E}[\|\Delta f_k\|^2]-(\frac{\lambda_k}{2}-\lambda_k^2)\mathbb{E}[\|\Delta f_k\|^2]\notag\\
&\hspace{20pt}+\frac{4}{\lambda_k}\mathbb{E}[\left\|F(\theta_k,\omega_k)-F(\theta_{k+1},\omega_{k+1})\right\|^2]+2B\lambda_k^2.\label{lem:Delta_f:proof_eq2:nonconvex}
\end{align}

By Lemma~\ref{lem:Lipschitz:nonconvex} and the Lipschitz condition of $F$,
\begin{align*}
    &\|F(\theta_k,\omega_k)-F(\theta_{k+1},\omega_{k+1})\|^2\notag\\
    & \leq 2L^2\|\theta_{k+1}-\theta_k\|^2+ 2L^2\|\omega_{k+1}-\omega_k\|^2\notag\\
    &\leq 6L^2\alpha_k^2(\|\Delta f_k\|^2+L^2 y_k+\|\nabla h(\theta_k)\|^2)+4L^2\beta_k^2(\|\Delta g_k\|^2+L^2 y_k)\notag\\
    &\leq 6L^2\beta_k^2\|\Delta f_k\|^2 + 6L^2\beta_k^2\|\Delta g_k\|^2+10L^4\beta_k^2 y_k + 6L^2 \alpha_k^2 \|\nabla h(\theta_k)\|^2.
\end{align*}
Plugging this bound into \eqref{lem:Delta_f:proof_eq2:nonconvex}, we get
\begin{align*}
&\mathbb{E}[\|\Delta f_{k+1}\|^2]\notag\\
&\leq (1-\lambda_k)\mathbb{E}[\|\Delta f_k\|^2]-(\frac{\lambda_k}{2}-\lambda_k^2)\mathbb{E}[\|\Delta f_k\|^2]\notag\\
&\hspace{20pt}+\frac{4}{\lambda_k}\mathbb{E}\left[6L^2\beta_k^2\|\Delta f_k\|^2 + 6L^2\beta_k^2\|\Delta g_k\|^2+10L^4\beta_k^2 y_k + 6L^2\alpha_k^2 \|\nabla h(\theta_k)\|^2\right] + 2B\lambda_k^2\notag\\
&\leq (1-\lambda_k)\mathbb{E}[\|\Delta f_k\|^2]-(\frac{\lambda_k}{4}-\frac{24L^2\beta_k^2}{\lambda_k})\mathbb{E}[\|\Delta f_k\|^2]\notag\\
&\hspace{20pt}+\frac{4}{\lambda_k}\mathbb{E}\left[6L^2\beta_k^2\|\Delta g_k\|^2+10L^4\beta_k^2 y_k + 6L^2\alpha_k^2 \|\nabla h(\theta_k)\|^2\right] +2B\lambda_k^2,
\end{align*}
where the last inequality follows from $\lambda_k \leq 1/4$.

\qed

\subsection{Proof of Lemma~\ref{lem:Delta_g:nonconvex}}

By the update rule in \eqref{eq:update_auxiliary}, we have
\begin{align*}
\Delta g_{k+1}&=g_{k+1}-G(\theta_{k+1},\omega_{k+1})\notag\\
&=(1-\lambda_k)g_k+\lambda_k G(\theta_k,\omega_k,X_k)-G(\theta_{k+1},\omega_{k+1})\notag\\
&=(1-\lambda_k)g_k+\lambda_k G(\theta_k,\omega_k)-G(\theta_{k+1},\omega_{k+1})+\lambda_k\left(G(\theta_k,\omega_k)-G(\theta_k,\omega_k,X_k)\right)\notag\\
&=(1-\lambda_k)\Delta g_k + G(\theta_k,\omega_k)-G(\theta_{k+1},\omega_{k+1})+\lambda_k\left(G(\theta_k,\omega_k)-G(\theta_k,\omega_k,X_k)\right).
\end{align*}
This leads to 
\begin{align}
&\|\Delta g_{k+1}\|^2\notag\\
&=(1-\lambda_k)^2\|\Delta g_k\|^2 +\Big\|\left(G(\theta_k,\omega_k)-G(\theta_{k+1},\omega_{k+1})\right)+\lambda_k\left(G(\theta_k,\omega_k)-G(\theta_k,\omega_k,X_k)\right)\Big\|^2\notag\\
&\hspace{20pt}+2(1-\lambda_k)\Delta g_k^{\top}\left(G(\theta_k,\omega_k)-G(\theta_{k+1},\omega_{k+1})\right)\notag\\
&\hspace{20pt}+2(1-\lambda_k)\lambda_k\Delta g_k^{\top}\left(G(\theta_k,\omega_k)-G(\theta_k,\omega_k,X_k)\right)\notag\\
&\leq (1-\lambda_k)^2\|\Delta g_k\|^2 +2\|G(\theta_k,\omega_k)-G(\theta_{k+1},\omega_{k+1})\|^2+2\lambda_k^2\|G(\theta_k,\omega_k)-G(\theta_k,\omega_k,X_k)\|^2\notag\\
&\hspace{20pt}+\frac{\lambda_k}{2}\|\Delta g_k\|^2+\frac{2}{\lambda_k}\left\|G(\theta_k,\omega_k)-G(\theta_{k+1},\omega_{k+1})\right\|^2\notag\\
&\hspace{20pt}+2(1-\lambda_k)\lambda_k\Delta g_k^{\top}\left(G(\theta_k,\omega_k)-G(\theta_k,\omega_k,X_k)\right)\notag\\
&\leq (1-\lambda_k)\|\Delta g_k\|^2 +2\lambda_k^2\|G(\theta_k,\omega_k)-G(\theta_k,\omega_k,X_k)\|^2\notag\\
&\hspace{20pt}+(\lambda_k^2-\frac{\lambda_k}{2})\|\Delta g_k\|^2+\frac{4}{\lambda_k}\left\|G(\theta_k,\omega_k)-G(\theta_{k+1},\omega_{k+1})\right\|^2\notag\\
&\hspace{20pt}+2(1-\lambda_k)\lambda_k\Delta g_k^{\top}\left(G(\theta_k,\omega_k)-G(\theta_k,\omega_k,X_k)\right).\label{lem:Delta_g:proof_eq1:nonconvex}
\end{align}

The second term on the right hand side of \eqref{lem:Delta_g:proof_eq1:nonconvex} is bounded in expectation, i.e.
\[\mathbb{E}_{X_k\sim\xi}[\|G(\theta_k,\omega_k,X_k)-G(\theta_k,\omega_k)\|^2]\leq B.\]
% since $\mathbb{E}_{X_k\sim\xi}[\|G(\theta_k,\omega_k,X_k)-G(\theta_k,\omega_k)\|^2]\leq B$.
In addition, the last term of \eqref{lem:Delta_f:proof_eq1:nonconvex} is zero in expectation since $\mathbb{E}[\Delta g_k^{\top}(G(\theta_k,\omega_k)-G(\theta_k,\omega_k,X_k))\mid \Hcal_{k-1}]=\Delta g_k^{\top}\mathbb{E}[(G(\theta_k,\omega_k)-G(\theta_k,\omega_k,X_k))\mid \Hcal_{k-1}]=0$. As a result, we have
\begin{align}
\mathbb{E}[\|\Delta g_{k+1}\|^2]&\leq (1-\lambda_k)\mathbb{E}[\|\Delta g_k\|^2]-(\frac{\lambda_k}{2}-\lambda_k^2)\mathbb{E}[\|\Delta g_k\|^2]\notag\\
&\hspace{20pt}+\frac{4}{\lambda_k}\mathbb{E}[\left\|G(\theta_k,\omega_k)-G(\theta_{k+1},\omega_{k+1})\right\|^2]+2B\lambda_k^2.\label{lem:Delta_g:proof_eq2:nonconvex}
\end{align}

By Lemma~\ref{lem:Lipschitz:nonconvex} and the Lipschitz condition of $G$,
\begin{align*}
    &\|G(\theta_k,\omega_k)-G(\theta_{k+1},\omega_{k+1})\|^2\notag\\
    & \leq 2L^2\|\theta_{k+1}-\theta_k\|^2+ 2L^2\|\omega_{k+1}-\omega_k\|^2\notag\\
    &\leq 6L^2\alpha_k^2(\|\Delta f_k\|^2+L^2 y_k+\|\nabla h(\theta_k)\|^2)+4L^2\beta_k^2(\|\Delta g_k\|^2+L^2 y_k)\notag\\
    &\leq 6L^2\beta_k^2\|\Delta f_k\|^2 + 6L^2\beta_k^2\|\Delta g_k\|^2+10L^4\beta_k^2 y_k + 6L^2 \alpha_k^2 \|\nabla h(\theta_k)\|^2.
\end{align*}
Plugging this bound into \eqref{lem:Delta_g:proof_eq2:nonconvex}, we get
\begin{align*}
&\mathbb{E}[\|\Delta g_{k+1}\|^2]\notag\\
&\leq (1-\lambda_k)\mathbb{E}[\|\Delta g_k\|^2]-(\frac{\lambda_k}{2}-\lambda_k^2)\mathbb{E}[\|\Delta g_k\|^2]\notag\\
&\hspace{20pt}+\frac{4}{\lambda_k}\mathbb{E}\left[6L^2\beta_k^2\|\Delta f_k\|^2 + 6L^2\beta_k^2\|\Delta g_k\|^2+10L^4\beta_k^2 y_k + 6L^2\alpha_k^2 \|\nabla h(\theta_k)\|^2\right]+2B\lambda_k^2\notag\\
&\leq (1-\lambda_k)\mathbb{E}[\|\Delta g_k\|^2]-(\frac{\lambda_k}{4}-\frac{24L^2\beta_k^2}{\lambda_k})\mathbb{E}[\|\Delta g_k\|^2]\notag\\
&\hspace{20pt}+\frac{4}{\lambda_k}\mathbb{E}\left[6L^2\beta_k^2\|\Delta f_k\|^2+10L^4\beta_k^2 y_k + 6L^2\alpha_k^2 \|\nabla h(\theta_k)\|^2\right]+2B\lambda_k^2,
\end{align*}
where the last inequality follows from $\lambda_k \leq 1/4$.

\qed

\subsection{Proof of Lemma~\ref{lem:x_k:nonconvex}}

From the update rule \eqref{alg:update_decision} and the Lipschitz gradient condition in Assumption~\ref{assump:Lipschitz}, we have
\begin{align}
    h(\theta_{k+1})&\leq h(\theta_{k})+\langle\nabla h(\theta_k),\theta_{k+1}-\theta_k\rangle+\frac{L}{2}\|\theta_{k+1}-\theta_k\|^2\notag\\
    &= h(\theta_{k})-\alpha_k\langle\nabla h(\theta_k),f_k\rangle+\frac{L \alpha_k^2}{2}\|f_k\|^2\notag\\
    % &= h(\theta_{k})-\alpha_k\langle\nabla h(\theta_k),\Delta f_k\rangle-\alpha_k\langle\nabla h(\theta_k),F(\theta_k,\omega_k)\rangle+\frac{L \alpha_k^2}{2}\|f_k\|^2\notag\\
    &= h(\theta_{k})-\alpha_k\langle\nabla h(\theta_k),\Delta f_k\rangle-\alpha_k\langle\nabla h(\theta_k),F(\theta_k,\omega^{\star}(\theta_k))\rangle\notag\\
    &\hspace{20pt}+\alpha_k\langle\nabla h(\theta_k),F(\theta_k,\omega^{\star}(\theta_k))-F(\theta_k,\omega_k)\rangle+\frac{L \alpha_k^2}{2}\|f_k\|^2\notag\\
    &=h(\theta_k)-\alpha_k\|\nabla h(\theta_k)\|^2-\alpha_k\langle\nabla h(\theta_k),\Delta f_k\rangle\notag\\
    &\hspace{20pt}+\alpha_k\langle\nabla h(\theta_k),F(\theta_k,\omega^{\star}(\theta_k))-F(\theta_k,\omega_k)\rangle+\frac{L \alpha_k^2}{2}\|f_k\|^2.\label{lem:x_k:proof_eq1:nonconvex}
\end{align}

To bound the second term on the right hand side of \eqref{lem:x_k:proof_eq1:nonconvex}, 
\begin{align}
-\alpha_k\langle\nabla h(\theta_k),\Delta f_k\rangle
&\leq \alpha_k\|\nabla h(\theta_k)\|\|\Delta f_k\|\leq\frac{\alpha_k}{4}\|\nabla h(\theta_k)\|^2+\alpha_k\|\Delta f_k\|^2,
\label{lem:x_k:proof_eq2:nonconvex}
\end{align}
where the final inequality follows from the fact that $2\|\va\|\|\vb\|\leq c\|\va\|^2+\frac{1}{c}\|\vb\|^2$ for any vector $\va,\vb$ and scalar $c>0$.

The third term on the right hand side of \eqref{lem:x_k:proof_eq1:nonconvex} can be treated in a similar manner
\begin{align}
\alpha_k\langle\nabla h(\theta_k),F(\theta_k,\omega^{\star}(\theta_k))-F(\theta_k,\omega_k)\rangle &\leq \alpha_k\|\nabla h(\theta_k)\|\|F(\theta_k,\omega^{\star}(\theta_k))-F(\theta_k,\omega_k)\|\notag\\
&\leq L\alpha_k\|\nabla h(\theta_k)\|\|\omega^{\star}(\theta_k)-\omega_k\|\notag\\
&\leq\frac{\alpha_k}{4}\|\nabla h(\theta_k)\|^2+L^2\alpha_k y_k.
\label{lem:x_k:proof_eq3:nonconvex}
\end{align}

Plugging \eqref{lem:x_k:proof_eq2:nonconvex} and \eqref{lem:x_k:proof_eq3:nonconvex} into \eqref{lem:x_k:proof_eq1:nonconvex}, we get
\begin{align*}
&h(\theta_{k+1})\\
&\leq h(\theta_k)-\alpha_k\|\nabla h(\theta_k)\|^2-\alpha_k\langle\nabla h(\theta_k),\Delta f_k\rangle\notag\\
&\hspace{20pt}+\alpha_k\langle\nabla h(\theta_k),F(\theta_k,\omega^{\star}(\theta_k))-F(\theta_k,\omega_k)\rangle+\frac{L \alpha_k^2}{2}\|f_k\|^2\notag\\
&\leq h(\theta_k)-\alpha_k\|\nabla h(\theta_k)\|^2+\frac{\alpha_k}{4}\|\nabla h(\theta_k)\|^2+\alpha_k\|\Delta f_k\|^2+\frac{\alpha_k}{4}\|\nabla h(\theta_k)\|^2\notag\\
&\hspace{20pt}+L^2\alpha_k y_k+\frac{L \alpha_k^2}{2}\|f_k\|^2\notag\\
&\leq h(\theta_k)-\frac{\alpha_k}{2}\|\nabla h(\theta_k)\|^2+\alpha_k\|\Delta f_k\|^2+L^2\alpha_k y_k+\frac{3L \alpha_k^2}{2}\left(\|\Delta f_k\|^2 + L^2 y_k + \|\nabla h(\theta_k)\|^2\right)\notag\\
&\leq h(\theta_k)-\frac{\alpha_k}{4}\|\nabla h(\theta_k)\|^2+\frac{5\alpha_k}{4}\|\Delta f_k\|^2+\frac{5L^2\alpha_k}{4} y_k,
\end{align*}
where the third inequality plugs in the result of Lemma~\ref{lem:Lipschitz:nonconvex}, and the final inequality uses $\alpha_k\leq\frac{1}{6L}$.

% This implies
% \begin{align*}
% \alpha_k\|\nabla h(\theta_k)\|^2&\leq 4\left(h(\theta_k)-h(\theta_{k+1})\right)+5\alpha_k\alpha_k\|\Delta f_k\|^2+5L^2\alpha_k y_k.
% \end{align*}

\qed

\subsection{Proof of Lemma~\ref{lem:y_k:nonconvex}}

By the update rule in \eqref{alg:update_decision},
\begin{align*}
\omega_{k+1}-\omega^{\star}(\theta_{k+1})&=\omega_{k+1}-\omega^{\star}(\theta_{k+1})\notag\\
&=\left(\omega_{k}-\omega^{\star}(\theta_{k})\right)-\beta_k g_k - \left(\omega^{\star}(\theta_{k+1})-\omega^{\star}(\theta_k)\right)\notag\\
&=\left(\omega_{k}-\omega^{\star}(\theta_{k})\right)-\beta_k G(\theta_k,\omega_k)-\beta_k \Delta g_k - \left(\omega^{\star}(\theta_{k+1})-\omega^{\star}(\theta_k)\right).
\end{align*}
Taking the norm yields
\begin{align}
y_{k+1}&=\|\omega_{k+1}-\omega^{\star}(\theta_{k+1})\|^2\notag\\
&= \underbrace{\|\omega_{k}-\omega^{\star}(\theta_{k})-\beta_k G(\theta_k,\omega_k)\|^2}_{T_1}+\beta_k^2\|\Delta g_k\|^2+\underbrace{\|\omega^{\star}(\theta_{k+1})-\omega^{\star}(\theta_k)\|^2}_{T_2}\notag\\
&\hspace{20pt} \underbrace{-2\left(\omega_{k}-\omega^{\star}(\theta_{k})-\beta_k G(\theta_k,\omega_k)\right)^{\top}\left(\omega^{\star}(\theta_{k+1})-\omega^{\star}(\theta_k)\right)}_{T_3}\notag\\
&\hspace{20pt}\underbrace{-2\beta_k\left(\omega_{k}-\omega^{\star}(\theta_{k})-\beta_k G(\theta_k,\omega_k)\right)^{\top}\Delta g_k}_{T_4} + \underbrace{2\beta_k \Delta g_k^{\top}\left(\omega^{\star}(\theta_{k+1})-\omega^{\star}(\theta_k)\right)}_{T_5}.\label{lem:y_k:proof_eq1:nonconvex}
\end{align}

We bound each term of \eqref{lem:y_k:proof_eq1:nonconvex} individually. First, since by definition $G(\theta,\omega^{\star}(\theta))=0$ for all $\theta$, we have
\begin{align*}
T_1&=\|\omega_{k}-\omega^{\star}(\theta_{k})\|^2 - 2\beta_k\left(\omega_{k}-\omega^{\star}(\theta_{k})\right)^{\top}G(\theta_k,\omega_k)+\beta_k^2\|G(\theta_k,\omega_k)\|^2\notag\\
&=y_k - 2\beta_k\left(\omega_{k}-\omega^{\star}(\theta_{k})\right)^{\top}G(\theta_k,\omega_k)+\beta_k^2\|G(\theta_k,\omega_k)-G(\theta_k,\omega^{\star}(\theta_k))\|^2\notag\\
&\leq y_k - 2\mu_G\beta_k\|\omega_k-\omega^{\star}(\theta_k)\|^2+L^2\beta_k^2\|\omega_k-\omega^{\star}(\theta_k)\|^2\notag\\
&=(1-2\mu_G\beta_k+L^2\beta_k^2)y_k,
\end{align*}
where the first inequality follows from Assumption~\ref{assump:stronglymonotone_G} and the Lipschitz continuity of $G$.

The term $T_2$ can be simply treated using the Lipschitz condition of $\omega^{\star}$
\begin{align*}
T_2 
&=\|\omega^{\star}(\theta_{k+1})-\omega^{\star}(\theta_k)\|^2
\leq L^2\|\theta_{k+1}-\theta_k\|^2\leq L^2\alpha_k^2\left(\|\Delta f_k\| + L\sqrt{y_k} + \|\nabla h(\theta_k)\|\right)^2\notag\\
&\leq 3L^2\alpha_k^2\left(\|\Delta f_k\|^2 + L^2 y_k + \|\nabla h(\theta_k)\|^2\right).
\end{align*}

By the Cauchy-Schwarz inequality, 
\begin{align*}
&T_3\notag\\
&\leq 2\|\omega_{k}-\omega^{\star}(\theta_{k})-\beta_k G(\theta_k,\omega_k)\|\|\omega^{\star}(\theta_{k+1})-\omega^{\star}(\theta_k)\|\notag\\
&\leq 2\|\omega_{k}\hspace{-2pt}-\hspace{-2pt}\omega^{\star}(\theta_{k})\|\|\omega^{\star}(\theta_{k+1})-\omega^{\star}(\theta_k)\|+2\beta_k\|G(\theta_k,\omega_k)-G(\theta_k, \omega^{\star}(\theta_k))\|\|\omega^{\star}(\theta_{k+1})-\omega^{\star}(\theta_k)\|\notag\\
&\leq (2L+2L^2\beta_k)\sqrt{y_k}\|\theta_{k+1}-\theta_k\|\notag\\
&\leq 4L\alpha_k\sqrt{y_k}\left(\|\Delta f_k\| + L\sqrt{y_k} + \|\nabla h(\theta_k)\|\right).
\end{align*}
where the second inequality follows from $G(\theta,\omega^{\star}(\theta))=0$ for all $\theta$, and the fourth inequality applies Lemma~\ref{lem:Lipschitz:nonconvex} and the step size rule $\beta_k\leq1/L$.

We can bound $T_4$ following a similar line of analysis
\begin{align*}
T_4 &\leq 2\beta_k\|\omega_{k}-\omega^{\star}(\theta_{k})-\beta_k G(\theta_k,\omega_k)\|\|\Delta g_k\|\leq 4\beta_k\sqrt{y_k}\|\Delta g_k\|.
\end{align*}

Finally, we treat $T_5$ again by the Cauchy-Schwarz inequality
\begin{align*}
T_5&\leq2\beta_k\|\Delta g_k\|\|\omega^{\star}(\theta_{k+1})-\omega^{\star}(\theta_k)\|\notag\\
&\leq 2L\alpha_k\beta_k\|\Delta g_k\| \Big(\|\Delta f_k\| + L\sqrt{y_k} + \|\nabla h(\theta_k)\|\Big),
\end{align*}
where the second inequality again uses Lemma~\ref{lem:Lipschitz:nonconvex}.

Applying the bounds on $T_1$-$T_5$ to \eqref{lem:y_k:proof_eq1:nonconvex}, we get
\begin{align*}
&y_{k+1}\notag\\
& \leq (1-2\mu_G\beta_k+L^2\beta_k^2)y_k+\beta_k^2\|\Delta g_k\|^2+3L^2\alpha_k^2\left(\|\Delta f_k\|^2 + L^2 y_k + \|\nabla h(\theta_k)\|^2\right)\notag\\
&\hspace{20pt}+4L\alpha_k\sqrt{y_k}\left(\|\Delta f_k\| + L\sqrt{y_k} + \|\nabla h(\theta_k)\|\right) + 4\beta_k\sqrt{y_k}\|\Delta g_k\|\notag\\
&\hspace{20pt}+2L\alpha_k\beta_k\|\Delta g_k\| \Big(\|\Delta f_k\| + L\sqrt{y_k} + \|\nabla h(\theta_k)\|\Big)\\
&\leq (1-2\mu_G\beta_k+L^2\beta_k^2)y_k+\beta_k^2\|\Delta g_k\|^2+3L^2\alpha_k^2\left(\|\Delta f_k\|^2 + L^2 y_k + \|\nabla h(\theta_k)\|^2\right)\notag\\
&\hspace{20pt}+\hspace{-2pt}2L\alpha_k y_k+2L\alpha_k\|\Delta f_k\|^2\hspace{-2pt}+\hspace{-2pt}4L^2\alpha_k y_k+96L^2 \alpha_k y_k\hspace{-2pt}+\hspace{-2pt}\frac{\alpha_k}{24}\|\nabla h(\theta_k)\|^2\hspace{-2pt}+\hspace{-2pt} \frac{\mu_G\beta_k}{2} y_k\hspace{-2pt}+\hspace{-2pt}\frac{8\beta_k}{\mu_G}\|\Delta g_k\|^2\notag\\
&\hspace{20pt}+L\alpha_k\beta_k\|\Delta g_k\|^2+3L\alpha_k\beta_k\left(\|\Delta f_k\|^2 + L^2 y_k + \|\nabla h(\theta_k)\|^2\right)\\
&\leq (1-\mu_G\beta_k)y_k-\left(\frac{\mu_G}{2}\beta_k -108L^4\alpha_k-L^2\beta_k^2\right)y_k\notag\\
&\hspace{20pt}+\frac{\alpha_k}{8}\|\nabla h(\theta_k)\|^2+\left(\frac{8\beta_k}{\mu_G}+2L\beta_k^2\right)\|\Delta g_k\|^2+8L\alpha_k\|\Delta f_k\|^2,
\end{align*}
where the second inequality follows from the fact that $2\|\va\|\|\vb\|\leq c\|\va\|^2+\frac{1}{c}\|\vb\|^2$ for any vector $\va,\vb$ and scalar $c>0$. We have also simplified terms using the conditions $\alpha_k\leq\frac{1}{72L^2}$ and $\beta_k\leq\frac{1}{72L}$.

\qed
\section{Simulations Details}\label{sec:appendix_simulation}

We generate completely random MDPs with transition probability matrix $\Pcal\in\mathbb{R}^{|\Scal|\times|\Acal|\times |\Scal|}$ drawn i.i.d. from $\text{Unif}(0,1)$ and then normalized such that
\begin{align*}
    \sum_{s'\in\Scal}\Pcal(s'\mid s,a)=1,\quad\forall s,a.
\end{align*}
The behavior policy $\pi_b$ is chosen to be uniform for all states. The policy $\pi$ to be evaluated is generated randomly by first drawing $\psi\in\mathbb{R}^{|\Scal|\times|\Acal|}$ such that $\psi_{s,a}\sim N(0,1)$, which produces the policy $\pi$ through the softmax function
\begin{align*}
    \pi(a\mid s)=\frac{\exp(\psi_{s,a})}{\sum_{a'}\exp(\psi_{s,a'})},\quad\forall s,a.
\end{align*}

We draw the feature matrix $\Phi\in\mathbb{R}^{|\Scal|\times dx}$ entry-wise i.i.d. from $N(0,1)$. The value function $V\in\mathbb{R}^{|\Scal|}$ is the solution to the Bellman equation 
\begin{align*}
    V = R + \gamma P^{\pi} V,
\end{align*}
where the matrix $P^{\pi}\in\mathbb{R}^{|\Scal|\times|\Scal|}$ is the state transition matrix under policy $\pi$ 
\[P_{s,s'}=\sum_{a}\Pcal(s'\mid s,a)\pi(a\mid s)\]

We would like to design the experiments such that the value function $V\in\mathbb{R}^{|\Scal|}$ lies in the span of the feature matrix. To achieve, we determine the optimal parameter $\theta^{\star}$ first, by sampling it entry-wise i.i.d. from $N(0,1)$. This produces the value function $V=\Phi\theta$, from which we reverse engineer the reward function
\[R=(I-\gamma P^{\pi})V.\]
The discount factor $\gamma$ is 0.5.

All algorithms start with zero initialization, i.e. $\theta_0=\omega_0=f_0=g_0=0$. The step sizes $\alpha_k$ and $\beta_k$ are set to constant values $5e-4$ and $2e-3$ for simplicity, while the step size $\lambda_k$ is selected as
\[\lambda_k=\frac{4}{5(k+10)}.\]

\end{document}